\documentclass[11pt]{article}
\usepackage[T1]{fontenc} 
\usepackage{lmodern} 
\usepackage{amsmath,amssymb,amsfonts,amsthm,mathtools}
\usepackage{mathrsfs,bm,cases,empheq}

\usepackage{indentfirst} % 首行缩进
\usepackage{xcolor}      % 颜色支持
\usepackage{graphicx,graphics}
\usepackage{etoolbox,xstring,mfirstuc,textcase}

\usepackage{cite} % 自动排序参考文献引用

\usepackage[colorlinks=true,linkcolor=blue,citecolor=blue,urlcolor=blue]{hyperref}

\usepackage[capitalize,nameinlink]{cleveref}

\usepackage{changes}
\usepackage{lineno}
\definechangesauthor[name={R1 cusse}, color=orange]{r1}
\definechangesauthor[name={Ty cusse}, color=purple]{ty}

\allowdisplaybreaks[4] % 允许公式跨页显示

\newtheorem{theorem}{\textbf{Theorem}}[section]
\newtheorem{lemma}{\textbf{Lemma}}[section]
\newtheorem{proposition}{\textbf{Proposition}}[section]
\newtheorem{corollary}{\textbf{Corollary}}[section]
\newtheorem{remark}{\textbf{Remark}}[section]
\newtheorem{definition}{\textbf{Definition}}[section]
\newtheorem{assumption}{\textbf{Assumption}}[section]

\numberwithin{equation}{section}

\def\be{\begin{equation}}
\def\ee{\end{equation}}
\def\bea{\begin{eqnarray}}
\def\eea{\end{eqnarray}}
\def\bt{\begin{theorem}}
\def\et{\end{theorem}}
\def\bl{\begin{lemma}}
\def\el{\end{lemma}}
\def\bnum{\begin{numcases}{}}
\def\enum{\end{numcases}}
\def\br{\begin{remark}}
\def\er{\end{remark}}
\def\bp{\begin{proposition}}
\def\ep{\end{proposition}}
\def\bc{\begin{corollary}}
\def\ec{\end{corollary}}
\def\bd{\begin{definition}}
\def\ed{\end{definition}}
\def\ba{\begin{assumption}}
\def\ea{\end{assumption}}

\begin{document}
\title{Strong solutions to the Cahn--Hilliard--Darcy--Forchheimer system for incompressible two-phase flows with a soluble surfactant}

\author{
Maurizio Grasselli
\thanks{Dipartimento di Matematica,
Politecnico di Milano, Via E. Bonardi 9, 20133 Milano, Italy, Email: \texttt{maurizio.grasselli@polimi.it}.}
\and
Bohan Ouyang
\thanks{Corresponding author. School of Mathematical Sciences,
Fudan University, Handan Road 220, 200433 Shanghai, P.R. China, Email: \texttt{22110180033@m.fudan.edu.cn}.}
\and
Hao Wu
\thanks{School of Mathematical Sciences
and Shanghai Key Laboratory for Contemporary Applied Mathematics, Fudan University, Handan Road 220, 200433 Shanghai, P.R. China, Email: \texttt{haowufd@fudan.edu.cn}.}
}
%%%%%%%%%%%%%%%%%%%%%%%%%%%%%%%%

\date{\today}

\maketitle
%%%%%%%%%%%%%%%%%%%%%%%%%%%%%%%%
\begin{abstract}
\noindent 
We study a diffuse-interface model for incompressible viscous two-phase flows with a soluble surfactant in a bounded porous medium in two or three dimensions. 
The resulting system consists of a Darcy--Forchheimer equation for the seepage velocity $\bm{u}$ coupled with two Cahn--Hilliard equations, one for the difference in volume fractions of the two fluids denoted by $\phi$, and the other for the local concentration of the surfactant $\psi$. 
In the general scenario of non-constant coefficients and mobilities, we establish the existence of a local (resp. global) strong solution in three (resp. two) dimensions to the corresponding initial-boundary value problem with the no-flux boundary conditions for the phase-field variables and chemical potentials and a no-penetration condition for the velocity.
In the case of constant coefficients in the Darcy--Forchheimer equation and constant mobilities in the Cahn--Hilliard equations, with a suitably large Forchheimer exponent, we establish the existence of a global strong solution in three dimensions. 
Furthermore, using the {\L}ojasiewicz--Simon approach, we prove the global existence of a unique strong solution when the initial data are sufficiently close to a local energy minimizer of the free energy.


\end{abstract}

%%%%%%%%%%%%%%%%%%%%%%%%%%%%%%%%

\section{Introduction} 
\setcounter{equation}{0}

In \cite{GOPW2026}, the authors introduced and analyzed a generalization of the Cahn--Hilliard--Hele--Shaw system proposed in \cite{LLG2002} to describe the dynamics of incompressible two-phase flows in a Hele--Shaw cell, which itself is a simplification of the classical ``Model H'' that couples the Navier--Stokes equations with a convective Cahn--Hilliard equation \cite{AMW1998,HH1977}. Compared with the Cahn--Hilliard--Hele--Shaw framework, the model in \cite{GOPW2026} incorporates a soluble surfactant through the additional phase-field variable $\psi$, allows mass transfer through the source term in the $\phi$-equation, and includes a Darcy--Forchheimer law that describes the fluid motion in a porous medium. 
This law contains an additional nonlinear drag term known as the Forchheimer term \cite{Bear70,W1996}. 
The existence of a global weak solution and its convergence to a single equilibrium are the known theoretical results so far (see \cite{GOPW2026}; see also \cite{BF2025} for a related result). In this study, we deal with the existence and uniqueness of strong solutions. This is a delicate issue, especially in three dimensions, even for the standard Cahn--Hilliard--Hele--Shaw system (see Section \ref{lit} below). As we shall see, the presence of the Forchheimer term makes the analysis more complicated when dealing with strong solutions, but it helps in proving the existence of a global strong solution in dimension three also for initial data that are not necessarily close to an energy minimizer, at least in the case of constant coefficients and mobilities. 
This is still an open question for the original Cahn--Hilliard--Hele--Shaw system. In a similar spirit, it is worth recalling the so-called Navier--Stokes--Brinkman--Forchheimer equation, a ``tamed'' version of the Navier--Stokes system. In this case, thanks to the additional regularity given by the Forchheimer term, it is possible to establish a number of results in dimension three which are unknown for the usual Navier--Stokes system (see, e.g.,  \cite{SZ2025} and references therein, see also \cite{BDM2026} for the Navier--Stokes--Cahn--Hilliard system).

\subsection{Model description}
Let $\Omega \subset \mathbb{R}^d$ ($d \in \{ 2,3 \}$) be a bounded domain with a smooth boundary $\partial\Omega$, and let $T>0$ be a given final time. We consider the following diffuse-interface model for an incompressible two-phase flow with a soluble surfactant (see \cite{GOPW2026}):
\begin{numcases}{}
    \nu (\phi,\psi) \bm{u} + \eta (\phi,\psi) |\bm{u}|^{r-2} \bm{u} + \nabla \pi = \mu_{\phi} \nabla \phi + \mu_{\psi} \nabla \psi, \label{eq:chdf1}
    \\
    \mathrm{div} \, \bm{u} =0, \label{eq:chdf2}
    \\
    \partial_t \phi + \bm{u} \cdot \nabla \phi + \sigma_1 \big( \overline{\phi} -c \big) = \mathrm{div} (m_{\phi} (\phi) \nabla \mu_\phi), \label{eq:chdf3}
    \\
    \mu_{\phi} = -\Delta \phi + \sigma_2 \mathcal{N} \big( \phi - \overline{\phi} \big) + F_{\phi}'(\phi) + \partial_{\phi} G(\phi,\psi), \label{eq:chdf4}
    \\
    \partial_t \psi + \bm{u} \cdot \nabla \psi = \mathrm{div} (m_{\psi} (\psi) \nabla \mu_\psi), \label{eq:chdf5}
    \\
    \mu_{\psi} = -\beta \Delta \psi + F_{\psi}'(\psi) + \partial_\psi G(\phi,\psi), \label{eq:chdf6}
\end{numcases}
in $\Omega\times(0,T)$, 
subject to the boundary and initial conditions
\begin{numcases}{}
    \bm{u} \cdot \mathbf{n} = \partial_{\mathbf{n}} \phi = \partial_{\mathbf{n}} \mu_{\phi} = \partial_{\mathbf{n}} \psi = \partial_{\mathbf{n}} \mu_{\psi} = 0, \label{eq:chdfb}& on $\partial \Omega \times (0,T)$, 
    \\
    \phi|_{t=0} = \phi_0(x), \quad \psi|_{t=0} = \psi_0(x), \label{eq:chdfi} & in $\Omega$.
\end{numcases}
Here, $\mathbf n$ denotes the outward unit normal vector to $\partial\Omega$, and $\partial_{\mathbf n}$ is the corresponding outward normal derivative.

The state variables consist of the seepage velocity $\bm{u}: \Omega \times [0,T) \to \mathbb{R}^d$, the pressure $\pi:\Omega\times[0,T)\to\mathbb R$, the phase-field variable $\phi:\Omega\times[0,T)\to[-1,1]$, and the surfactant concentration  $\psi:\Omega\times[0,T)\to[0,1]$. The order parameter $\phi$ represents the volume-fraction difference between the two fluid components, whereas $\psi$ denotes the volume fraction of the soluble surfactant. In the Darcy--Forchheimer equation \eqref{eq:chdf1}, $\nu$ and $\eta$ denote the kinematic viscosity and the Forchheimer coefficient, respectively. Both of them are allowed to depend on the concentrations $\phi, \psi$ in this study. The nonlinear term $\eta|\bm{u}|^{r-2}\bm{u}$ is referred to as the Forchheimer drag term with $r \in [2,\infty)$ being a prescribed exponent (the case $r=2$ reduces to the Darcy law). 
The Forchheimer law (typically $r=3$) provides a nonlinear correction to the classical Darcy law for transport processes in porous media with relatively high velocities, where inertial drag effects cannot be neglected \cite{AGO2007,NB1992,SR2014,W1996}. 

The free energy associated with the Cahn--Hilliard system \eqref{eq:chdf3}--\eqref{eq:chdf6} for the phase-field variables $\phi$, $\psi$ is 
given by
\[
\begin{aligned}
    E_{\mathrm{free}}(\phi,\psi)
    &\stackrel{\rm def}{=} \int_\Omega \left(\frac12|\nabla\phi|^2 + F_\phi(\phi)
        +\frac{\sigma_2}{2}
        \left|\nabla\mathcal N(\phi-\overline\phi)\right|^2
        \right) \, \mathrm{d} x 
      \\
    &\quad +\int_\Omega
    \left(  \frac{\beta}{2}|\nabla\psi|^2
         +F_\psi(\psi) \right) \, \mathrm{d} x
    +\int_\Omega G(\phi,\psi) \, \mathrm{d} x,
\end{aligned}
\]
where $\beta>0$ is a regularizing parameter, cf. \cite{La1992,KK1997,EDAT2013}. The chemical potentials $\mu_\phi,\,\mu_\psi: \Omega \times [0,T)\to \mathbb{R}$ in \eqref{eq:chdf4} and \eqref{eq:chdf6} are defined by the variational derivatives of $E_{\mathrm{free}}$ with respect to $\phi$ and $\psi$, respectively. 
The mobility coefficients $m_\phi$ and $m_\psi$ characterize the strength of diffusion. They are assumed to be non-degenerate and concentration-dependent. In $E_{\mathrm{free}}$, the nonlinear function $F_\phi$ (resp. $F_\psi$) denotes the bulk potential associated with the binary fluid mixture (resp. the surfactant). Throughout this paper, both $F_\phi$ and $F_\psi$ are assumed to be singular, with the Flory--Huggins potential serving as a prototypical example (see \cite[Introduction]{OGW2025} and references therein). The singular nature of these bulk potentials guarantees that the phase-field variables remain within their corresponding physical domains along the evolution, that is, $\phi\in [-1,1]$ and $\psi\in[0,1]$, cf. \cite{AW2007,B1999}. 
When $\sigma_2\neq 0$, the free energy associated with $\phi$ is of Ohta--Kawasaki type, which also includes possible nonlocal interactions between the two fluid components (see \cite{NO1995}).
The operator $\mathcal{N}$ denotes the inverse of the minus Laplacian with a homogeneous Neumann boundary condition and $\overline{f}$ denotes the spatial average of an integrable function $f$ over $\Omega$, i.e., $\overline{f} = | \Omega |^{-1} \int_{\Omega} f \, \mathrm{d}x$, $\vert \Omega\vert$ being the $d$-dimensional Lebesgue measure of $\Omega$. 
The source term $\sigma_1(\overline{\phi} -c)$, where $\sigma_1\geq 0$ and $c\in (-1,1)$, accounts for possible chemical reactions (see \cite{GGM2017} and references therein). Finally, the bivariate function $G$ represents the density of the interaction energy between the binary fluids and the surfactant. A typical choice is
\begin{equation*}
G(\phi,\psi)
=
\frac{\gamma_1}{2}\psi\phi^2
-
\frac{\gamma_2}{4}\psi(\phi^2-1)^2,
\end{equation*}
where $\gamma_1,\gamma_2\in\mathbb R_+$ are given constants. 
The first term penalizes the presence of free surfactant in the bulk phases, 
while the second one favors the accumulation of surfactant at the free interfaces between the two fluids. 
In other words, the first term globally counteracts the occurrence of free surfactants, while the second one locally attracts surfactants to the interface; see \cite{EDAT2013} for further discussions on $G$.

The complex dynamics of the multi-component system is driven by the total energy, i.e., the free energy $E_{\mathrm{free}}$.
For sufficiently smooth solutions to the hydrodynamic system \eqref{eq:chdf1}--\eqref{eq:chdf6} subject to boundary conditions \eqref{eq:chdfb}, a formal calculation yields the following energy identity:
\begin{align}
&\frac{\mathrm d}{\mathrm dt}
E_{\mathrm{free}}(\phi,\psi)
+\int_\Omega
\left(
\nu(\phi,\psi)|\bm{u}|^2
+\eta(\phi,\psi)|\bm{u}|^r
\right)\, \mathrm{d} x 
\nonumber \\
&\quad
+\int_\Omega
\left(
m_\phi(\phi)|\nabla\mu_\phi|^2
+m_\psi(\psi)|\nabla\mu_\psi|^2
\right)\, \mathrm{d} x
+\sigma_1(\overline\phi-c)
\int_\Omega\mu_\phi\, \mathrm{d} x
=0,
\label{basic_energy_identity}
\end{align}
in $(0,T)$.
In addition, integrating equations \eqref{eq:chdf3}, \eqref{eq:chdf5} over $\Omega$ gives the mass balance relations
\begin{align}
\frac{\mathrm d}{\mathrm dt}\overline\phi
+\sigma_1(\overline\phi-c)=0,
\qquad
\frac{\mathrm d}{\mathrm dt}\overline\psi=0,
\qquad \text{in } (0,T).
\label{mass-relation}
\end{align}
The above energy identity and mass balance relations are the cornerstones of the mathematical analysis of the initial-boundary value problem \eqref{eq:chdf1}--\eqref{eq:chdfi}.   

%%%%%%%%%%%%%%%%%%%%%%%%%%%%%%%%
\subsection{Related literature}
\label{lit}

For incompressible two-phase flows with a soluble surfactant, related diffuse-interface models with the full Navier--Stokes system in place of the Darcy--Forchheimer equation \eqref{eq:chdf1} have been studied, for instance, in \cite{DGW2023,GLS2014,OGW2025,GOW2026}
(see also the references therein). 
Instead, the mathematical analysis of the Cahn--Hilliard--Darcy--Forchheimer system \eqref{eq:chdf1}--\eqref{eq:chdf6} remains largely unexplored. As mentioned in the introduction, in \cite{GOPW2026}, we established the existence of global weak solutions to problem \eqref{eq:chdf1}--\eqref{eq:chdfi} and demonstrated that every global weak solution converges to a single equilibrium as $t\to \infty$, by an extended {\L}ojasiewicz--Simon approach developed in \cite{GrasPoi2025}. 
We also point out that, in \cite{GOPW2026}, we accounted for the presence of an inertial term $\alpha\partial_t\bm{u}$, $\alpha>0$, and we supposed that $\sigma_1$ may depend on $\phi$, but we did not consider the case $\sigma_2<0$.

Here, our aim is to develop a strong solution theory for this problem.
If the Forchheimer drag term and all surfactant-related equations and terms are removed, system \eqref{eq:chdf1}--\eqref{eq:chdf6} reduces to the Cahn--Hilliard--Hele--Shaw system (or the Cahn--Hilliard--Darcy system), 
for which a more developed theory on strong solutions is available. 
For example, in the case of constant coefficients and singular potentials including the Flory--Huggins ones, the existence and uniqueness of global-in-time strong solutions for initial data sufficiently close to a local energy minimizer of the free energy were established in \cite{GGW2018}. 
Subsequently, for the generalized system with concentration-dependent kinematic viscosity, constant mobility, and singular potentials, in \cite{G2020} the author proved the existence and uniqueness of global-in-time strong solutions in two dimensions, 
while in three dimensions he showed the existence and uniqueness of strong solutions,
which are local in time for large data or global in time for appropriate small data. 
These works provide the primary mathematical background for our current analysis.

Another related extension of the Cahn--Hilliard--Darcy framework can be obtained by incorporating viscous friction through a Brinkman-type correction, leading to the so-called Cahn--Hilliard--Brinkman system \cite{VVM2001}. In \cite{DG2016}, the authors studied a nonlocal version of this system and established the existence and uniqueness of weak solutions, as well as the inviscid limit toward the nonlocal Cahn--Hilliard--Darcy system. 
Concerning the nonlocal version, it is worth noting that, assuming degenerate mobility, the existence of a global strong solution can be proven for unmatched viscosities without the assumption of an initial datum close to an energy minimizer (see \cite{CFG2022}). We also refer to \cite{CG2020} for the well-posedness analysis of a Cahn--Hilliard--Brinkman system with logarithmic free energy density and unmatched viscosities.

Despite these developments, a strong well-posedness theory for the present Cahn--Hilliard--Darcy--Forchheimer system remains unavailable. 
In particular, the simultaneous presence of the Forchheimer nonlinearity, singular potentials, and concentration-dependent coefficients leads to additional difficulties that do not occur in the corresponding Cahn--Hilliard--Hele--Shaw or Cahn--Hilliard--Brinkman systems. 
The strong solution theory developed in this work aims to overcome these obstacles.

\subsection{Summary of results and strategy of the proof}

We establish the existence of strong solutions to the initial-boundary value problem \eqref{eq:chdf1}--\eqref{eq:chdfi} and prove uniqueness under the strict separation and regularity assumptions specified below.
In the general scenario of non-constant coefficients and mobilities, we prove the existence of local-in-time strong solutions in three dimensions and global-in-time strong solutions in two dimensions. A distinguishing feature of our analysis is the crucial regularizing role played by the Forchheimer term. In particular, when the coefficients and mobilities are positive constants, its superlinear dissipation provides sufficient integrability for the velocity field to close the higher-order energy estimates, yielding the existence of global-in-time strong solutions in three dimensions for $r\in[3,\infty)$. 
To the best of our knowledge, such a global result remains unavailable for the classical Cahn--Hilliard--Darcy system. 
Furthermore, using the {\L}ojasiewicz--Simon approach, we establish the global-in-time existence of strong solutions for initial data sufficiently close to an arbitrary local energy minimizer. 
Our main results are precisely stated in Theorems \ref{ss_chdf_exist_3d}, \ref{ss_chdf_exist_global_small}, and \ref{ss_chdf_exist_2d}.

A primary obstacle in proving the existence of strong solutions is that the regularity of global weak solutions obtained in \cite[Theorem 2.1]{GOPW2026} is insufficient to rigorously justify the required higher-order \textit{a priori} estimates. 
We overcome this difficulty by constructing suitable smooth approximate solutions and deriving uniform higher-order estimates. 
To circumvent limitations of the time-difference approach in the presence of the mass source term $\sigma_1(\overline{\phi}-c)$ and singular potentials, we exploit a suitable compactness argument. 
Specifically, we replace the singular potential functions with their quadratic approximations and apply a Galerkin scheme to the entire system, which preserves the monotonicity structure of the momentum equation \eqref{eq:chdf1}.

The convective terms in the Cahn--Hilliard equations further complicate the derivation of higher-order estimates. 
Differentiating the higher-order energy functional naturally generates cross-terms involving $\partial_t\mu_\phi$, $\bm{u}\cdot\nabla\phi$,  $\partial_t\mu_\psi$, and $\bm{u}\cdot\nabla\psi$. 
Using the momentum equation, we rewrite these terms in terms of $\bm{u}$, $\partial_t\phi$, $\partial_t\psi$, and the gradients of chemical potentials. 
In three dimensions, bounding these terms requires enhanced spatial integrability of the velocity field. 
For the classical Cahn--Hilliard--Darcy system, the basic energy law provides only $L^2$-control of the velocity, which fails to close the estimates without introducing problematic higher-order terms. 
In contrast, the Forchheimer dissipation intrinsically supplies an additional $L^r$-control of the velocity. 
For $r\ge3$, this enhanced integrability is precisely what is needed to reduce the power of the estimate of the corresponding terms and close the higher-order energy estimate, thereby ensuring global-in-time strong solutions for the constant-coefficient scenario.

When the coefficients depend on the concentrations, additional difficulties arise from the nonlinear terms in both the momentum equation and the Cahn--Hilliard equations. Controlling these terms requires an $H^1$-estimate for the velocity field. Formally, this estimate can be obtained by applying the curl operator to the momentum equation and testing the resulting equation by $\mathrm{curl} \,\bm{u}$, together with a div--curl estimate. However, the Darcy--Forchheimer equation contains no viscous regularization, and the basic energy estimate does not provide the spatial regularity needed to justify this calculation.
To make this argument rigorous, we construct a Galerkin approximation whose finite-dimensional spaces are invariant under the double-curl operator. This property ensures that $\mathrm{curl} \,\mathrm{curl} \,\bm{u}_m$ is an admissible test function in the approximate momentum equation. The choice of the Galerkin basis is essential here, since the double curl of a smooth divergence-free vector field with vanishing normal component need not belong to $\bm{L}^2_\sigma(\Omega)$ (see Section \ref{pre} for the notation). We therefore introduce the auxiliary Stokes operator $\bm A=\mathrm{curl} \,\mathrm{curl} \,$ on $\bm{L}^2_\sigma(\Omega)$ with domain
\[
\mathcal D(\bm A)
=
\left\{
\bm v\in H^2(\Omega)^d\cap\bm{L}^2_\sigma(\Omega):
(\mathrm{curl} \,\bm v)\times\bm n=\bm 0
\ \text{a.e. on }\partial\Omega
\right\}.
\]
We recall that the boundary conditions 
\[
    (\mathrm{curl} \, \bm v) \times \bm n = \bm 0, \quad \bm v \cdot \bm n = 0
\]
are called vorticity-slip boundary conditions.
By \cite[Theorem 3.5]{XX2013} (taking the slip parameter therein to be zero), the operator $I+\bm A$ has a compact inverse, and its eigenfunctions form a complete orthonormal basis of $\bm{L}^2_\sigma(\Omega)$. The shift by the identity accommodates the possible nontrivial harmonic kernel of $\bm A$, so no topological restriction on $\Omega$ is required.
The additional vorticity boundary condition is used to construct the approximation spaces and is not imposed on the solution of the original Darcy--Forchheimer problem.
A further difficulty arises from a boundary integral generated by integration by parts in the Forchheimer term. This integral involves the second fundamental form of $\partial\Omega$ and need not have a definite sign on a nonconvex domain. We control it using a trace interpolation inequality applied to $|\bm{u}_m|^{\tfrac{r}{2}}$. The resulting gradient contribution can be absorbed into the positive bulk term, while the remaining contribution is controlled by the $L^r$-bound from the higher-order energy. 

Global-in-time existence can also be obtained for initial data sufficiently close to an arbitrary local energy minimizer. 
In this case, the \L ojasiewicz--Simon approach provides a uniform-in-time control of the higher-order energy, allowing us to extend the local strong solution globally in time. 
Extending local strong solutions also requires uniqueness, which we prove by showing continuous dependence on the initial data. 
We note that the nonlinear Forchheimer term introduces an additional difficulty in this argument.
For the full range of exponents considered here, an $L^\infty(0,T;H^1(\Omega))$-bound for the velocity field alone does not provide the integrability needed to control the difference of the Forchheimer terms corresponding to two strong solutions. 
The higher-order estimate, however, yields an additional $L^\infty_t L^{3r}_x$-bound for the velocity field in three dimensions. 
This bound follows from the $L^\infty_t H^1_x$-control of $|\bm{u}|^{\tfrac{r}{2}}$ in the \textit{a priori} estimate and the Sobolev embedding theorem, and is sufficient to control the $L^\infty_t L^{2r-2}_x$-norm of the velocity field arising in the continuous-dependence estimate.

The paper is organized as follows. 
In Section \ref{preliminary_and_result}, after introducing some preliminaries and notation, we state our main results. 
Section \ref{ss_approximate} is devoted to the  solvability of an auxiliary problem, whose solutions serve as approximate solutions in the subsequent analysis. 
In Section \ref{ss_3d_proof}, we establish the existence of local strong solutions to problem \eqref{eq:chdf1}--\eqref{eq:chdfi} and the existence of global strong solutions for constant coefficients and mobilities when $r \ge 3$ in three dimensions.
We also establish uniqueness under the stated additional assumptions. 
In Section \ref{exist_global_small}, we further prove the existence and uniqueness of a global strong solution when the initial data are close to a local minimizer of the free energy. Finally, in Section \ref{ss_2d}, we show the existence of global-in-time strong solutions for arbitrarily large initial data in two dimensions. In the appendix, we present an auxiliary lemma that is used in this study. 

%%%%%%%%%%%%%%%%%%%%%%%%%%%%%%%% 
\section{Main Results}\label{preliminary_and_result}
\setcounter{equation}{0}

\subsection{Preliminaries and notation}\label{pre}
Let $\mathcal{X}$ be a (real) Banach space with norm $\left\|\cdot\right\|_{\mathcal{X}}$. We denote by $\mathcal{X}^*$ its dual space and by $\langle \cdot,\cdot \rangle_{\mathcal{X}^*,\mathcal{X}}$ the associated duality pairing.
The direction map $\bm{e}_{\mathcal{X}}:\mathcal{X} \to \mathcal{X}$ is defined by
\[
    \bm{e}_{\mathcal{X}}(\bm{r})
    \stackrel{\rm{def}}{=} \left\{
        \begin{aligned}
            & \frac{\bm{r}}{\| \bm{r} \|_{\mathcal{X}}} , &&\ \text{for}\ \bm{r} \in \mathcal{X} \backslash \{ \bm{0} \}, \\
            & \bm{0}, &&\ \text{if}\ \bm{r} = \bm{0}.
        \end{aligned}
        \right.
\]
Given a (real) Hilbert space $\mathcal{H}$, its inner product is denoted by $(\cdot,\cdot)_{\mathcal{H}}$.
Throughout this paper, we assume that $\Omega \subset \mathbb{R}^d$ ($d\in\{2,3\}$) is a bounded domain with the boundary $\partial \Omega$ of class $C^\infty$.
For the standard Lebesgue and Sobolev spaces in $\Omega$, we use the notation $L^{p}(\Omega)$, $W^{k,p}(\Omega)$ for any $p \in [1,\infty]$ and $k\in \mathbb{N}$,
equipped with the corresponding norms
$\|\cdot\|_{L^{p}(\Omega)}$, $\|\cdot\|_{W^{k,p}(\Omega)}$, respectively. The space $W^{k,p}_0(\Omega)$ denotes the closure of $C^{\infty}_0 ( \Omega )$ in $W^{k,p}(\Omega)$ and its dual space is denoted by $ W^{-k,p'}(\Omega) \stackrel{\mathrm{def}}{=}  \big(  W^{k,p}_0(\Omega)\big)^*$,
where $\frac{1}{p}+\frac{1}{p'}=1$, $p\in (1,\infty)$. 
For $p = 2$, we use the convention $H^{k}(\Omega) \stackrel{\mathrm{def}}{=}  W^{k,2}(\Omega)$.
The $L^2$-Bessel potential spaces are denoted by $H^{s}(\Omega)$, $s\in \mathbb{R}$, which are defined by the restriction of distributions in $H^s( \mathbb{R}^{d} )$ to $\Omega$. For $s >0$, $H^{s}_0 (\Omega)$ denotes the closure of $C^{\infty}_0( \Omega )$ in $H^{s}(\Omega)$.
For simplicity, the norm and the inner product of $L^2(\Omega)$ will be denoted by $\|\cdot\|$ and $(\cdot,\cdot)$, respectively, while the duality pairing between $H^1(\Omega)$ and $(H^1(\Omega))^*$ will be indicated by $\langle \cdot,\cdot \rangle$. 
In order to avoid confusion, we also use the notation $(\cdot,\cdot)_{\mathcal{D}}$ to point out the integration domain $\mathcal{D}$.
Bold letters will be used for vector-valued spaces, 
for example, $\bm{L}^p(\Omega)= L^p(\Omega;\mathbb{R}^d)$, $p \in [1,\infty]$. When there is no ambiguity, for a vector-valued function $\bm{f}:\Omega\to \mathbb{R}^d$, the expression $\bm{f} \in L^p (\Omega)$ should be understood as $\bm{f} \in \bm{L}^p (\Omega) $.

Given a measurable set $I$ of $\mathbb{R}$, we introduce the space $L^p(I;\mathcal{X})$ with $p\in [1,\infty]$, which consists of Bochner measurable $p$-integrable functions (if $p \in [1,\infty)$) or essentially bounded functions (if $p =\infty$) with values in a given Banach space $\mathcal{X}$. If $I=(a,b)$, we simply write $L^p(a,b;\mathcal{X})$. In addition, $f \in L^{p}_{\mathrm{loc}} ( [0,\infty);\mathcal{X} )$ if and only if $f \in L^{p} \big( 0,T;\mathcal{X} \big)$ for every $T\in (0,\infty)$. The space $L^p_{\mathrm{uloc}}([0,\infty);\mathcal{X})$ denotes the uniformly local variant of $L^p(0,\infty;\mathcal{X})$ consisting of all strongly measurable $f:[0,\infty)\to \mathcal{X}$ such that
\begin{equation*}
    \|f\|_{L^p_{\mathrm{uloc}}([0,\infty);\mathcal{X})} 
    \stackrel{\mathrm{def}}{=}  \mathop\mathrm{sup}\limits_{t\ge0} \|f\|_{L^p(t,t+1;\mathcal{X})} < \infty.
\end{equation*}
If $T\in(0,\infty)$, we set $L^p_{\mathrm{uloc}}([0,T);\mathcal{X})= L^p(0,T;\mathcal{X})$.
For every $k\in \mathbb{N}_+$, $f\in W^{k,p} ( 0,T;\mathcal{X} )$ if and only if $ \{ (\frac{\mathrm{d}}{\mathrm{d}t})^jf \}_{j=0}^k \subset L^{p} ( 0,T;\mathcal{X} )$, where $\frac{\mathrm{d}}{\mathrm{d t}} f$ denotes the vector-valued distributional derivative of $f$ and $(\frac{\mathrm{d}}{\mathrm{d}t})^0 f \stackrel{\mathrm{def}}{=} f$. 
Furthermore, $W^{k,p}_{\mathrm{uloc}} ( [0,\infty);\mathcal{X})$ is defined by replacing $L^{p} ( 0,T;\mathcal{X} )$ with $L^{p}_{\mathrm{uloc}} \big( [0,\infty);\mathcal{X} \big)$, and we set
\[ 
H^{k} ( 0,T;\mathcal{X} ) =W^{k,2} ( 0,T;\mathcal{X} ), 
\quad 
H^{k}_{\mathrm{uloc}} ( [0,\infty);\mathcal{X} ) = W^{k,2}_{\mathrm{uloc}} ( [0,\infty);\mathcal{X} ). 
\]
Let $I=[0,T]$ if $T \in (0,\infty)$ or $I=[0,\infty)$ if $T=\infty$. Then $BC( I;\mathcal{X} )$ denotes the Banach space of all bounded and continuous functions $f:I \to \mathcal{X}$ equipped with the supremum norm, and $BUC( I;\mathcal{X} )$ is the subspace of all bounded and uniformly continuous functions.
Furthermore, for every $k \in \mathbb{N}_+$, $BC^k( I;\mathcal{X} )$ (resp. $BUC^k( I;\mathcal{X} )$) denotes all functions $f \in BC( I;\mathcal{X} )$ (resp. $f \in BUC( I;\mathcal{X} )$) whose Fr\'echet derivatives of order no larger than $k$ exist and belong to $BC( I;\mathcal{X} )$ (resp. $BUC( I;\mathcal{X} )$).
Finally, we denote by $BC_{w}( I;\mathcal{X} )$ the topological vector space of all bounded and weakly continuous functions $f:I \to \mathcal{X}$.

For convenience, we introduce the following shorthand notation:
\begin{align*}
    & H \stackrel{\mathrm{def}}{=}  L^2\left(\Omega\right),
    \quad V \stackrel{\mathrm{def}}{=}  H^1\left(\Omega\right), 
    \quad  W \stackrel{\mathrm{def}}{=} \big\{u\in H^2(\Omega)\ \big|\ \partial_\mathbf{n}u=0\ \text{a.e. on}\ \partial\Omega\big\}.
\end{align*}
$H$ is identified with its dual, and we have the following continuous, dense, and compact embeddings:
$W \hookrightarrow V  \hookrightarrow H \hookrightarrow V^*$. In addition, $V=D \big((I-\Delta_N)^\frac{1}{2}\big)$, $W=D (I-\Delta_N)$, where $I$ is the identity operator on $H$ and $-\Delta_N: W \subset H \to H$ is the minus Laplace operator on $\Omega$ subject to the homogeneous Neumann boundary condition. For all $k \in \mathbb{N}_+$, we define $H_N^k(\Omega)$ as $ D \big((I-\Delta_N)^\frac{k}{2}\big)$. 
For every $f\in V^*$, $\overline{f}$ denotes its generalized mean value over $\Omega$, given by
$\overline{f}=|\Omega|^{-1}\langle f,1\rangle$. 
We recall the well-known Poincar\'{e}--Wirtinger inequality:
\begin{equation*}
\big\|f-\overline{f}\big\|\leq C_P \|\nabla f\|,\quad \forall\,
f\in V,
\end{equation*}
where the positive constant $C_P$ depends only on $\Omega$.
Let $A_N: V \to V^*$ denote the extension of $-\Delta_N$ on $V$, that is,
\begin{equation*}
    \langle A_N f,g \rangle 
    = \int_{\Omega} \nabla f \cdot \nabla g \,\ \mathrm{d}x, \quad \forall\,f,g \in V.
\end{equation*}
Setting
\begin{align*}
    & L^p_{(0)}(\Omega) \stackrel{\mathrm{def}}{=} \big\{f\in L^{p}(\Omega) \ \big|\ \overline{f} =0\big\},
    \quad H_{(0)} \stackrel{\mathrm{def}}{=}  L^2_{(0)}(\Omega), 
    \\
    & V_{(0)} \stackrel{\mathrm{def}}{=}  V\cap H_{(0)}, 
    \quad V^*_{(0)} \stackrel{\mathrm{def}}{=}  \big\{ L \in V^* \ \big|\   \langle L,1 \rangle = 0 \big\},
\end{align*}
it is easy to verify that the restriction of $A_N$ to $V_{(0)}$ is a linear isomorphism between $V_{(0)}$ and $V^*_{(0)}$. 
Thus, we can define the inverse operator $\mathcal{N} \stackrel{\mathrm{def}}{=}  \big( A_N|_{V_{(0)}} \big)^{-1}:V_{(0)}^* \to V_{(0)}$, and obtain the following identities
\begin{align*}
&\left\langle A_N u,\mathcal{N}L \right\rangle 
= \left\langle L,u \right\rangle,\quad \forall\, u \in V_{(0)},\ L \in V_{(0)}^*,
\\
&\left\langle L_1,\mathcal{N}L_2 \right\rangle = ( \nabla \left(\mathcal{N}L_1\right), \nabla \left(\mathcal{N}L_2\right) ),
\quad \forall\, L_1,\,L_2\ \in V_{(0)}^*.
\end{align*}
Define
\begin{align*}
    &\|L\|_* \stackrel{\mathrm{def}}{=}  \| \nabla \left(\mathcal{N}L\right)\| 
    = \sqrt{\langle L,\mathcal{N}L \rangle},\quad \forall\, L\in V^*_{(0)}, \quad \|L\|_{-1}^2 \stackrel{\mathrm{def}}{=}  \|L-\overline{L}\|_*^2 + |\overline{L}|^2, 
    \quad \forall\, L\in V^*.
\end{align*}
We observe that $\left\|\cdot\right\|_*$ and $\left\|\cdot\right\|_{-1}$ are equivalent to the standard dual norms on $V^*_{(0)}$ and $V^*$, respectively (cf. \cite{MZ2004}). 
In order to handle non-constant mobilities, we define the following weighted Laplace operators and their inverses:
\begin{align}
    & \mathcal{A}_{\phi, h} : V \to V^*, \ \langle \mathcal{A}_{\phi, h} (f),g \rangle = \int_{\Omega} m_\phi( h ) \nabla f \cdot \nabla g  \, \mathrm{d}x, \ \mathcal{N}_{\phi, h} \stackrel{\mathrm{def}}{=}  \big(  \mathcal{A}_{\phi, h} |_{V_{(0)}}  \big)^{-1} : V_{(0)}^* \to V_{(0)}, 
    \label{ahi_definition} 
    \\
    & \mathcal{A}_{\psi, h} : V \to V^*, \ \langle \mathcal{A}_{\psi, h} (f),g \rangle = \int_{\Omega} m_\psi( h ) \nabla f \cdot \nabla g  \, \mathrm{d}x, \ \mathcal{N}_{\psi, h} \stackrel{\mathrm{def}}{=}  \big(  \mathcal{A}_{\psi, h} |_{V_{(0)}}  \big)^{-1} : V_{(0)}^* \to V_{(0)}, 
    \label{asi_definition}
\end{align}
where $h:\Omega\to\mathbb{R}$ is Lebesgue measurable.
Since $m_\phi$ and $m_\psi$ are bounded and uniformly positive by (\textbf{H2}), the inverses $\mathcal{N}_{\phi,h}$ and $\mathcal{N}_{\psi,h}$ are well defined.
For more properties of the four operators above, see \cite{BB1999,CGGG2025,GOW2026}.

The following interpolation inequality will be useful in the subsequent analysis: 
\begin{equation*}
    \|f\|^2 \le \xi \|\nabla f\|^2 + C_\xi \|f\|^2_{V^*}, \quad \forall\, f \in V,
\end{equation*}
where $\xi \in (0,1)$ is an arbitrary constant and $C_\xi>0$ depends only on $\xi$ and $\Omega$.
In addition, we recall the classical Gagliardo--Nirenberg (GN for short) inequality:
\begin{equation}\label{gns_inequality}
    \| D^j f \|_{L^{p}(\Omega)} \le C \| D^m f \|^a_{L^{r}(\Omega)} \| f \|^{1-a}_{L^{q}(\Omega)} 
    + C' \| f \|_{L^{q}(\Omega)}, 
    \ \ \forall\, f \in W^{m,r}(\Omega) \cap L^{q}(\Omega),
\end{equation}
where $a\in [0,1]$, $p,q,r \in [1,\infty]$, and $j,m \in \mathbb{N}$ satisfy
\begin{align*}
    & 0 \le j < m, \quad  \frac{j}{m} \le a \le 1, \quad  \frac{1}{p} - \frac{j}{d} = a \Big( \frac{1}{r}-\frac{m}{d} \Big) + ( 1-a ) \frac{1}{q},
\end{align*}
and when $r>1$ as well as $m-j-(d/r)\in \mathbb{N}$, the constant $a$ is additionally required to be in $[j/m,1)$. In \eqref{gns_inequality}, $C$ and $C'$ are positive constants only depending on $a, \ q, \ r, \ j, \ m, \ d$ and $\Omega$.
Agmon's inequality says that for all $f \in H^{2}(\Omega)$, 
\begin{equation}\label{Agmon_inequality}
    \| f \|_{L^{\infty}(\Omega)}^2 \le \left\{
    \begin{aligned}
     & C \| f \| \| f \|_{H^{2}(\Omega)}, 
    && \mathrm{if} \ d=2,
    \\
    & C \| f \|_{H^{1}(\Omega)} \| f \|_{H^{2}(\Omega)}, && \mathrm{if} \ d=3.   
    \end{aligned}
    \right.
\end{equation}
In addition, we need the following Brezis--Gallou\"et--Wainger inequality (see 
\cite[Theorem 1]{E1989} for an open subset with the cone property and \cite[Theorem 1]{BW1980} for the whole space).
For $k,l \in \mathbb{N}_+$ and $p \in (1,\infty)$, $q \in [1,\infty]$, with $kp=d$ and $lq >d$, there exists some $C>0$ depending on $p$, $q$, $d$ and $\Omega$, such that for all $f \in W^{l,q}(\Omega) \cap W^{k,p}(\Omega)$, 
\begin{align}\label{BGW}
    \| f \|_{L^\infty(\Omega)} \le C (1 + \| f \|_{W^{k,p}(\Omega)}) \mathrm{ln}^\frac{1}{p'} \big(  e + \| f \|_{W^{l,q}(\Omega)}  \big), \quad \frac{1}{p} + \frac{1}{p'} = 1.
\end{align}

Concerning the velocity field, we define 
\[
    \bm{C}^{\infty}_{0,\sigma}(\Omega) \stackrel{\mathrm{def}}{=}  \{ \bm{f} \in C^{\infty}_0(\Omega)^d \ | \ \mathrm{div}\, \bm{f} =0\ \text{in}\ \Omega\}
\]
and 
\begin{align*}
    &\bm{L}^{p}_\sigma(\Omega) \stackrel{\mathrm{def}}{=}  \overline{\bm{C}^{\infty}_{0,\sigma}(\Omega)}^{\bm{L}^p(\Omega)}, 
    \quad 
    \bm{H}^{s}_{0,\sigma}(\Omega) \stackrel{\mathrm{def}}{=}  \overline{\bm{C}^{\infty}_{0,\sigma}(\Omega)}^{\bm{H}^s(\Omega)}, 
    \quad 
    \bm{W}^{k,p}_{0,\sigma}(\Omega) \stackrel{\mathrm{def}}{=}  \overline{\bm{C}^{\infty}_{0,\sigma}(\Omega)}^{\bm{W}^{k,p}(\Omega)},
\end{align*}
for all $p \in (1,\infty)$, $k \in \mathbb{N}_+$ and $s > 0$.
The Helmholtz projection corresponding to $\bm{L}^{2}_\sigma(\Omega)$ is denoted by $\mathbb{P}_{\sigma}$ (see, for example, \cite{S2001}). Then $\mathbb{P}_{\sigma} \bm{f} = \bm{f} - \nabla \pi$, where $\pi \in V_{(0)}$ is the solution to the (weak) Neumann problem
\[ 
( \nabla \pi, \nabla \phi )_{\Omega} 
= ( \bm{f}, \nabla \phi )_{\Omega}, 
\quad  \forall\, \phi \in C^{\infty} \big( \overline{\Omega} \big). 
\]
In addition, since $\partial \Omega$ is sufficiently smooth, we can also define the Helmholtz--Weyl decomposition of $\bm{L}^{p}(\Omega)$ for all $p \in (1,\infty)$ 
(see \cite[Section III.1]{G2011}).
When there is no confusion, $\mathbb{P}_{\sigma}$ also denotes the projection from $\bm{L}^{p}(\Omega)$ to $\bm{L}^{p}_\sigma(\Omega)$, whose null space is
\[
\big\{  
\bm{w} \in L^{p}(\Omega):\bm{w} = \nabla \pi, \ \mathrm{for} \ \mathrm{some} \ \pi \in W^{1,p}(\Omega)   \big\}.
\]
The $\mathrm{curl} \,$ operator is defined by $\nabla \times$ in three dimensions.
In two dimensions, for vectors $\bm{a}=(a_1,a_2), \bm{c}=(c_1,c_2) \in \mathbb{R}^2$ and a scalar $b \in \mathbb{R}$, we set
\[
    \bm{a}\times b = (b a_2, -b a_1), \qquad b\times\bm{a} = -\bm{a}\times b, \qquad \bm{a}\times\bm{c} = a_1 c_2 - a_2 c_1.
\]
For sufficiently smooth vector fields $\bm{v}:\Omega\to\mathbb{R}^2$ and scalar fields $b:\Omega\to\mathbb{R}$, we define the scalar curl and vector curl operators, respectively, by
\[
    \mathrm{curl} \,\bm{v}=\partial_1v_2-\partial_2v_1,
    \qquad
    \mathrm{curl} \, b=(\partial_2b,-\partial_1b),
\]
and we also write $\nabla\times\bm{v}=\mathrm{curl} \,\bm{v}$. 
For simplicity, in the two-dimensional case, we use the same notation for the scalar curl and vector curl operators, which are formal adjoints of each other.
The corresponding meaning depends on the input field.
Under these conventions, for any $\bm{v}\in \bm{H}^1(\Omega)$ and $b\in H^1(\Omega)$, the integration by parts formula holds in the same form as in three dimensions:
\begin{align*}
    \int_\Omega \bm{v}\cdot(\mathrm{curl} \, b)\,\mathrm{d}x & = \int_\Omega (\mathrm{curl} \,\bm{v})\,b\,\mathrm{d}x + \int_{\partial\Omega} b (\bm{v}\times\bm{n})\,\mathrm{d}\Gamma \\
    & = \int_\Omega (\mathrm{curl} \,\bm{v})\,b\,\mathrm{d}x + \int_{\partial\Omega} \bm{v}\cdot(\bm{n} \times b)\,\mathrm{d}\Gamma.
\end{align*}
Moreover, formally it holds that $\mathrm{curl} \,\mathrm{curl} \,\bm{v} = \nabla \mathrm{div} \, \bm{v}-\Delta\bm{v}$, and that $(\mathrm{curl} \,\bm{v})\times\bm{n}=\bm{0}$ is equivalent to $\mathrm{curl} \,\bm{v}=0$ on $\partial\Omega$.

The capital letter $C$ denotes a generic positive constant that depends on the problem's structural data. Its meaning may change from line to line and even within the same chain of computations. Specific dependence will be pointed out if necessary.
Finally, for simplicity, we set $Q_{(s,t)}=\Omega \times (s,t), \ Q_{t}=Q_{(0,t)}$, $Q=Q_{(0,\infty)}$,
and analogously $S_{(s,t)} = \partial \Omega \times (s,t)$, $S_t = S_{(0,t)}$, and $S=S_{(0,\infty)}$.

%%%%%%%%%%%%%%%%%%%%%%%%%%%%%%%%
\subsection{Statements of the main results} \label{main_result}
\label{sectionassumptions}

First, we introduce some basic assumptions (cf. \cite{ADG2013ndm,AGG2012,B1999}).

\begin{assumption} \label{assumptions} \rm
We assume that $\Omega \subset \mathbb{R}^{d}\ (d\in \{ 2,3\})$ is a bounded domain with a smooth boundary $\partial\Omega$ and $\beta>0$, $\sigma_1\geq 0$, $\sigma_2\in \mathbb{R}$, $c\in (-1,1)$, $r \in [2,\infty)$ are given constants.
In addition, we impose the following assumptions.
\medskip 
\begin{itemize}
    \item[($\mathbf{H1}$)] The potential function
    $F_{\phi}\in C([-1,1]) \cap C^2((-1,1))$ satisfies
    \begin{equation*}
        \mathop{\mathrm{lim}}\limits_{s \to -1^+} F_{\phi}'(s) = -\infty,\quad
        \mathop{\mathrm{lim}}\limits_{s \to 1^-} F_{\phi}'(s) = + \infty, \quad
        F_{\phi}''(s) \ge \theta_{\phi} > 0, \quad \forall\, s \in (-1,1),
    \end{equation*}
    where $\theta_\phi$ is a given constant.
    Similarly, $F_{\psi}\in C([0,1]) \cap C^2((0,1))$ satisfies
    \begin{equation*}
        \mathop{\mathrm{lim}}\limits_{s \to 0^+} F_{\psi}'(s) = -\infty,\quad
        \mathop{\mathrm{lim}}\limits_{s \to 1^-} F_{\psi}'(s) = + \infty, \quad
        F_{\psi}''(s) \ge \theta_{\psi} > 0, \quad \forall\, s \in (0,1).
    \end{equation*}
    We set $F_{\phi}(s) = +\infty $ whenever $|s| > 1$ and $F_{\psi}(s) = +\infty $ whenever $s \in (-\infty,0) \cup (1,+\infty)$. Without loss of generality, we take
     \begin{equation*}
         F_{\phi}(0)= F_{\phi}'(0) = F_{\psi}\left(\frac12\right) = F_{\psi}'\left(\frac12\right) = 0.
     \end{equation*}
    Besides, we assume that $G(\cdot,\cdot)\in BUC^2(\mathbb{R}^{2})$.
    \item[($\mathbf{H2}$)]
    We assume $m_{\phi}(\cdot),m_{\psi}(\cdot) \in BUC^2(\mathbb{R})$ such that 
    \[ m_{\phi,*} \stackrel{\rm{def}}{=} \inf\limits_{s \in \mathbb{R}} \{ m_{\phi}(s) \} >0 \quad  \mathrm{and} \quad m_{\psi,*} \stackrel{\rm{def}}{=} \inf\limits_{s \in \mathbb{R}} \{ m_{\psi}(s) \} >0. \]
    Their upper bounds are denoted by \[m_{\phi}^* \stackrel{\rm{def}}{=} \sup\limits_{s \in \mathbb{R}} \{ m_{\phi}(s) \}\quad \text{and}\quad m_{\psi}^* \stackrel{\rm{def}}{=} \sup\limits_{s \in \mathbb{R}} \{ m_{\psi}(s) \},\] 
    respectively. Moreover, the kinematic viscosity $\nu$ and the Forchheimer coefficient $\eta$ satisfy $\nu, \eta \in BUC^1( \mathbb{R}^2 )$ with
    \[ \nu_* \stackrel{\rm{def}}{=} \inf\limits_{(s_1,s_2) \in \mathbb{R}^2} \{ \nu(s_1,s_2) \} >0, \quad \eta_* \stackrel{\rm{def}}{=} \inf\limits_{(s_1,s_2) \in \mathbb{R}^2} \{ \eta(s_1,s_2) \} >0.\]
    Their upper bounds are
    denoted by
    \[
        \nu^* \stackrel{\rm{def}}{=} \sup\limits_{(s_1,s_2) \in \mathbb{R}^2} \{ \nu(s_1,s_2) \}, \quad \eta^* \stackrel{\rm{def}}{=} \sup\limits_{(s_1,s_2) \in \mathbb{R}^2} \{ \eta(s_1,s_2) \}.
    \]
\end{itemize}
\end{assumption}

\begin{remark} \rm
\label{mob}
Following \cite[Rem.~2.1]{OGW2025},
for the sake of convenience, in $(\mathbf{H1})$ we have absorbed the smooth (possibly non-convex) part of $F_\phi$ and $F_\psi$ into $G$.
 For instance, suppose we are dealing with the Flory--Huggins potentials
\begin{align}
\label{FH3}
    & F_\phi(s) = \frac{\theta_1}{2}  [(1+s) \ln(1+s) + (1-s) \ln(1-s)] + \frac{\widetilde{\theta}_1}{2} (1-s^2), \quad s\in(-1,1),\\
\label{FH1}
    & F_\psi(s) = \frac{\theta_2}{2}  [s \ln s + (1-s) \ln(1-s)] + \frac{\widetilde{\theta}_2}{2} s(1-s), \quad s\in(0,1).
\end{align}
  Then, we absorb the smooth concave parts (i.e. the demixing terms) into $G$ and we only keep the convex mixing entropy densities.
  Since the solution $(\phi,\psi)$ takes values in $[-1,1]\times[0,1]$, we only need the functions $m_{\phi}(\cdot)$ in $[-1,1]$, $m_{\psi}(\cdot)$ in $[0,1]$ and $\nu(\cdot,\cdot)$, $\eta(\cdot,\cdot)$, $G(\cdot,\cdot)$ in $[-1,1]\times[0,1]$. We can thus easily extend these functions to $\mathbb{R}$ or $\mathbb{R}^{2}$ in such a way that the conditions in Assumption \ref{assumptions} are fulfilled. In particular, any bivariate polynomial is admissible as an interaction energy density $G$.
\end{remark}

%%%%%%%%%%%%%%%%%%%%%%%%%%%%%%%%
Next, we introduce the notion of strong solutions to problem \eqref{eq:chdf1}--\eqref{eq:chdfi}.

\begin{definition} \label{ss_chdf_def}
Let $\Omega \subset \mathbb{R}^d$ be a bounded domain with a smooth boundary. 
Let $T \in (0,\infty]$.
Set $I=[0,\infty)$ if $T= \infty$ or $I = [0,T]$ if $T<\infty$. 
Suppose that Assumption \ref{assumptions} holds, and $\phi_0,\,\psi_0 \in W$ satisfy $\overline{\phi_0} \in (-1,1)$, $\overline{\psi_0} \in (0,1)$, $-\Delta \phi_0 + F_\phi'(\phi_0)$, $-\beta \Delta \psi_0 + F_\psi'(\psi_0) \in V$.
We call $(\bm{u},\pi,\phi,\psi,\mu_\phi,\mu_\psi)$ a strong solution to problem \eqref{eq:chdf1}--\eqref{eq:chdfi} if it satisfies 
\begin{align*}
    & \bm{u} \in L^\infty (I;\bm{L}^{3r}_\sigma(\Omega) \cap \bm{H}^1(\Omega)) \ \textit{with} \ | \bm{u} |^\frac{r}{2} \in L^\infty(I;H^1(\Omega)), \\
    & \pi \in L^\infty (I;W^{1,3r'}_{(0)}(\Omega)), 
     \\
    &\phi,\psi \in L^{\infty} (I;W^{2,p}(\Omega)\cap W ), \ \ \partial_t \phi, \partial_t \psi \in L^{\infty} (I;V^*) \cap L^2_{\mathrm{uloc}} (I;V), 
    \\
    &(\phi,\psi) \in (-1,1) \times (0,1) \ \mathrm{a.e.} \ \mathrm{in} \ \Omega \times I, 
    \\
    &\mu_\phi,\mu_\psi \in L^\infty (I;V) \cap L^2_{\mathrm{uloc}} (I;H^3(\Omega)\cap W),\quad F_\phi'(\phi),\,F_\psi'(\psi) \in L^\infty (I;L^p(\Omega)),
\end{align*}
where $p = 6 $ if $d=3$ and $p \in [2,\infty)$ if $d=2$, and it satisfies the equations \eqref{eq:chdf1}--\eqref{eq:chdf6} almost everywhere in $\Omega \times I$, the boundary conditions \eqref{eq:chdfb} almost everywhere on $\partial \Omega \times I$ and the initial conditions \eqref{eq:chdfi} almost everywhere in $\Omega$.   
\end{definition}

Our first result reads as follows.
\begin{theorem} \label{ss_chdf_exist_3d}
Let Assumption \ref{assumptions} hold. 
Suppose that $\phi_0,\psi_0 \in W$ satisfy $\overline{\phi_0} \in (-1,1)$, $\overline{\psi_0} \in (0,1)$ and $-\Delta \phi_0 + F_\phi'(\phi_0), -\beta \Delta \psi_0 + F_\psi'(\psi_0) \in V$.

(1) There exists $T_e>0$, depending on the norms of the initial data, such that problem \eqref{eq:chdf1}--\eqref{eq:chdfi} admits a strong solution in $I = [0,T_e]$ in the sense of Definition \ref{ss_chdf_def}. If additionally $\nu$, $\eta$ and $m_\phi$, $m_\psi$ are positive constants and $r\ge3$, then problem \eqref{eq:chdf1}--\eqref{eq:chdfi} admits a global-in-time strong solution in $I=[0,\infty)$.

(2) If   
\[
\| \phi_0 \|_{L^\infty(\Omega)} \le 1 - \delta_0, \quad \Big\|\psi_0 - \frac{1}{2}\Big\|_{L^\infty(\Omega)} \le \frac{1}{2} - \delta_0,
\quad \text{for some}\ \delta_0 \in \Big(0,\frac12\Big),
\]
then there exists a positive time $T_s \in I$, depending on $\delta_0$, the $H^3(\Omega)$-norms of $\phi_0$ and $\psi_0$, $\overline{\phi_0}$, $\overline{\psi_0}$ and coefficients of the system \eqref{eq:chdf1}--\eqref{eq:chdf6}, such that the strong solution satisfies
\[
\| \phi(t) \|_{L^\infty(\Omega)} \le 1 - \frac12 \delta_0, \quad \| \psi(t) - \frac{1}{2} \|_{L^\infty(\Omega)} \le \frac{1}{2} - 
\frac12 \delta_0, \quad \forall\, t \in [0,T_s].
\]
In addition, if $F_\phi \in C^3((-1,1))$, $F_\psi \in C^3((0,1))$, and $G \in C^3((-1,1) \times (0,1))$, then a strong solution with the strict separation property is unique in the associated time interval. 
\end{theorem}
\begin{remark} \label{pressure_term} \rm 
The pressure $\pi$ in \eqref{eq:chdf1} can be recovered by the de Rham theorem, as in \cite{GOPW2026}. Its regularity follows directly from the elliptic regularity theory.
\end{remark}

Next, we show that the strong solution exists globally if its initial data are small perturbations of a minimizer of the free energy $E_{\mathrm{free}}$. 
\begin{definition} \label{minimizer_def}
(Energy minimizer) Let $m_1 \in (-1,1)$ and $m_2 \in (0,1)$ be two given constants. Define the admissible set
\[
\mathcal{Z}_{m_1,m_2} \stackrel{\mathrm{def}}{=} \Big\{ \phi,\psi \in V \cap L^\infty(\Omega) \ \Big| \ \| \phi \|_{L^\infty(\Omega)} \le 1, \ \Big\|\psi -\frac{1}{2}\Big\|_{L^\infty(\Omega)} \le \frac{1}{2}, \ \overline{\phi} = m_1, \ \overline{\psi} = m_2 \Big\}.
\]
A pair $(\phi_*,\psi_*) \in \mathcal{Z}_{m_1,m_2}$ is called a local minimizer of $E_{\mathrm{free}}$ in $\mathcal{Z}_{m_1,m_2}$ if there exists a positive constant $\lambda_*$ such that for any $(\phi,\psi) \in \mathcal{Z}_{m_1,m_2}$ with $\| \phi-\phi_* \|_V + \| \psi-\psi_* \|_V \le \lambda_*$, it holds
\begin{equation} \label{eq:minimum_condition}
    E_{\mathrm{free}}(\phi_*,\psi_*) \le E_{\mathrm{free}}(\phi,\psi).
\end{equation}
If \eqref{eq:minimum_condition} is satisfied for all $(\phi,\psi) \in \mathcal{Z}_{m_1,m_2}$, we call $(\phi_*,\psi_*)$ a global minimizer of $E_{\mathrm{free}}$ in $\mathcal{Z}_{m_1,m_2}$.
\end{definition}

Proposition \ref{minimizer_existence_regularity} establishes the existence and regularity of energy minimizers. 
Then we have the following result.
\begin{theorem} \label{ss_chdf_exist_global_small}
Let Assumption \ref{assumptions} hold. 
In addition, we assume that $F_\phi$ is real analytic in $(-1,1)$, $F_\psi$ is real analytic in $(0,1)$, and $G$ is real analytic in $(-1,1) \times (0,1)$. 
Let $(\phi_*,\psi_*)$ be a local minimizer of $E_{\mathrm{free}}$ in $\mathcal{Z}_{m_1,m_2}$ for given $m_1 \in (-1,1)$, with $\sigma_1(m_1 - c)=0$, and $m_2 \in (0,1)$. 
Assume that the initial data $(\phi_0,\psi_0)$ satisfy $\phi_0, \psi_0 \in H^3(\Omega) \cap W$, $\overline{\phi_0}\in (-1,1)$ (resp. $\overline{\phi_0} = m_1$) if $\sigma_1>0$   (resp. $\sigma_1=0$), $\overline{\psi_0} = m_2$, and $\|\phi_0\|_{H^3(\Omega)}\leq M$, $\|\psi_0\|_{H^3(\Omega)}\leq M$ for some $M>0$. For any $\epsilon >0$, there exist small constants $r_\phi,\, r_\psi \in (0,1)$ depending on $\epsilon$ and $M$, such that, if
\begin{align*}
    & \|  \phi_0 - \phi_*  \|_V \le r_\phi, \quad \|  \psi_0 - \psi_*  \|_V \le r_\psi, 
\end{align*}
then problem \eqref{eq:chdf1}--\eqref{eq:chdfi} admits
a unique global strong solution $(\bm{u},\pi,\phi,\psi,\mu_\phi,\mu_\psi)$ with the following additional properties
\begin{align*}
    & \phi,\psi \in BC([0,\infty);H^3(\Omega)\cap W) \cap L^2_{\mathrm{uloc}} (0,\infty;H^5(\Omega)), 
    \\
    & \mu_\phi,\mu_\psi \in BC ([0,\infty);V),
\end{align*}
and satisfying the conditions
\begin{align*}
    &\| \phi(t)-\phi_* \|_W \le \epsilon, \quad \| \psi(t)-\psi_* \|_W \le \epsilon, \quad  \forall\, t \ge 0.
\end{align*}
\end{theorem}
\begin{remark} \rm 
Theorem \ref{ss_chdf_exist_global_small} yields the Lyapunov stability of local energy minimizers. The analyticity of the nonlinear functions $F_\phi, F_\psi, G$ is crucial for applying the {\L}ojasiewicz--Simon approach (cf. \cite{AW2007,GrasPoi2025,O2025}), which is key to estimating the distance between the strong solution and the given energy minimizer. Moreover, in view of Proposition \ref{minimizer_existence_regularity}, for sufficiently small $r_\phi,\, r_\psi \in (0,1)$, we find  $0<\delta_1<\delta_0<1/2$ 
such that 
\begin{align*}
&\| \phi_0 \|_{L^\infty(\Omega)} \le 1 - \delta_0, \quad \big\| \psi_0-\frac{1}{2} \big\|_{L^\infty(\Omega)} \le \frac{1}{2} - \delta_0, 
\\
&\| \phi(t) \|_{L^\infty(\Omega)} \le 1 - \delta_1, \quad \big\| \psi(t)-\frac{1}{2} \big\|_{L^\infty(\Omega)} \le \frac{1}{2} - \delta_1,\quad \forall\, t\ge 0,
\end{align*}
namely, $\phi$ and $\psi$ satisfy the uniform strict separation property for all times. This property plays an essential role in the proof of Theorem \ref{ss_chdf_exist_global_small}. 
\end{remark}

When $d=2$, we prove the existence of a global-in-time strong solution for arbitrarily large initial data. Inspired by  \cite{GalPoi2025}, we impose the following structural assumptions on $F_\phi$ and $F_\psi$ that ensure the strict separation property of $\phi$ and $\psi$:
\begin{itemize}
    \item[($\mathbf{H1*}$)] There exists a constant $\rho> 1/2$ such that as $\delta \to 0^+$, it holds
    \[
        \frac{1}{F_\phi'(1-2\delta)} = O\left(\frac{1}{| \ln \delta|^{\rho}}\right),
        \quad \frac{1}{F_\phi'(-1+2\delta)} = O\left(\frac{1}{| \ln \delta|^{\rho}}\right),
    \]
    and
    \[
        \frac{1}{F_\psi'(1-2\delta)} = O\left(\frac{1}{| \ln \delta|^{\rho}}\right),\quad \frac{1}{F_\psi'(2\delta)} = O\left(\frac{1}{| \ln \delta|^{\rho}}\right).
    \]
\end{itemize}
Then, our last result reads as follows:

\begin{theorem} \label{ss_chdf_exist_2d}
Let $d=2$, and let Assumption \ref{assumptions} hold. 
For any initial data $\phi_0,\psi_0 \in W$ satisfying $\overline{\phi_0} \in (-1,1)$, $\overline{\psi_0} \in (0,1)$ and $-\Delta \phi_0 + F_\phi'(\phi_0)$, $-\beta \Delta \psi_0 + F_\psi'(\psi_0) \in V$, 
problem \eqref{eq:chdf1}--\eqref{eq:chdfi} admits
a global strong solution $(\bm{u},\pi,\phi,\psi,\mu_\phi,\mu_\psi)$
in $[0,\infty)$ in the sense of Definition \ref{ss_chdf_def}.
If $(\mathbf{H1*})$ holds, then we also have 
\begin{align*}
    & \phi, \psi \in L^\infty(0,\infty;H^3(\Omega))\ \textit{and}\ \phi \in [-1+\delta_{\phi},1-\delta_{\phi}], \ \psi \in [\delta_{\psi},1-\delta_{\psi}] \ \mathrm{a.e. \ in} \ Q, \\
    & \mu_\phi, \mu_\psi \in BC ([0,\infty);V), \quad \partial_t \mu_\phi, \partial_t \mu_\psi \in L^2 (0,\infty; V^*),
\end{align*}
where $\delta_{\phi}, \delta_{\psi} \in (0,1/2)$ depend on $\Omega$, $\| G \|_{C^2( [-1,1]\times[0,1] )}$, $\beta$, $F_\phi$, $F_\psi$, $\| \mu_\phi \|_{L^\infty ( 0,\infty;V )}$ and $\| \mu_\psi \|_{L^\infty ( 0,\infty;V )}$.
If, in addition to $(\mathbf{H1*})$, $F_\phi \in C^3( (-1,1) )$, $F_\psi \in C^3( (0,1) )$ and $G \in C^3( (-1,1) \times (0,1) )$, 
then the strong solution is unique.
\end{theorem}
\begin{remark} \label{CHD} \rm 
Note that the Forchheimer term plays a basic role only in Theorem \ref{ss_chdf_exist_3d} for the global existence. The other results also hold for Darcy's law (i.e. $r=2$).
\end{remark}

%%%%%%%%%%%%%%%%%%%%%%%%%%%%%%%%
\section{Construction of Approximate Solutions}
\label{ss_approximate}
\setcounter{equation}{0}

In this section, we construct (local-in-time) approximate solutions by solving an auxiliary problem. 
Let us define the Stokes operator $\bm{A}$  (see \cite{CDG2002,XX2007,XX2013})
    \begin{align*}
        & \mathcal{D} (\bm{A}) \stackrel{\rm{def}}{=} \{  \bm{u} \in \bm{H}^2(\Omega) \cap \bm{L}^2_\sigma(\Omega): \mathrm{curl} \,\bm{u} \times \bm{n} = \bm{0} \ \text{ a.e. on $\partial \Omega$}  \}, 
        \\
        & \bm{A} \bm{u} \stackrel{\rm{def}}= - \Delta \bm{u} = \mathrm{curl} \, \mathrm{curl} \, \bm{u}  \in \bm{L}^2_\sigma(\Omega), \ \forall\, \bm{u} \in \mathcal{D} (\bm{A}).
    \end{align*} 
As in \cite{GOW2026,G2022}, we consider the family of eigenvalues $\{ 1 + \lambda_{\bm{A}}^j \}_{j=1}^\infty$ and the corresponding eigenfunctions $\{  \bm{w}_j  \}_{j=1}^\infty$ of the operator $I + \bm{A}$.
For every $m\in \mathbb{N}_+$, we set the finite dimensional space $\bm{V}_m \stackrel{\rm{def}}{=} \mathrm{span} \{ \bm{w}_1, \dots, \bm{w}_m \}$, endowed with the  $\bm{L}^2_\sigma (\Omega)$-norm. 
We denote the orthogonal projection on $\bm{V}_m$ with respect to the inner product in $\bm{L}^2_\sigma (\Omega)$ by $\mathbb{P}^{\bm{V}}_m$.

Our approximate system is the following
\begin{numcases}{}
    (  \nu (\phi_{m,n},\psi_{m,n}) \bm{u}_{m,n} + \eta (\phi_{m,n},\psi_{m,n}) |\bm{u}_{m,n}|^{r-2} \bm{u}_{m,n} , \bm{\theta} ) \nonumber \\
    \quad = (  \mu_{\phi,m,n} \nabla \phi_{m,n} + \mu_{\psi,m,n} \nabla \psi_{m,n} , \bm{\theta} ) \label{eq:semigalerkin_chdf1}, & $\forall \bm{\theta} \in \bm{V}_m$, in $[0,T]$, \\
    \partial_t \phi_{m,n} + \bm{u}_{m,n} \cdot \nabla \phi_{m,n} + \sigma_1 \big( \overline{\phi_{m,n}} -c \big) = \mathrm{div} (m_{\phi} (\phi_{m,n}) \nabla \mu_{\phi,m,n}) \label{eq:semigalerkin_chdf2}, & a.e. in $\Omega \times (0,T)$, \\
    \mu_{\phi,m,n} = -\Delta \phi_{m,n} + \sigma_2 \mathcal{N} \big( \phi_{m,n} - \overline{\phi_{m,n}} \big) + F_{\phi,n,\kappa}'(\phi_{m,n}) + \partial_{\phi} G_\kappa(\phi_{m,n},\psi_{m,n}) \label{eq:semigalerkin_chdf3}, & a.e. in $\Omega \times (0,T)$, \\
    \partial_t \psi_{m,n} + \bm{u}_{m,n} \cdot \nabla \psi_{m,n} = \mathrm{div} (m_{\psi} (\psi_{m,n}) \nabla \mu_{\psi,m,n}) \label{eq:semigalerkin_chdf4}, & a.e. in $\Omega \times (0,T)$, \\
    \mu_{\psi,m,n} = -\beta \Delta \psi_{m,n} + F_{\psi,n,\kappa}'(\psi_{m,n}) + \partial_\psi G_\kappa(\phi_{m,n},\psi_{m,n}) \label{eq:semigalerkin_chdf5}, & a.e. in $\Omega \times (0,T)$.
\end{numcases}
The system is subject to the boundary conditions
\begin{equation*}
    \partial_{\mathbf{n}} \phi_{m,n} = \partial_{\mathbf{n}} \psi_{m,n} = \partial_{\mathbf{n}} \mu_{\phi,m,n} = \partial_{\mathbf{n}} \mu_{\psi,m,n} =0, \quad \text{ a.e. on } \partial\Omega \times (0,T),
\end{equation*}
and to the initial conditions:
\begin{equation}
    \phi_{m,n}|_{t=0} = \phi_{0,k}(x), \quad \psi_{m,n}|_{t=0} = \psi_{0,k}(x), \quad \text{a.e. in} \ \Omega. \label{eq:semigalerkin_chdfi}
\end{equation}
Here the approximate initial data $\phi_{0,k} \in H^3(\Omega)\cap W$ and $\psi_{0,k}\in H^3(\Omega)\cap W$ are defined as in \cite[Section 3.1]{GOW2026}.
In particular, they are  strictly separated from pure phases, namely, $\phi_{0,k} \in [-1+\delta_k,1-\delta_k]$ and $\psi_{0,k} \in [\delta_k,1-\delta_k]$ for some $\delta_k \in (0,\tfrac{1}{2}]$. 
For their regularity and convergence properties as $k \to \infty$, we refer the reader to \cite[Lemma 3.1]{GOW2026}.
The approximate energy density functions $F_{\phi,n,\kappa}$, $F_{\psi,n,\kappa}$ and $G_\kappa$ are defined by
\begin{align*}
    & F_{\phi,n,\kappa}(s) \stackrel{\mathrm{def}}{=} F_{\phi,n} * \xi_1^\kappa (s) - F_{\phi,n} * \xi_1^\kappa (0) - F_{\phi,n}' * \xi_1^\kappa (0) s, \quad \forall s \in \mathbb{R}, \\ & F_{\psi,n,\kappa}(s) \stackrel{\mathrm{def}}{=} F_{\psi,n} * \xi_1^\kappa (s) - F_{\psi,n} * \xi_1^\kappa \Big( \frac{1}{2} \Big) - F_{\psi,n}' * \xi_1^\kappa \Big( \frac{1}{2} \Big) \Big( s-\frac{1}{2} \Big),  \quad \forall s \in \mathbb{R}, \\ & G_\kappa \stackrel{\mathrm{def}}{=} G * \xi_2^\kappa,
\end{align*}
where $\xi_1^\kappa$ (resp. $\xi_2^\kappa$) is a family of mollifiers on $\mathbb{R}$ (resp. $\mathbb{R}^2$) with the support in $[-\kappa,\kappa]$ (resp. $[-\kappa,\kappa]^2$).
The approximate potential functions $F_{\phi,n}$, $F_{\psi,n}$ without mollification are defined by
\[
    F_{\phi,n}(s)
    = \left\{
        \begin{aligned}
            & F_\phi(s_{-1,n}^\phi) + F_\phi'(s_{-1,n}^\phi) (s- s_{-1,n}^\phi) + \frac{F_\phi''(s_{-1,n}^\phi)}{2} (s- s_{-1,n}^\phi)^2 , &&\ s < s_{-1,n}^\phi, \\
            & F_\phi(s), &&\ s \in \left[ s_{-1,n}^\phi,s_{1,n}^\phi \right],\\
            & F_\phi(s_{1,n}^\phi) + F_\phi'(s_{1,n}^\phi) (s - s_{1,n}^\phi) + \frac{F_\phi''(s_{1,n}^\phi)}{2} (s-s_{1,n}^\phi)^2, &&\ s > s_{1,n}^\phi,
        \end{aligned}
        \right.
\]
\[
    F_{\psi,n}(s)
    = \left\{
        \begin{aligned}
            & F_\psi(s_{0,n}^\psi) + F_\psi'(s_{0,n}^\psi) (s- s_{0,n}^\psi) + \frac{F_\psi''(s_{0,n}^\psi)}{2} (s- s_{0,n}^\psi)^2 , &&\ s < s_{0,n}^\psi, \\
            & F_\psi(s), &&\ s \in \left[ s_{0,n}^\psi,s_{1,n}^\psi \right],\\
            & F_\psi(s_{1,n}^\psi) + F_\psi'(s_{1,n}^\psi) (s-s_{1,n}^\psi) + \frac{F_\psi''(s_{1,n}^\psi)}{2} (s-s_{1,n}^\psi)^2, &&\ s > s_{1,n}^\psi,
        \end{aligned}
        \right.
\]
\begin{align*}
    & s_{-1,n}^\phi \stackrel{\mathrm{def}}{=} \sup \{  s \in (-1,0]: F_\phi''(s) \ge n  \}, \quad s_{1,n}^\phi \stackrel{\mathrm{def}}{=} \inf \{  s \in [0,1): F_\phi''(s) \ge n  \}, \\
    & s_{0,n}^\psi \stackrel{\mathrm{def}}{=} \sup \{  s \in (0,\tfrac{1}{2}]: F_\psi''(s) \ge n  \}, \quad s_{1,n}^\psi \stackrel{\mathrm{def}}{=} \inf \{  s \in [\tfrac{1}{2},1): F_\psi''(s) \ge n  \}.
\end{align*}
For $n \ge \theta_\phi$, it follows that 
\[
    F_{\phi,n}(s)
    \ge \left\{
        \begin{aligned}
            & \frac{n}{8} s^2  , &&\ s \in (-\infty,-2] \cup [2,+\infty), \\
            & \frac{\theta_\phi}{2} s^2, &&\ s \in \left(-2,2\right),
        \end{aligned}
        \right.
\]
where we use the facts $F_{\phi,n}(0) = F_{\phi,n}'(0)=0$ as well as
\begin{align*}
    & F_{\phi,n}''(s) \ge \theta_\phi, \quad \forall\, s \in \mathbb{R} \quad \text{and}\ \  F_{\phi,n}''(s) \ge n, \quad \forall\, s \in (-\infty,-1] \cup [1,+\infty).
\end{align*}
Without loss of generality, we assume $\kappa \le \min \{ 1-s_{1,n}^\phi,s_{-1,n}^\phi+1 \}$.
Then, it also holds
\[
    F_{\phi,n,\kappa}(s)
    \ge \left\{
        \begin{aligned}
            & \frac{n}{8} s^2  , &&\ s \in (-\infty,-2] \cup [2,+\infty), \\
            & \frac{\theta_\phi}{2} s^2, &&\ s \in \left(-2,2\right).
        \end{aligned}
        \right.
\]
By the embedding $V \hookrightarrow H \hookrightarrow V^*$ and the Poincar\'e--Wirtinger inequality, there exists a positive constant $C_P$ only depending on $\Omega$ such that
\[
    \| f - \overline{f} \|_* \le C_P \| f \|, \quad \forall f \in H,
\]
which implies
\[
    | \sigma_2 \| f - \overline{f} \|_*^2 | \le | \sigma_2 | C_P^2 \| f \|^2.
\]
Therefore, for $n \ge \theta_\phi + 8 | \sigma_2 | C_P^2$, we have
\begin{equation} \label{potential_control_ohtakawasaki}
    | \sigma_2 \| f - \overline{f} \|_*^2 | \le \int_\Omega F_{\phi,n,\kappa}(f) \, \mathrm{d} x + 4 | \sigma_2 | C_P^2 | \Omega |.
\end{equation}
Indeed, without loss of generality, we assume $n \ge \theta_\phi + 8 | \sigma_2 | C_P^2$ and 
\[
\max \{  1 - s_{1,n}^\phi, 1 + s_{-1,n}^\phi, 1 - s_{1,n}^\psi, s_{0,n}^\psi  \} \le \delta_k,
\]
so that
\[
    F_{\phi,n}^{(i)} (\phi_{0,k}) = F_\phi^{(i)}(\phi_{0,k}), \ F_{\psi,n}^{(i)} (\psi_{0,k}) = F_\psi^{(i)}(\psi_{0,k}),
\]
for all $i \in \{ 0,1,2 \}$ and the estimate \eqref{potential_control_ohtakawasaki} holds.

We now prove the existence and uniqueness of a (local-in-time) strong solution to the above approximate problem.
\begin{proposition} \label{approximate_solution_exists}
Assume that $\phi_{0,k},\psi_{0,k} \in H^3(\Omega) \cap W$.
Then there exists some $T_e' >0$, which may depend on the initial data, $n$ and $m$, such that system \eqref{eq:semigalerkin_chdf1}--\eqref{eq:semigalerkin_chdf5} admits a unique solution $(\bm{u}_{m,n},\phi_{m,n},\psi_{m,n},\mu_{\phi,m,n},\mu_{\psi,m,n})$ in $[0,T_e')$ satisfying the corresponding boundary conditions and the initial conditions, with the following regularity properties
\begin{align*}
    & \bm{u}_{m,n} \in H^1\big(  0,T_0;\bm{V}_m  \big), \\
    & \phi_{m,n}, \psi_{m,n} \in L^\infty \big(  0,T_0;H^3(\Omega)\cap W  \big) \cap H^1 \big(  0,T_0;V  \big), \\
    & \mu_{\phi,m,n}, \mu_{\psi,m,n} \in L^2 \big(  0,T_0;H^3(\Omega)\cap W \big) \cap C\big( [0,T_0];V \big) \cap H^1 \big(  0,T_0;V^*  \big),
\end{align*}
for all $T_0 \in (0,T_e')$.
\end{proposition}

\noindent \textbf{Proof.} 
The proof consists of three steps.

\noindent \textbf{Step 1: uniqueness.}
Suppose that $(\bm{u}_i,\phi_i,\psi_i,\mu_{\phi,i},\mu_{\psi,i})$, $i \in \{ 1,2 \}$, are solutions to the approximate system \eqref{eq:semigalerkin_chdf1}--\eqref{eq:semigalerkin_chdf5}, originating from $(\phi_0^i,\psi_0^i)$, $i \in \{ 1,2 \}$, with the listed regularities.
Setting
\[\big(  \widetilde{\bm{u}},\widetilde{\phi},\widetilde{\psi},\widetilde{\mu_\phi},\widetilde{\mu_\psi}  \big) = (\bm{u}_1,\phi_1,\psi_1,\mu_{\phi,1},\mu_{\psi,1}) - (\bm{u}_2,\phi_2,\psi_2,\mu_{\phi,2},\mu_{\psi,2}),\]
we have
\begin{align*}
    & \frac{1}{2} \frac{\mathrm{d}}{\mathrm{d}t} \big(  \big\| \nabla \widetilde{\phi} \big\|^2 + \beta \big\| \nabla \widetilde{\psi} \big\|^2  \big) + \int_\Omega m_\phi(\phi_1) \big| \nabla \Delta \widetilde{\phi} \big|^2 + \beta^2 m_\psi(\psi_1) \big| \nabla \Delta \widetilde{\psi} \big|^2 \ \mathrm{d} x \\
    & \qquad + \int_\Omega \nu(\phi_1,\psi_1) \big| \widetilde{\bm{u}} \big|^2 + \eta(\phi_1,\psi_1) \big(  | \bm{u}_1 |^{r-2} \bm{u}_1 - | \bm{u}_2 |^{r-2} \bm{u}_2  \big) \cdot \widetilde{\bm{u}} \ \mathrm{d} x \\
    & =  \int_\Omega \Delta \widetilde{\phi} \ \bm{u}_2 \cdot \nabla \widetilde{\phi} + \beta \Delta \widetilde{\psi} \ \bm{u}_2 \cdot \nabla \widetilde{\psi} \ \mathrm{d} x \\
    & \qquad - \int_\Omega (\nu(\phi_1,\psi_1)-\nu(\phi_2,\psi_2)) \bm{u}_2 \cdot \widetilde{\bm{u}} + (\eta(\phi_1,\psi_1)-\eta(\phi_2,\psi_2)) | \bm{u}_2 |^{r-2} \bm{u}_2 \cdot \widetilde{\bm{u}} \ \mathrm{d} x \\
    & \qquad + \int_\Omega \big(  F_{\phi,n,\kappa}'(\phi_1) - F_{\phi,n,\kappa}'(\phi_2) + \partial_\phi G_\kappa(\phi_1,\psi_1) - \partial_\phi G_\kappa(\phi_2,\psi_2) + \sigma_2 \mathcal{N} \big( \widetilde{\phi} - \overline{\widetilde{\phi}} \big)  \big) \nabla \phi_1 \cdot \widetilde{\bm{u}} \ \mathrm{d} x \\
    & \qquad + \int_\Omega \big(  F_{\psi,n,\kappa}'(\psi_1) - F_{\psi,n,\kappa}'(\psi_2) + \partial_\psi G_\kappa(\phi_1,\psi_1) - \partial_\psi G_\kappa(\phi_2,\psi_2)  \big) \nabla \psi_1 \cdot \widetilde{\bm{u}} \ \mathrm{d} x \\
    & \qquad + \int_\Omega \mu_{\phi,2} \nabla \widetilde{\phi} \cdot \widetilde{\bm{u}} + \mu_{\psi,2} \nabla \widetilde{\psi} \cdot \widetilde{\bm{u}} \ \mathrm{d} x \\
    & \qquad + \int_\Omega \big(  m_\phi(\phi_1)-m_\phi(\phi_2)  \big) \nabla \mu_{\phi,2} \cdot \nabla \Delta \widetilde{\phi} + \beta \big(  m_\psi(\psi_1)-m_\psi(\psi_2)  \big) \nabla \mu_{\psi,2} \cdot \nabla \Delta \widetilde{\psi} \ \mathrm{d} x \\
    & \qquad + \int_\Omega m_\phi(\phi_1) \big(  F_{\phi,n,\kappa}''(\phi_1) \nabla \phi_1 - F_{\phi,n,\kappa}''(\phi_2) \nabla \phi_2  \big) \cdot \nabla \Delta \widetilde{\phi} \ \mathrm{d} x \\
    & \qquad + \beta \int_\Omega m_\psi(\psi_1) \big(  F_{\psi,n,\kappa}''(\psi_1) \nabla \psi_1 - F_{\psi,n,\kappa}''(\psi_2) \nabla \psi_2  \big) \cdot \nabla \Delta \widetilde{\psi} \ \mathrm{d} x \\
    & \qquad + \int_\Omega \sigma_2 m_\phi(\phi_1) \nabla \mathcal{N} \Big(  \widetilde{\phi} - \overline{\widetilde{\phi}}  \Big) \cdot \nabla \Delta \widetilde{\phi} \ \mathrm{d} x \\
    & \qquad + \int_\Omega m_\phi(\phi_1) \big(  \nabla (\partial_\phi G_\kappa(\phi_1,\psi_1)) - \nabla (\partial_\phi G_\kappa(\phi_2,\psi_2))  \big) \cdot \nabla \Delta \widetilde{\phi} \ \mathrm{d} x \\
    & \qquad + \beta \int_\Omega m_\psi(\psi_1) \big(  \nabla (\partial_\psi G_\kappa(\phi_1,\psi_1)) - \nabla (\partial_\psi G_\kappa(\phi_2,\psi_2))  \big) \cdot \nabla \Delta \widetilde{\psi} \ \mathrm{d} x \\
    & \quad \le C \| \bm{u}_2 \|_V \big(  \| \nabla \widetilde{\phi} \|^{\frac{5}{4}} \|  \nabla \Delta \widetilde{\phi}  \|^{\frac{3}{4}} + \| \nabla \widetilde{\psi} \|^{\frac{5}{4}} \|  \nabla \Delta \widetilde{\psi}  \|^{\frac{3}{4}}  \big) \\
    & \qquad + C \big(  \| \bm{u}_2 \| + \| \bm{u}_2 \|_{\bm{L}^{2r-2}(\Omega)}^{r-1}  \big) \| \widetilde{\bm{u}} \| \big(  \|  \nabla \widetilde{\phi}  \|^{\frac{3}{4}} \|  \nabla \Delta \widetilde{\phi}  \|^{\frac{1}{4}} + \|  \nabla \widetilde{\psi}  \|^{\frac{3}{4}} \|  \nabla \Delta \widetilde{\psi}  \|^{\frac{1}{4}} + \big| \overline{\widetilde{\phi}} \big| + \big| \overline{\widetilde{\psi}} \big|  \big) \\
    & \qquad + C \big(  \| \nabla \phi_1 \|_{L^\infty(\Omega)} + \| \nabla \psi_1 \|_{L^\infty(\Omega)}  \big) \| \widetilde{\bm{u}} \| \big(  \big| \overline{\widetilde{\phi}} \big| + \big| \overline{\widetilde{\psi}} \big| + \| \nabla \widetilde{\phi} \| + \|  \nabla \widetilde{\psi}  \|  \big) \\
    & \qquad + C \big(  \| \mu_{\phi,2} \|_{L^6(\Omega)} + \| \mu_{\psi,2} \|_{L^6(\Omega)}  \big) \| \widetilde{\bm{u}} \| \big(  \|  \nabla \widetilde{\phi}  \|^{\frac{3}{4}} \|  \nabla \Delta \widetilde{\phi}  \|^{\frac{1}{4}} + \|  \nabla \widetilde{\psi}  \|^{\frac{3}{4}} \|  \nabla \Delta \widetilde{\psi}  \|^{\frac{1}{4}}  \big) \\
    & \qquad + C \big(  \| \nabla \mu_{\phi,2} \| + \| \nabla \mu_{\psi,2} \|  \big) \big(  \| \nabla \Delta \widetilde{\phi} \| + \| \nabla \Delta \widetilde{\psi} \|  \big) \big(  \|  \nabla \widetilde{\phi}  \|^{\frac{3}{4}} \|  \nabla \Delta \widetilde{\phi}  \|^{\frac{1}{4}} + \|  \nabla \widetilde{\psi}  \|^{\frac{3}{4}} \|  \nabla \Delta \widetilde{\psi}  \|^{\frac{1}{4}} + \big| \overline{\widetilde{\phi}} \big| + \big| \overline{\widetilde{\psi}} \big|  \big) \\
    & \qquad + C \| \nabla \phi_1 \|_{L^\infty(\Omega)} \| \nabla \Delta \widetilde{\phi} \| \big(  \big| \overline{\widetilde{\phi}} \big| + \| \nabla \widetilde{\phi} \|  \big) + C \| \nabla \psi_1 \|_{L^\infty(\Omega)} \| \nabla \Delta \widetilde{\psi} \| \big(  \big| \overline{\widetilde{\psi}} \big| + \| \nabla \widetilde{\psi} \|  \big) \\
    & \qquad + C \| \nabla \widetilde{\phi} \| \| \nabla \Delta \widetilde{\phi} \| + C \| \nabla \widetilde{\psi} \| \| \nabla \Delta \widetilde{\psi} \| \\
    & \quad \le C \Big(  \big\| \nabla \widetilde{\phi} \big\|^2 + \big\|  \nabla \widetilde{\psi} \big\|^2 + \big| \overline{\widetilde{\phi}} \big|^2 + \big| \overline{\widetilde{\psi}} \big|^2  \Big) + \frac{1}{2} \int_\Omega m_\phi(\phi_1) \big| \nabla \Delta \widetilde{\phi} \big|^2 + \beta^2 m_\psi(\psi_1) \big| \nabla \Delta \widetilde{\psi} \big|^2 + \nu(\phi_1,\psi_1) \big| \widetilde{\bm{u}} \big|^2 \ \mathrm{d} x.
\end{align*}
This implies the continuous dependence of the solutions on the initial data and, in particular, the uniqueness of a strong solution to the approximate problem \eqref{eq:semigalerkin_chdf1}--\eqref{eq:semigalerkin_chdfi} with the required properties. 
Here we use the $C^3$-regularity of $F_{\phi,n,\kappa}$, $F_{\psi,n,\kappa}$, and $G_\kappa$.
\smallskip

\noindent \textbf{Step 2: solving an auxiliary system.}
To begin with, for the Laplace operator $-\Delta: W \subset H \to H$, its eigenvalues and the corresponding eigenfunctions will be denoted by $\{ \lambda_{\Delta}^j \}_{j=1}^\infty$ and $\{ w_j \}_{j=1}^\infty$.
For every $l\in \mathbb{N}_+$, the subspace $X_l$ of $H$ is given by $X_l \stackrel{\rm{def}}{=} \mathrm{span} \{ w_1, \dots, w_l \}$, endowed with the norm of $H$. 
The orthogonal projection on $X_l$ with respect to the inner product in $H$ is denoted by $\mathbb{P}^X_l$.
Then, we introduce an auxiliary system as the approximation of the problem \eqref{eq:semigalerkin_chdf1}--\eqref{eq:semigalerkin_chdfi}. 
This is given by
\begin{numcases}{}
    (  \alpha \partial_t \bm{u}^{l,\alpha} + \nu (\phi^{l,\alpha},\psi^{l,\alpha}) \bm{u}^{l,\alpha} + \eta (\phi^{l,\alpha},\psi^{l,\alpha}) \big| \bm{u}^{l,\alpha} \big|^{r-2} \bm{u}^{l,\alpha} , \bm{\theta} ) \nonumber \\
    \quad = (  \mu_\phi^{l,\alpha} \nabla \phi^{l,\alpha} + \mu_\psi^{l,\alpha} \nabla \psi^{l,\alpha} , \bm{\theta} ) \label{eq:second_approximate_chdf1}, & $\forall \bm{\theta} \in \bm{V}_m$, in $ [0,T]$, \\
    \big(  \partial_t \phi^{l,\alpha} + \bm{u}^{l,\alpha} \cdot \nabla \phi^{l,\alpha} + \sigma_1 \big( \overline{\phi^{l,\alpha}} -c \big) , \zeta_1 \big) \nonumber \\
    \quad = - \big(  m_{\phi} (\phi^{l,\alpha}) \nabla \mu_\phi^{l,\alpha} , \nabla \zeta_1 \big) \label{eq:second_approximate_chdf2}, & $\forall \zeta_1 \in X_l$, in $ [0,T]$, \\
    ( \mu_\phi^{l,\alpha}, \zeta_2 ) = ( \nabla \phi^{l,\alpha}, \nabla \zeta_2 ) + \sigma_2 \big(  \mathcal{N} \big( \phi^{l,\alpha} - \overline{\phi^{l,\alpha}} \big) , \zeta_2 \big) \nonumber \\
    \qquad \qquad \quad + \big(  F_{\phi,n,\kappa}'(\phi^{l,\alpha}) + \partial_{\phi} G_\kappa(\phi^{l,\alpha},\psi^{l,\alpha}) , \zeta_2 \big) \label{eq:second_approximate_chdf3}, & $\forall \zeta_2 \in X_l$, in $ [0,T]$, \\
    \big(  \partial_t \psi^{l,\alpha} + \bm{u}^{l,\alpha} \cdot \nabla \psi^{l,\alpha} , \zeta_3 \big) = - \big(  m_{\psi} (\psi^{l,\alpha}) \nabla \mu_\psi^{l,\alpha} , \nabla \zeta_3 \big) \label{eq:second_approximate_chdf4}, & $\forall \zeta_3 \in X_l$, in $ [0,T]$, \\
    \big( \mu_\psi^{l,\alpha}, \zeta_4 \big) = \beta ( \nabla \psi^{l,\alpha}, \nabla \zeta_4 ) + \big(  F_{\psi,n,\kappa}'(\psi^{l,\alpha}) + \partial_\psi G_\kappa(\phi^{l,\alpha},\psi^{l,\alpha}) , \zeta_4 \big) \label{eq:second_approximate_chdf5}, & $\forall \zeta_4 \in X_l$, in $ [0,T]$,
\end{numcases}
where $\alpha>0$. 
The system is supplemented with the initial conditions
\begin{equation} \label{eq:second_approximate_chdfi}
    \bm{u}^{l,\alpha}|_{t=0} = \bm{u}_{0,l}, \quad \phi^{l,\alpha}|_{t=0} = \mathbb{P}^X_l (\phi_{0,k}), \quad \psi^{l,\alpha}|_{t=0} = \mathbb{P}^X_l (\psi_{0,k}).
\end{equation}
Here $\bm{u}_{0,l} \in \bm{V}_m$ is the unique solution in $\bm{V}_m$ to the following problem
\begin{align*}
    & (  \nu (\mathbb{P}^X_l (\phi_{0,k}),\mathbb{P}^X_l (\psi_{0,k})) \bm{u} + \eta (\mathbb{P}^X_l (\phi_{0,k}),\mathbb{P}^X_l (\psi_{0,k})) |\bm{u}|^{r-2} \bm{u} , \bm{\theta}_0 ) \nonumber \\
    & \quad = (  \mu_\phi^{l,\alpha}(0) \nabla \mathbb{P}^X_l (\phi_{0,k}) + \mu_\psi^{l,\alpha}(0) \nabla \mathbb{P}^X_l (\psi_{0,k}) , \bm{\theta}_0 ), \ \forall \bm{\theta}_0 \in \bm{V}_m, \\
    & \mu_\phi^{l,\alpha}(0) = -\Delta \mathbb{P}^X_l (\phi_{0,k}) + \sigma_2 \mathcal{N} \big( \mathbb{P}^X_l (\phi_{0,k}) - \overline{\mathbb{P}^X_l (\phi_{0,k})} \big) + \mathbb{P}^X_l F_{\phi,n,\kappa}'(\mathbb{P}^X_l (\phi_{0,k}))\\
    & \qquad \qquad + \mathbb{P}^X_l \partial_{\phi} G_\kappa(\mathbb{P}^X_l (\phi_{0,k}),\mathbb{P}^X_l (\psi_{0,k})), \\
    & \mu_\psi^{l,\alpha}(0) = - \beta \Delta \mathbb{P}^X_l (\psi_{0,k}) + \mathbb{P}^X_l F_{\psi,n,\kappa}'(\mathbb{P}^X_l (\psi_{0,k})) + \mathbb{P}^X_l \partial_{\psi} G_\kappa(\mathbb{P}^X_l (\phi_{0,k}),\mathbb{P}^X_l (\psi_{0,k})).
\end{align*}
Testing the second equation with $-\Delta \mu_\phi^{l,\alpha}(0)$, we have
\[ \big\| \nabla \mu_\phi^{l,\alpha}(0) \big\| \le \big\| \nabla \Delta \phi_{0,k} \big\| + C_n \big(  \big\| \nabla \phi_{0,k} \big\| + \big\| \nabla \psi_{0,k} \big\|  \big), \]
and, similarly for $\mu_\psi^{l,\alpha}(0)$, it holds
\[ \big\| \nabla \mu_\psi^{l,\alpha}(0) \big\| \le \beta \big\| \nabla \Delta \psi_{0,k} \big\| + C_n \big(  \big\| \nabla \phi_{0,k} \big\| + \big\| \nabla \psi_{0,k} \big\|  \big). \]
So, taking $\bm{\theta}_0 = \bm{u}_{0,l}$, we obtain
\begin{align*}
    & \nu_* \| \bm{u}_{0,l} \|^2 + \eta_* \| \bm{u}_{0,l} \|^r_{L^r(\Omega)} \\
    & \quad \le C \big\| \mathbb{P}^X_l (\phi_{0,k}) \big\|_{L^\infty(\Omega)}^2 \| \nabla \mu_\phi^{l,\alpha}(0) \|^2 + C \big\| \mathbb{P}^X_l (\psi_{0,k}) \big\|_{L^\infty(\Omega)}^2 \| \nabla \mu_\psi^{l,\alpha}(0) \|^2 \\
    & \quad \le C \big(  1 + \big\| \Delta \phi_{0,k} \big\|^2 + \beta \big\| \Delta \psi_{0,k} \big\|^2  \big) \big(  \big\| \nabla \Delta \phi_{0,k} \big\|^2 + \beta^2 \big\| \nabla \Delta \psi_{0,k} \big\|^2 + C_n \big(  \big\| \nabla \phi_{0,k} \big\|^2 + \big\| \nabla \psi_{0,k} \big\|^2  \big)  \big),
\end{align*}
which is uniformly bounded with respect to both $\alpha \in (0,1]$ and $l \in \mathbb{N}_+$.
The time-relaxation term in \eqref{eq:second_approximate_chdf1} allows us to write the Galerkin system as a system of ordinary differential equations to which the Cauchy--Lipschitz theorem applies.
It follows that the above problem admits a unique solution 
\[\big(  \bm{u}^{l,\alpha}, \phi^{l,\alpha}, \mu_\phi^{l,\alpha}, \psi^{l,\alpha}, \mu_\psi^{l,\alpha}  \big) \in C^1 \big(  [0,T];\bm{V}_m \times X_l^4  \big),\]
for some $T>0$.
Moreover, the choice of $\bm{u}_{0,l}$ gives $\partial_t \bm{u}^{l,\alpha} |_{t=0}=\bm 0$.
This eliminates the initial contribution $\tfrac{\alpha}{2}\| \partial_t\bm{u}^{l,\alpha}(0) \|^2$ in the time-integrated estimate \eqref{time_estimate_sp_with_alpha}, allowing a bound uniform in $\alpha\in(0,1]$ with the other approximation parameters fixed.
Then, we take $\bm{\theta} = \bm{u}^{l,\alpha}$, $\zeta_1 = \mu_\phi^{l,\alpha}$, $\zeta_2 = \partial_t \phi^{l,\alpha}$, $\zeta_3 = \mu_\psi^{l,\alpha}$ and $\zeta_4 = \partial_t \psi^{l,\alpha}$ to obtain
\begin{align}
    & \frac{\mathrm{d}}{\mathrm{d}t} \big(  E_{\mathrm{free}}^{n,\kappa} (\phi^{l,\alpha},\psi^{l,\alpha}) + \frac{\alpha}{2} \big\| \bm{u}^{l,\alpha} \big\|^2  \big) \nonumber \\
    & \qquad + \int_\Omega \nu (\phi^{l,\alpha},\psi^{l,\alpha}) \big| \bm{u}^{l,\alpha} \big|^2 + \eta (\phi^{l,\alpha},\psi^{l,\alpha}) \big| \bm{u}^{l,\alpha} \big|^r + m_\phi (\phi^{l,\alpha}) \big| \nabla \mu_\phi^{l,\alpha} \big|^2 + m_\psi (\psi^{l,\alpha}) \big| \nabla \mu_\psi^{l,\alpha} \big|^2 \ \mathrm{d} x \nonumber \\
    & \quad = - \sigma_1 |\Omega| \big( \overline{\phi^{l,\alpha}}-c \big) \overline{\mu_\phi^{l,\alpha}} \nonumber \\
    & \quad \le C_{\mathrm{source}} \sigma_1 \big|  \overline{\phi^{l,\alpha}}-c  \big| \big(  1 + \big\| \nabla \phi^{l,\alpha} \big\|^2  \big) + \frac{1}{2} \int_\Omega m_\phi (\phi^{l,\alpha}) \big| \nabla \mu_\phi^{l,\alpha} \big|^2 \ \mathrm{d} x, \label{basic_energy_estimate_approximate} \\
    & E_{\mathrm{free}}^{n,\kappa} (\phi^{l,\alpha},\psi^{l,\alpha}) \stackrel{\mathrm{def}}{=} \frac{1}{2} \big\| \nabla \phi^{l,\alpha} \big\|^2 + \frac{\beta}{2} \big\| \nabla \psi^{l,\alpha} \big\|^2 + \int_\Omega F_{\phi,n,\kappa}(\phi^{l,\alpha}) + F_{\psi,n,\kappa}(\psi^{l,\alpha}) + G_\kappa(\phi^{l,\alpha},\psi^{l,\alpha}) \ \mathrm{d} x \nonumber \\
    & \qquad \qquad \qquad \qquad + \frac{\sigma_2}{2} \big\|  \phi^{l,\alpha} - \overline{\phi^{l,\alpha}}  \big\|^2_*, \nonumber
\end{align}
where $C_{\mathrm{source}}$ is a positive constant only depending on $\overline{\phi_0}$,  $\Omega$, and the structural data (including $G$, $F_\phi$, $F_\psi$, $c$, $m_{\phi,*}$).
It follows from \eqref{potential_control_ohtakawasaki} that
\[
    E_{\mathrm{free}}^{n,\kappa} (\phi^{l,\alpha},\psi^{l,\alpha}) + 2 | \sigma_2 | C_P^2 | \Omega | + | \Omega | \| G \|_{C(\mathbb{R}^2)} \ge \frac{1}{2} \big\| \nabla \phi^{l,\alpha} \big\|^2 + \frac{\beta}{2} \big\| \nabla \psi^{l,\alpha} \big\|^2.
\]
Here we use the estimate
\begin{align*}
    \big| \overline{\mu_\phi^{l,\alpha}} \big| & = \big|  \overline{ F_{\phi,n,\kappa}'(\phi^{l,\alpha}) } + \overline{\partial_\phi G_\kappa(\phi^{l,\alpha},\psi^{l,\alpha})}  \big| \\
    & \le C \big( 1 + \big(  F_{\phi,n,\kappa}'(\phi^{l,\alpha}) , \phi^{l,\alpha} - \overline{\phi^{l,\alpha}} \big) \big) \\
    & \le C \big( 1+ \big( \mu_\phi^{l,\alpha} - \partial_\phi G_\kappa(\phi^{l,\alpha},\psi^{l,\alpha}) - \sigma_2 \mathcal{N}(\phi^{l,\alpha} - \overline{\phi^{l,\alpha}}) , \phi^{l,\alpha} - \overline{\phi^{l,\alpha}} \big)  \big) \\
    & \le C \big( 1 + (1 + \big\| \nabla \mu_\phi^{l,\alpha} \big\| + \big\| \nabla \phi^{l,\alpha} \big\|) \big\|  \nabla \phi^{l,\alpha}  \big\|  \big).
\end{align*}
Since it holds
\[ \sigma_1 \big|  \overline{\phi^{l,\alpha}}(t)-c  \big| = \sigma_1 e^{-\sigma_1 t} \big|  \overline{\phi^{l,\alpha}}(0)-c  \big|, \]
we have the uniform boundedness of $\bm{u}^{l,\alpha}$, $\big( \phi^{l,\alpha}, \psi^{l,\alpha}, \mu_\phi^{l,\alpha}, \mu_\psi^{l,\alpha} \big)$ in $\bm{V}_m$ and $X_l^4$, respectively, with respect to $T$.
So the solution $\big(  \bm{u}^{l,\alpha}, \phi^{l,\alpha}, \psi^{l,\alpha}, \mu_\phi^{l,\alpha}, \mu_\psi^{l,\alpha}  \big)$ can be extended to a global solution defined on $\mathbb{R}_+$.
In addition, $\big( \mu_\phi^{l,\alpha}, \mu_\psi^{l,\alpha} \big)$, $\bm{u}^{l,\alpha}$ are uniformly bounded in $L^2_{\mathrm{uloc}}\big( 0,\infty;X_l^2 \big)$ and $L^2\big( 0,\infty;\bm{V}_m \big) \cap L^r\big( 0,\infty;\bm{V}_m \big)$ with respect to $\alpha$.
Moreover, $\nabla \mu_\phi^{l,\alpha},\nabla \mu_\psi^{l,\alpha}$ are uniformly bounded in $L^2\big( 0,\infty;H \big)$ with respect to $\alpha$.
Therefore, by comparison in equations \eqref{eq:second_approximate_chdf2}--\eqref{eq:second_approximate_chdf5}, 
we obtain the uniform boundedness of $\partial_t \phi^{l,\alpha}, \partial_t \psi^{l,\alpha}$ in $L^2\big( 0,\infty;X_l \big)$ and then the uniform boundedness of $\partial_t \mu_\phi^{l,\alpha}, \partial_t \mu_\psi^{l,\alpha}$ in $L^2\big( 0,\infty;X_l \big)$ with respect to $\alpha$.
So by comparison in equation \eqref{eq:second_approximate_chdf1}, we have $\partial_t \bm{u}^{l,\alpha} \in \mathrm{Lip}\big( 0,T;\bm{V}_m \big)$ for all $T \ge0$.
Next, we take $\bm{\theta} =\partial_t \bm{u}^{l,\alpha}$ in \eqref{eq:second_approximate_chdf1} to obtain
\begin{align*}
    & \frac{\mathrm{d}}{\mathrm{d}t} \Big(  \int_\Omega \frac{1}{2} \nu (\phi^{l,\alpha},\psi^{l,\alpha}) \big| \bm{u}^{l,\alpha} \big|^2 + \frac{1}{r} \eta (\phi^{l,\alpha},\psi^{l,\alpha}) \big| \bm{u}^{l,\alpha} \big|^r \ \mathrm{d} x  \Big) + \alpha \big\| \partial_t \bm{u}^{l,\alpha} \big\|^2 \\
    & \quad = \big(  \mu_\phi^{l,\alpha} \nabla \phi^{l,\alpha} + \mu_\psi^{l,\alpha} \nabla \psi^{l,\alpha} , \partial_t \bm{u}^{l,\alpha} \big)  + \int_\Omega \frac{1}{2} \partial_t \nu (\phi^{l,\alpha},\psi^{l,\alpha}) \big| \bm{u}^{l,\alpha} \big|^2 + \frac{1}{r} \partial_t \eta (\phi^{l,\alpha},\psi^{l,\alpha}) \big| \bm{u}^{l,\alpha} \big|^r \ \mathrm{d} x .
\end{align*}
Observe now that
\begin{align*}
    & \big(  \mu_\phi^{l,\alpha} \nabla \phi^{l,\alpha} + \mu_\psi^{l,\alpha} \nabla \psi^{l,\alpha} , \partial_t \bm{u}^{l,\alpha} \big) \\
    & \quad = \frac{\mathrm{d}}{\mathrm{d}t} \Big(  \int_\Omega \mu_\phi^{l,\alpha} \bm{u}^{l,\alpha} \cdot \nabla \phi^{l,\alpha} + \mu_\psi^{l,\alpha} \bm{u}^{l,\alpha} \cdot \nabla \psi^{l,\alpha} \ \mathrm{d} x  \Big) \\
    & \qquad - \big(  \partial_t \mu_\phi^{l,\alpha} \nabla \phi^{l,\alpha} + \partial_t \mu_\psi^{l,\alpha} \nabla \psi^{l,\alpha} + \mu_\phi^{l,\alpha} \nabla \partial_t \phi^{l,\alpha} + \mu_\psi^{l,\alpha} \nabla \partial_t \psi^{l,\alpha} , \bm{u}^{l,\alpha} \big).
\end{align*}
Thus, we have
\begin{align*}
    & \frac{\mathrm{d}}{\mathrm{d}t} \Big(  \int_\Omega \frac{1}{2} \nu (\phi^{l,\alpha},\psi^{l,\alpha}) \big| \bm{u}^{l,\alpha} \big|^2 + \frac{1}{r} \eta (\phi^{l,\alpha},\psi^{l,\alpha}) \big| \bm{u}^{l,\alpha} \big|^r - \mu_\phi^{l,\alpha} \bm{u}^{l,\alpha} \cdot \nabla \phi^{l,\alpha} - \mu_\psi^{l,\alpha} \bm{u}^{l,\alpha} \cdot \nabla \psi^{l,\alpha} \ \mathrm{d} x  \Big) \\
    & \quad \le - \big(  \partial_t \mu_\phi^{l,\alpha} \nabla \phi^{l,\alpha} + \partial_t \mu_\psi^{l,\alpha} \nabla \psi^{l,\alpha} + \mu_\phi^{l,\alpha} \nabla \partial_t \phi^{l,\alpha} + \mu_\psi^{l,\alpha} \nabla \partial_t \psi^{l,\alpha} , \bm{u}^{l,\alpha} \big) \\
    & \qquad + \int_\Omega \frac{1}{2} \partial_t \nu (\phi^{l,\alpha},\psi^{l,\alpha}) \big| \bm{u}^{l,\alpha} \big|^2 + \frac{1}{r} \partial_t \eta (\phi^{l,\alpha},\psi^{l,\alpha}) \big| \bm{u}^{l,\alpha} \big|^r \ \mathrm{d} x \\
    & \quad \le C \big(  \big\| \partial_t \mu_\phi^{l,\alpha} \big\| \big\| \nabla \phi^{l,\alpha} \big\| + \big\| \partial_t \mu_\psi^{l,\alpha} \big\| \big\| \nabla \psi^{l,\alpha} \big\| + \big\| \nabla \mu_\phi^{l,\alpha}  \big\| \big\|  \partial_t \phi^{l,\alpha} \big\| + \big\| \nabla \mu_\psi^{l,\alpha} \big\| \big\|  \partial_t \psi^{l,\alpha} \big\|  \big) \big\| \bm{u}^{l,\alpha} \big\| \\
    & \qquad + C \big(  \big\|  \partial_t \phi^{l,\alpha} \big\| + \big\|  \partial_t \psi^{l,\alpha} \big\|  \big) \big(  \big\|  \bm{u}^{l,\alpha}  \big\|^2 + \big\|  \bm{u}^{l,\alpha}  \big\|^r  \big).
\end{align*}
Hence, the uniform $L^\infty(0,\infty;\bm{V}_m)$-boundedness of $\bm{u}^{l,\alpha}$ with respect to $\alpha$ follows from the above inequality, the Gronwall inequality and the uniform Gronwall inequality.
As a result, we also obtain the uniform $L^\infty(0,\infty;X_l)$-boundedness of $\partial_t \phi^{l,\alpha}$, $\partial_t \psi^{l,\alpha}$, $\partial_t \mu_\phi^{l,\alpha}$ and $\partial_t \mu_\psi^{l,\alpha}$.
Then, by Lemma \ref{chain_rule}, if we differentiate equation \eqref{eq:second_approximate_chdf1} with respect to $t$ and take $\bm{\theta}=\partial_t \bm{u}^{l,\alpha}$ in the resulting equation, then we obtain
\begin{align*}
    & \frac{\alpha}{2} \frac{\mathrm{d}}{\mathrm{d}t} \big\| \partial_t \bm{u}^{l,\alpha} \big\|^2 + \int_\Omega \nu (\phi^{l,\alpha},\psi^{l,\alpha}) \big| \partial_t \bm{u}^{l,\alpha} \big|^2 + \eta (\phi^{l,\alpha},\psi^{l,\alpha}) \big| \bm{u}^{l,\alpha} \big|^{r-2} \big| \partial_t \bm{u}^{l,\alpha} \big|^2 \ \mathrm{d} x \\
    & \qquad + \int_\Omega (r-2) \eta (\phi^{l,\alpha},\psi^{l,\alpha}) \big| \bm{u}^{l,\alpha} \big|^{r-2} \big|  \bm{e}_{\mathbb{R}^d}(\bm{u}^{l,\alpha}) \cdot \partial_t \bm{u}^{l,\alpha} \big|^2 \ \mathrm{d} x \\
    & \quad = \big(  \partial_t \mu_\phi^{l,\alpha} \nabla \phi^{l,\alpha} + \partial_t \mu_\psi^{l,\alpha} \nabla \psi^{l,\alpha} + \mu_\phi^{l,\alpha} \nabla \partial_t \phi^{l,\alpha} + \mu_\psi^{l,\alpha} \nabla \partial_t \psi^{l,\alpha} , \partial_t \bm{u}^{l,\alpha} \big) \\
    & \qquad - \big(  \partial_t \nu (\phi^{l,\alpha},\psi^{l,\alpha}) \bm{u}^{l,\alpha} + \partial_t \eta (\phi^{l,\alpha},\psi^{l,\alpha}) \big| \bm{u}^{l,\alpha} \big|^{r-2} \bm{u}^{l,\alpha} , \partial_t \bm{u}^{l,\alpha} \big).
\end{align*}
It also follows from Lemma \ref{chain_rule} that
\begin{align*}
    & \big| \bm{u}^{l,\alpha} \big|^{r-2} \bm{u}^{l,\alpha} \cdot \partial_t \bm{u}^{l,\alpha} = \frac{2}{r} \big| \bm{u}^{l,\alpha} \big|^\frac{r}{2} \partial_t \big| \bm{u}^{l,\alpha} \big|^\frac{r}{2}, \\
    %%%%%%%%%%%%%%%%%%%%%%%%%%%%%%%%
    & \big| \bm{u}^{l,\alpha} \big|^{r-2} \big|  \bm{e}_{\mathbb{R}^d}(\bm{u}^{l,\alpha}) \cdot \partial_t \bm{u}^{l,\alpha} \big|^2 = \frac{4}{r^2} \Big|  \partial_t \big| \bm{u}^{l,\alpha} \big|^\frac{r}{2}  \Big|^2.
\end{align*}
In addition, note that
\begin{align*}
    & \big| \bm{u}^{l,\alpha} \big|^{r-2} \big| \partial_t \bm{u}^{l,\alpha} \big|^2 \ge \big| \bm{u}^{l,\alpha} \big|^{r-2} \big|  \bm{e}_{\mathbb{R}^d}(\bm{u}^{l,\alpha}) \cdot \partial_t \bm{u}^{l,\alpha} \big|^2.
\end{align*}
Therefore, we get
\begin{align}
    & \frac{\alpha}{2} \frac{\mathrm{d}}{\mathrm{d}t} \big\| \partial_t \bm{u}^{l,\alpha} \big\|^2 + \nu_* \big\| \partial_t \bm{u}^{l,\alpha} \big\|^2 + \frac{4 (r-1) \eta_*}{r^2} \Big\|  \partial_t \big| \bm{u}^{l,\alpha} \big|^\frac{r}{2}  \Big\|^2 \nonumber \\
    & \quad \le C \big(  \big\| \partial_t \mu_\phi^{l,\alpha} \big\| \big\| \phi^{l,\alpha} \big\| + \big\| \partial_t \mu_\psi^{l,\alpha} \big\| \big\| \psi^{l,\alpha} \big\| + \big\| \mu_\phi^{l,\alpha}  \big\| \big\|  \partial_t \phi^{l,\alpha} \big\| + \big\| \mu_\psi^{l,\alpha} \big\| \big\|  \partial_t \psi^{l,\alpha} \big\|  \big) \big\| \partial_t \bm{u}^{l,\alpha} \big\| \nonumber \\
    & \qquad + C \big(  \big\|  \partial_t \phi^{l,\alpha} \big\| + \big\|  \partial_t \psi^{l,\alpha} \big\|  \big) \Big(  \Big\| \big| \bm{u}^{l,\alpha} \big|^\frac{r}{2} \Big\| \Big\| \partial_t \big| \bm{u}^{l,\alpha} \big|^\frac{r}{2} \Big\| + \big\| \bm{u}^{l,\alpha} \big\| \big\| \partial_t \bm{u}^{l,\alpha} \big\|  \Big) \nonumber \\
    & \quad \le \frac{1}{2} \nu_* \big\| \partial_t \bm{u}^{l,\alpha} \big\|^2 + \frac{2 (r-1) \eta_*}{r^2} \Big\|  \partial_t \big| \bm{u}^{l,\alpha} \big|^\frac{r}{2}  \Big\|^2 \nonumber \\
    & \qquad + C \big(  \big\| \partial_t \mu_\phi^{l,\alpha} \big\| \big\| \phi^{l,\alpha} \big\| + \big\| \partial_t \mu_\psi^{l,\alpha} \big\| \big\| \psi^{l,\alpha} \big\| + \big\| \mu_\phi^{l,\alpha}  \big\| \big\|  \partial_t \phi^{l,\alpha} \big\| + \big\| \mu_\psi^{l,\alpha} \big\| \big\|  \partial_t \psi^{l,\alpha} \big\|  \big)^2 \nonumber \\
    & \qquad + C \big(  \big\|  \partial_t \phi^{l,\alpha} \big\|^2 + \big\|  \partial_t \psi^{l,\alpha} \big\|^2  \big) \big(  \big\| \bm{u}^{l,\alpha} \big\|^2 + \big\| \bm{u}^{l,\alpha} \big\|^r  \big). \label{time_estimate_sp_with_alpha}
\end{align}
Recall that our choice of $\bm{u}^{l,\alpha}(0)$ implies that $\partial_t \bm{u}^{l,\alpha} |_{t=0} = \bm 0$.
Hence, we obtain the uniform $H^1(0,\infty;\bm{V}_m)$-boundedness of $\bm{u}^{l,\alpha}$.
Now, letting $\alpha \to 0$, with a standard compactness argument, we find a solution $\big( \bm{u}^l,\phi^l,\psi^l,\mu_\phi^l,\mu_\psi^l \big)$ to problem \eqref{eq:second_approximate_chdf1}--\eqref{eq:second_approximate_chdfi} with $\alpha=0$, which possesses the following regularities:
\begin{align*}
    \phi^l - \overline{\phi^l},\psi^l - \overline{\psi^l} & \in W^{1,\infty} (0,\infty;X_l) \cap H^1_{\mathrm{uloc}} (0,\infty;X_l),\\
    \mu_\phi^l,\mu_\psi^l & \in W^{1,\infty} (0,\infty;X_l) \cap L^2_{\mathrm{uloc}} (0,\infty;X_l), \\
    \nabla \mu_\phi^l, \nabla \mu_\psi^l, \partial_t \mu_\phi^l,\partial_t \mu_\psi^l & \in L^2 (0,\infty;H),\\
    \bm{u}^l & \in H^1(0,\infty;\bm{V}_m).
\end{align*}

\noindent \textbf{Step 3: passing to the limit as $l \to \infty$.}
We first recall the basic energy estimate:
\begin{align*}
    & \frac{\mathrm{d}}{\mathrm{d}t} \Big(  \frac{1}{2} \big\| \nabla \phi^l \big\|^2 + \frac{\beta}{2} \big\| \nabla \psi^l \big\|^2 + \int_\Omega F_{\phi,n,\kappa}\big( \phi^l \big) + F_{\psi,n,\kappa}\big( \psi^l \big) + G_\kappa(\phi^l,\psi^l) \ \mathrm{d} x + \frac{\sigma_2}{2} \big\|  \phi^l - \overline{\phi^l}  \big\|^2_*  \Big) \\ 
    & \qquad + \frac{1}{2} \int_\Omega \nu (\phi^l,\psi^l) \big| \bm{u}^l \big|^2 + \eta (\phi^l,\psi^l) \big| \bm{u}^l \big|^r + m_\phi \big( \phi^l \big) \big| \nabla \mu_\phi^l \big|^2 + m_\psi \big( \psi^l \big) \big| \nabla \mu_\psi^l \big|^2 \ \mathrm{d} x \\
    & \quad \le C_{\mathrm{source}} \sigma_1 \big|  \overline{\phi^l}-c  \big| \big(  1 + \big\| \nabla \phi^l \big\|^2  \big).
\end{align*}
We also recall that $\phi_{0,k},\psi_{0,k} \in H^3(\Omega)$ and 
\[\mathbb{P}^X_l \phi_{0,k} \to \phi_{0,k}, \ \mathbb{P}^X_l \psi_{0,k} \to \psi_{0,k} \quad \mathrm{in} \ H^3(\Omega). \]
As a consequence, we deduce the following uniform boundedness with respect to $l$ in the corresponding spaces:
\begin{align*}
    & \phi^l,\psi^l \in L^\infty(0,\infty;V),\\
    & \mu_\phi^l,\mu_\psi^l \in L^2_{\mathrm{uloc}}(0,\infty;V) \ \text{with} \ \nabla \mu_\phi^l, \nabla \mu_\psi^l \in L^2(0,\infty;H) \\
    & \bm{u}^l \in L^2(0,\infty;\bm{V}_m) \cap L^r(0,\infty;\bm{V}_m).
\end{align*}
Let us now obtain some higher-order bounds.
Differentiating equation \eqref{eq:second_approximate_chdf1} with respect to time and taking $\bm{\theta}=\bm{u}^l$, we have
\begin{align}
    & \frac{\mathrm{d}}{\mathrm{d}t} \Big(  \int_\Omega \frac{1}{2} \nu \big( \phi^l,\psi^l \big) \big| \bm{u}^l \big|^2 + \frac{r-1}{r} \eta \big( \phi^l,\psi^l \big) \big| \bm{u}^l \big|^r \ \mathrm{d} x  \Big) \nonumber \\
    & \quad = - \int_\Omega \frac{1}{2} \partial_t \nu \big( \phi^l,\psi^l \big) \big| \bm{u}^l \big|^2 + \frac{1}{r} \partial_t \eta \big( \phi^l,\psi^l \big) \big| \bm{u}^l \big|^r \ \mathrm{d} x \nonumber \\
    & \qquad + \int_\Omega \big( \partial_t \mu_\phi^l \nabla \phi^l + \partial_t \mu_\psi^l \nabla \psi^l \big) \cdot \bm{u}^l \ \mathrm{d} x - \int_\Omega \big(  \partial_t \phi^l \nabla \mu_\phi^l + \partial_t \psi^l \nabla \mu_\psi^l  \big) \cdot \bm{u}^l \ \mathrm{d} x. \label{higher_order_estimate_sp}
\end{align}
Taking $\zeta_1 = \partial_t \mu_\phi^l$ (resp. $\zeta_3 = \partial_t \mu_\psi^l$) in \eqref{eq:second_approximate_chdf2} (resp. \eqref{eq:second_approximate_chdf4}) and differentiating \eqref{eq:second_approximate_chdf3} (resp. \eqref{eq:second_approximate_chdf5}) with respect to time, then taking $\zeta_2=\partial_t \phi^l$ (resp. $\zeta_4=\partial_t \psi^l$), we obtain
\begin{align*}
    & \frac{1}{2} \frac{\mathrm{d}}{\mathrm{d}t} \Big(  \int_\Omega m_\phi\big( \phi^l \big) \big| \nabla \mu_\phi^l \big|^2 + m_\psi\big( \psi^l \big) \big| \nabla \mu_\psi^l \big|^2 \ \mathrm{d} x  \Big) + \big\| \nabla \partial_t \phi^l \big\|^2 + \beta \big\| \nabla \partial_t \psi^l \big\|^2 \\
    & \qquad + \sigma_2 \big\| \partial_t \phi^l - \overline{\partial_t \phi^l} \big\|^2_* + \int_\Omega F_{\phi,n,\kappa}''\big( \phi^l \big) \big| \partial_t \phi^l \big|^2 + F_{\psi,n,\kappa}''\big( \psi^l \big) \big| \partial_t \psi^l \big|^2 \ \mathrm{d} x \\
    & \quad = \frac{1}{2} \int_\Omega m_\phi'\big( \phi^l \big) \partial_t \phi^l \big| \nabla \mu_\phi^l \big|^2 + m_\psi'\big( \psi^l \big) \partial_t \psi^l \big| \nabla \mu_\psi^l \big|^2 \ \mathrm{d} x \\
    & \qquad - \int_\Omega \bm{u}^l \cdot \big(  \partial_t \mu_\phi^l \nabla \phi^l + \partial_t \mu_\psi^l \nabla \psi^l  \big) \ \mathrm{d} x - \sigma_1 \big(  \overline{\phi^l} - c  \big) |\Omega| \overline{ \partial_t \mu_\phi^l } \\
    & \qquad - \int_\Omega \partial_\phi^2 G_\kappa (\phi^l,\psi^l) \big| \partial_t \phi^l \big|^2 + \partial_\psi^2 G_\kappa (\phi^l,\psi^l) \big| \partial_t \psi^l \big|^2 + 2 \partial_{\phi,\psi}^2 G_\kappa (\phi^l,\psi^l) \partial_t \phi^l \partial_t \psi^l \ \mathrm{d} x.
\end{align*}
Adding the above two equations together, we get
\begin{align}
    & \frac{\mathrm{d}}{\mathrm{d}t} \mathcal{E} \big( \bm{u}^l,\phi^l,\psi^l,\mu_\phi^l,\mu_\psi^l \big) \nonumber + \big\| \nabla \partial_t \phi^l \big\|^2 + \beta \big\| \nabla \partial_t \psi^l \big\|^2 \nonumber \\
    & \quad \le \frac{1}{2} \int_\Omega m_\phi'\big( \phi^l \big) \partial_t \phi^l \big| \nabla \mu_\phi^l \big|^2 + m_\psi'\big( \psi^l \big) \partial_t \psi^l \big| \nabla \mu_\psi^l \big|^2 \ \mathrm{d} x \nonumber \\
    & \qquad - \int_\Omega \frac{1}{2} \partial_t \nu \big( \phi^l,\psi^l \big) \big| \bm{u}^l \big|^2 + \frac{1}{r} \partial_t \eta \big( \phi^l,\psi^l \big) \big| \bm{u}^l \big|^r \ \mathrm{d} x \nonumber \\
    & \qquad - \int_\Omega \big(  \partial_t \phi^l \nabla \mu_\phi^l + \partial_t \psi^l \nabla \mu_\psi^l  \big) \cdot \bm{u}^l \ \mathrm{d} x - \sigma_1^2 |\Omega| \big(  \overline{\phi^l} - c  \big)\overline{ \mu_\phi^l } \nonumber \\
    & \qquad - \int_\Omega \partial_\phi^2 G_\kappa (\phi^l,\psi^l) \big| \partial_t \phi^l \big|^2 + \partial_\psi^2 G_\kappa (\phi^l,\psi^l) \big| \partial_t \psi^l \big|^2 + 2 \partial_{\phi,\psi}^2 G_\kappa (\phi^l,\psi^l) \partial_t \phi^l \partial_t \psi^l \ \mathrm{d} x \nonumber \\
    & \qquad + | \sigma_2 | \big\| \partial_t \phi^l - \overline{\partial_t \phi^l} \big\|^2_* \nonumber \\
    & \quad \le C \big(  \big\| \partial_t \phi^l \big\| \big\| \nabla \mu_\phi^l \big\|^2_{L^4(\Omega)} + \big\| \partial_t \psi^l \big\| \big\| \nabla \mu_\psi^l \big\|^2_{L^4(\Omega)}  \big) \nonumber \\
    & \qquad + C \big(  \big\| \partial_t \phi^l \big\| + \big\| \partial_t \psi^l \big\|  \big) \big(  \big\| \bm{u}^l \big\|^2_{L^4(\Omega)} + \big\| \bm{u}^l \big\|^r_{L^{2r}(\Omega)}  \big) \nonumber \\
    & \qquad + C \big\| \bm{u}^l \big\|_{L^3(\Omega)} \big(  \big\| \nabla \partial_t \phi^l \big\| \big\| \nabla \mu_\phi^l \big\| + \big\| \nabla \partial_t \psi^l \big\| \big\| \nabla \mu_\psi^l \big\|  \big) \nonumber \\
    & \qquad + C ( 1 + (1 + \big\| \nabla \mu_\phi^l \big\|) \| \nabla \phi^l \| ) + C \big(  \big| \overline{\partial_t \phi^l} \big|^2 + \big\| \partial_t \phi^l - \overline{\partial_t \phi^l} \big\|^2 + \big\| \partial_t \psi^l \big\|^2  \big), \label{higher_order_estimate_approximate}
\end{align}
where 
\begin{align}
    \mathcal{E}(\bm{u},\phi,\psi,\mu_\phi,\mu_\psi) & \stackrel{\rm{def}}{=} \int_\Omega \frac{1}{2} m_\phi\big( \phi \big) \big| \nabla \mu_\phi \big|^2 + \frac{1}{2} m_\psi\big( \psi \big) \big| \nabla \mu_\psi \big|^2 \ \mathrm{d} x + \sigma_1 |\Omega| \big( \overline{\phi}-c \big) \overline{ \mu_\phi } \nonumber \\
    & \quad + \int_\Omega \frac{1}{2} \nu \big( \phi,\psi \big) \big| \bm{u} \big|^2 + \frac{r-1}{r} \eta \big( \phi,\psi \big) \big| \bm{u} \big|^r \ \mathrm{d} x. \label{higher_order_energy}
\end{align}
It holds
\begin{align}
    \mathcal{E}(\bm{u},\phi,\psi,\mu_\phi,\mu_\psi) & \ge \int_\Omega \frac{1}{4} m_\phi\big( \phi \big) \big| \nabla \mu_\phi \big|^2 + \frac{1}{2} m_\psi\big( \psi \big) \big| \nabla \mu_\psi \big|^2 + \frac{1}{2} \nu \big( \phi,\psi \big) \big| \bm{u} \big|^2 \ \mathrm{d} x\nonumber \\
    & \quad + \frac{r-1}{r} \int_\Omega\eta \big( \phi,\psi \big) \big| \bm{u} \big|^r \ \mathrm{d} x - C_{\mathrm{source}} \sigma_1 \big| \overline{\phi}-c \big| (1 + \| \nabla \phi \|^2) , \label{higher_order_energy_lower_bound} \\
    \mathcal{E}(\bm{u},\phi,\psi,\mu_\phi,\mu_\psi) & \le \int_\Omega \frac{3}{4} m_\phi\big( \phi \big) \big| \nabla \mu_\phi \big|^2 + \frac{1}{2} m_\psi\big( \psi \big) \big| \nabla \mu_\psi \big|^2 + \frac{1}{2} \nu \big( \phi,\psi \big) \big| \bm{u} \big|^2\ \mathrm{d} x
    \nonumber\\
    & \quad + \frac{r-1}{r} \int_\Omega\eta \big( \phi,\psi \big) \big| \bm{u} \big|^r \ \mathrm{d} x + C_{\mathrm{source}} \sigma_1 \big| \overline{\phi}-c \big| (1 + \| \nabla \phi \|^2). \label{higher_order_energy_upper_bound}
\end{align}
Next, we estimate the right-hand side of the above inequality \eqref{higher_order_estimate_approximate} term by term.
Firstly, we have
\begin{align*}
    \big\| \partial_t \phi^l - \overline{\partial_t \phi^l} \big\| & \le C \big\| \nabla \partial_t \phi^l \big\|^\frac{1}{2} \big\| \partial_t \phi^l - \overline{\partial_t \phi^l}  \big\|_*^\frac{1}{2} \\
    & \le C \big\| \nabla \partial_t \phi^l \big\|^\frac{1}{2} \big(  \big\| \nabla \mu_\phi^l \big\|^\frac{1}{2} + \big\| \bm{u}^l \big\|_{L^3(\Omega)}^\frac{1}{2} \big\| \nabla \phi^l \big\|^\frac{1}{2}  \big) \\
    & \le C_m \big\| \nabla \partial_t \phi^l \big\|^\frac{1}{2} \big(  \big\| \nabla \mu_\phi^l \big\|^\frac{1}{2} + \big\| \bm{u}^l \big\|^\frac{1}{2}  \big) \\
    & \le \sqrt{\epsilon} \big\| \nabla \partial_t \phi^l \big\| + C_\epsilon \big( 1 + C_1 + \mathcal{E} \big( \bm{u}^l,\phi^l,\psi^l,\mu_\phi^l,\mu_\psi^l \big) \big)^\frac{1}{2}, \\
    \big\| \partial_t \psi^l \big\| & \le C \big\| \nabla \partial_t \psi^l \big\|^\frac{1}{2} \big\| \partial_t \psi^l \big\|_*^\frac{1}{2} \\
    & \le C_m \big\| \nabla \partial_t \psi^l \big\|^\frac{1}{2} \big(  \big\| \nabla \mu_\psi^l \big\|^\frac{1}{2} + \big\| \bm{u}^l \big\|^\frac{1}{2}  \big) \\
    & \le \sqrt{\epsilon} \big\| \nabla \partial_t \psi^l \big\| + C_\epsilon \big( 1 + C_1 + \mathcal{E} \big( \bm{u}^l,\phi^l,\psi^l,\mu_\phi^l,\mu_\psi^l \big) \big)^\frac{1}{2},
\end{align*}
where 
\begin{equation} \label{C_1_define}
    C_1 \stackrel{\mathrm{def}}{=} 4 \sigma_1 C_{\mathrm{source}} e^{4 C_{\mathrm{source}}} (E_{\mathrm{free}}(\phi_0,\psi_0) + | \Omega |\big( \| G \|_{C(\mathbb{R}^2)} + 2 | \sigma_2 | C_P^2 \big) + 1).
\end{equation}
Here we choose $k$, $l$, $n$ and $\kappa^{-1}$ sufficiently large, in this order, to ensure
\[
    E_{\mathrm{free}}^{n,\kappa} (\phi^l(0),\psi^l(0)) \le E_{\mathrm{free}}(\phi_0,\psi_0) + \frac{1}{4}.
\]
Then, taking $\zeta_1= \Delta \mu_\phi^l$ in \eqref{eq:second_approximate_chdf2}, we have
\begin{align*}
    & \int_\Omega m_\phi \big( \phi^l \big) \big| \Delta \mu_\phi^l \big|^2 \ \mathrm{d} x \\
    & \quad = \big(  \partial_t \phi^l - \overline{\partial_t \phi^l} + \bm{u}^l \cdot \nabla \phi^l - \nabla m_\phi\big( \phi^l \big) \cdot \nabla \mu_\phi^l , \Delta \mu_\phi^l \big)\\
    & \quad \le C \Big(  \big\| \partial_t \phi^l - \overline{\partial_t \phi^l} \big\|^\frac{1}{2}_* \big\| \nabla \partial_t \phi^l \big\|^\frac{1}{2} + \big\| \bm{u}^l \big\|_{L^3(\Omega)} \big\| \nabla \phi^l \big\|_{L^6(\Omega)} + \big\| \nabla \phi^l \big\|_{L^6(\Omega)} \big\| \nabla \mu_\phi^l \big\|_{L^3(\Omega)}  \Big) \big\| \Delta \mu_\phi^l \big\| \\
    & \quad \le C \big\| \nabla \partial_t \phi^l \big\|^\frac{1}{2} \big(  \big\| \bm{u}^l \big\|_{L^3(\Omega)} \big\| \nabla \phi^l \big\| + \big\| \nabla \mu_\phi^l \big\|  \big)^\frac{1}{2} \big\| \Delta \mu_\phi^l \big\| \\
    & \qquad + C_m \big\| \bm{u}^l \big\| \big(  \big\| \nabla \phi^l \big\|^\frac{1}{2} \big\| \nabla \mu_\phi^l \big\|^\frac{1}{2} + \big\| \nabla \phi^l \big\| + \big\| \nabla \psi^l \big\|  \big) \big\| \Delta \mu_\phi^l \big\| \\
    & \qquad + C \big(  \big\| \nabla \phi^l \big\|^\frac{1}{2} \big\| \nabla \mu_\phi^l \big\|^\frac{1}{2} + \big\| \nabla \phi^l \big\| + \big\| \nabla \psi^l \big\|  \big) \big\| \nabla \mu_\phi^l \big\|^\frac{1}{2} \big\| \Delta \mu_\phi^l \big\|^\frac{3}{2} \\
    & \quad \le \frac{1}{2} \int_\Omega m_\phi \big( \phi^l \big) \big| \Delta \mu_\phi^l \big|^2 \ \mathrm{d} x + \frac{1}{2} \big\| \nabla \partial_t \phi^l \big\|^\frac{3}{2} + C \big(  1 + \big\| \bm{u}^l \big\|^4 + \big\| \nabla \mu_\phi^l \big\|^4  \big),
\end{align*}
where we used the following estimates for $\phi^l$ by comparison in \eqref{eq:second_approximate_chdf3}
\begin{align*}
    & \big\| \Delta \phi^l \big\|^2 \le \big\| \nabla \mu_\phi^l \big\| \big\| \nabla \phi^l \big\| + C \big(  \big\| \nabla \phi^l \big\|^2 + \big\| \nabla \psi^l \big\|^2  \big), \\
    & \big\| \nabla \Delta \phi^l \big\| \le C_n \big(  \big\| \nabla \mu_\phi^l \big\| + \big\| \nabla \phi^l \big\| + \big\| \nabla \psi^l \big\|  \big) .
\end{align*}
This implies
\[ \int_\Omega m_\phi \big( \phi^l \big) \big| \Delta \mu_\phi^l \big|^2 \ \mathrm{d} x \le \big\| \nabla \partial_t \phi^l \big\|^\frac{3}{2} + C \big(  1 + \big\| \bm{u}^l \big\|^4 + \big\| \nabla \mu_\phi^l \big\|^4  \big),\]
and a similar argument yields
\[ \int_\Omega m_\psi \big( \psi^l \big) \big| \Delta \mu_\psi^l \big|^2 \ \mathrm{d} x \le \big\| \nabla \partial_t \psi^l \big\|^\frac{3}{2} + C \big(  1 + \big\| \bm{u}^l \big\|^4 + \big\| \nabla \mu_\psi^l \big\|^4  \big).\]
So it holds
\begin{align*}
    & \big\| \partial_t \phi^l \big\| \big\| \nabla \mu_\phi^l \big\|^2_{L^4(\Omega)} + \big\| \partial_t \psi^l \big\| \big\| \nabla \mu_\psi^l \big\|^2_{L^4(\Omega)} \\
    & \quad \le C_m \big(  \sigma_1 e^{-\sigma_1 t} + \big\| \nabla \partial_t \phi^l \big\|^\frac{1}{2} \big(  \big\| \nabla \mu_\phi^l \big\|^\frac{1}{2} + \big\| \bm{u}^l \big\|^\frac{1}{2}  \big)  \big) \big\| \nabla \mu_\phi^l \big\|^\frac{1}{2} \big\| \Delta \mu_\phi^l \big\|^\frac{3}{2} \\
    & \qquad + C_m \big\| \nabla \partial_t \psi^l \big\|^\frac{1}{2} \big(  \big\| \nabla \mu_\psi^l \big\|^\frac{1}{2} + \big\| \bm{u}^l \big\|^\frac{1}{2}  \big) \big\| \nabla \mu_\psi^l \big\|^\frac{1}{2} \big\| \Delta \mu_\psi^l \big\|^\frac{3}{2} \\
    & \quad \le C \big(  \big\| \nabla \mu_\phi^l \big\|^\frac{1}{2} \big\| \nabla \partial_t \phi^l \big\|^\frac{9}{8} + \big\| \nabla \mu_\phi^l \big\|^\frac{1}{2} + \big\| \bm{u}^l \big\|^\frac{7}{2} + \big\| \nabla \mu_\phi^l \big\|^\frac{7}{2}  \big) \\
    & \qquad + C_m \big\| \nabla \partial_t \phi^l \big\|^\frac{13}{8} \big(  1 + \big\| \bm{u}^l \big\| + \big\| \nabla \mu_\phi^l \big\|  \big) + C_m \big\| \nabla \partial_t \phi^l \big\|^\frac{1}{2} \big(  1 + \big\| \bm{u}^l \big\|^4 + \big\| \nabla \mu_\phi^l \big\|^4  \big) \\
    & \qquad + C_m \big\| \nabla \partial_t \psi^l \big\|^\frac{13}{8} \big(  1 + \big\| \bm{u}^l \big\| + \big\| \nabla \mu_\psi^l \big\|  \big) + C_m \big\| \nabla \partial_t \psi^l \big\|^\frac{1}{2} \big(  1 + \big\| \bm{u}^l \big\|^4 + \big\| \nabla \mu_\psi^l \big\|^4  \big) \\
    & \quad \le \epsilon \big(  \big\| \nabla \partial_t \phi^l \big\|^2 + \beta \big\| \nabla \partial_t \psi^l \big\|^2  \big) + C_\epsilon \big(  1 + C_1 + \mathcal{E} \big( \bm{u}^l,\phi^l,\psi^l,\mu_\phi^l,\mu_\psi^l \big)  \big)^4.
\end{align*}
For the second term on the right-hand side of \eqref{higher_order_estimate_approximate}, we have
\begin{align*}
    & \big(  \big\| \partial_t \phi^l \big\| + \big\| \partial_t \psi^l \big\|  \big) \big(  \big\| \bm{u}^l \big\|^2_{L^4(\Omega)} + \big\| \bm{u}^l \big\|^r_{L^{2r}(\Omega)}  \big) \\
    & \quad \le C_m \big(  \big\| \bm{u}^l \big\|^2 + \big\| \bm{u}^l \big\|^r  \big) \Big(  1 + \big\| \nabla \partial_t \phi^l \big\|^\frac{1}{2} \big(  \big\| \nabla \mu_\phi^l \big\|^\frac{1}{2} + \big\| \bm{u}^l \big\|^\frac{1}{2}  \big) + \big\| \nabla \partial_t \psi^l \big\|^\frac{1}{2} \big(  \big\| \nabla \mu_\psi^l \big\|^\frac{1}{2} + \big\| \bm{u}^l \big\|^\frac{1}{2}  \big)  \Big) \\
    & \quad \le \epsilon \big(  \big\| \nabla \partial_t \phi^l \big\|^2 + \beta \big\| \nabla \partial_t \psi^l \big\|^2  \big) + C_\epsilon \big(  1 + C_1 + \mathcal{E} \big( \bm{u}^l,\phi^l,\psi^l,\mu_\phi^l,\mu_\psi^l \big)  \big)^\frac{5}{3}.
\end{align*}
Then, for the third term in \eqref{higher_order_estimate_approximate}, it holds
\begin{align*}
    & \big\| \bm{u}^l \big\|_{L^3(\Omega)} \big(  \big\| \nabla \partial_t \phi^l \big\| \big\| \nabla \mu_\phi^l \big\| + \big\| \nabla \partial_t \psi^l \big\| \big\| \nabla \mu_\psi^l \big\|  \big) \\
    & \quad \le \epsilon \big(  \big\| \nabla \partial_t \phi^l \big\|^2 + \beta \big\| \nabla \partial_t \psi^l \big\|^2  \big) + C_\epsilon C_m \big(  1 + C_1 + \mathcal{E} \big( \bm{u}^l,\phi^l,\psi^l,\mu_\phi^l,\mu_\psi^l \big)  \big)^2.
\end{align*}
Adding the above estimates together, we obtain
    \[ \frac{\mathrm{d}}{\mathrm{d}t} \mathcal{E} \big( \bm{u}^l,\phi^l,\psi^l,\mu_\phi^l,\mu_\psi^l \big)  + \frac{1}{2} \big\| \nabla \partial_t \phi^l \big\|^2 + \frac{\beta}{2} \big\| \nabla \partial_t \psi^l \big\|^2 
     \le C \big( 1 + C_1 + \mathcal{E} \big( \bm{u}^l,\phi^l,\psi^l,\mu_\phi^l,\mu_\psi^l \big) \big)^4.
     \]
Observe now that
\[ \mathcal{E}(\bm{u},\phi,\psi,\mu_\phi,\mu_\psi) \le C_1 + m_\phi^* \big\| \nabla \mu_\phi \big\|^2 + m_\psi^* \big\| \nabla \mu_\psi \big\|^2 + \nu^* \big\| \bm{u} \big\|^2 + \eta^* \big\| \bm{u} \big\|^r_{L^r(\Omega)}. \]
Thus, in particular, we have
\begin{equation*}
     \mathcal{E}(\bm{u}_{0,l},\mathbb{P}^X_l \phi_{0,k},\mathbb{P}^X_l \psi_{0,k}, \mu_\phi^l(0),\mu_\psi^l(0)) \le K_{m,n,k}, 
\end{equation*}
where
\begin{align*}
    K_{m,n,k}  &\stackrel{\rm{def}}{=} C_2 \big(  1 + \big\| \Delta \phi_{0,k} \big\|^2 + \beta \big\| \Delta \psi_{0,k}  \big\|^2  \big) \big(  \big\| \nabla \Delta \phi_{0,k} \big\|^2 + \beta^2 \big\| \nabla \Delta \psi_{0,k} \big\|^2 \\
    & \quad + C_n \big(  \big\| \nabla \phi_{0,k} \big\|^2 + \big\| \nabla \psi_{0,k} \big\|^2  \big)  \big)  + C_3.
\end{align*}
Setting $T_e' = \frac{1}{3 C} \big(  1 + C_1 + K_{m,n,k}  \big)^{-3}$, we have, for any $t \in [0,T_e')$, that
\begin{align*}
    & 1 + C_1 + \mathcal{E}(\bm{u}^l(t),\phi^l(t),\psi^l(t),\mu_\phi^l(t),\mu_\psi^l(t)) + \int_0^t \frac{1}{2} \big\| \nabla \partial_t \phi^l \big\|^2 + \frac{\beta}{2} \big\| \nabla \partial_t \psi^l \big\|^2 \ \mathrm{d} \tau \\
    & \quad \le \Big(  \big( 1 + C_1 + K_{m,n,k} \big)^{-3} - 3C t  \Big)^{-\frac{1}{3}},
\end{align*}
which is uniform with respect to $l$.
Arguing as in the derivation of \eqref{time_estimate_sp_with_alpha}, we obtain
\begin{align}
    & \big\| \partial_t \bm{u}^l \big\|^2 + \Big\|  \partial_t \big| \bm{u}^l \big|^\frac{r}{2}  \Big\|^2 \nonumber \\
    & \quad \le C_m \big(  \big\| \partial_t \mu_\phi^l \big\|_{-1} \big\| \nabla \phi^l \big\|_V + \big\| \partial_t \mu_\psi^l \big\|_{-1} \big\| \nabla \psi^l \big\|_V + \big\| \nabla \mu_\phi^l  \big\| \big\|  \nabla \partial_t \phi^l \big\| + \big\| \nabla \mu_\psi^l \big\| \big\| \nabla \partial_t \psi^l \big\|  \big)^2 \nonumber \\
    & \qquad + C_m \big(  1 + \big\|  \nabla \partial_t \phi^l \big\|^2 + \big\|  \nabla \partial_t \psi^l \big\|^2  \big) \big(  \big\| \bm{u}^l \big\|^2 + \big\| \bm{u}^l \big\|^r  \big). \label{time_estimate_sp}
\end{align}
To control the $V^*$-norm of $\partial_t \mu_\phi$, we differentiate equation \eqref{eq:second_approximate_chdf3} with respect to time. This gives
\begin{align*}
    ( \partial_t \mu_\phi^l, \zeta ) & = ( \nabla \partial_t \phi^l, \nabla \zeta ) + \sigma_2 \big(  \mathcal{N} \big( \partial_t \phi^l - \overline{\partial_t \phi^l} \big) , \zeta \big) \\
    & \quad + \big(  F_{\phi,n,\kappa}''(\phi^l) \partial_t \phi^l + \partial_{\phi}^2 G_\kappa(\phi^l,\psi^l) \partial_t \phi^l + \partial_{\phi,\psi}^2 G_\kappa(\phi^l,\psi^l) \partial_t \psi^l , \zeta \big) \\
    & \le \big\| \nabla \partial_t \phi^l \big\| \big\| \nabla \zeta \big\| + C \big\|  \partial_t \phi^l - \overline{\partial_t \phi^l}  \big\|_* \big\| \zeta - \overline{\zeta} \big\|_* + C_n \big(  \big\| \partial_t \phi^l \big\| + \big\| \partial_t \psi^l \big\|  \big) \big\| \zeta \big\| \\
    & \le C_n \big(  \big| \overline{\partial_t \phi^l} \big| + \big\| \nabla \partial_t \phi^l \big\| + \big\| \nabla \partial_t \psi^l \big\|  \big) \big\| \zeta \big\|_V,
\end{align*}
for all $\zeta \in X_l$, which implies 
\[\big\| \partial_t \mu_\phi^l \big\|_{-1} \le C_n \big(  \sigma_1 e^{-\sigma_1 t} + \big\| \nabla \partial_t \phi^l \big\| + \big\| \nabla \partial_t \psi^l \big\|  \big).\]
Analogously, we have
\[ \big\| \partial_t \mu_\psi^l \big\|_{-1} \le C_n \big(  \sigma_1 e^{-\sigma_1 t} + \big\| \nabla \partial_t \phi^l \big\| + \big\| \nabla \partial_t \psi^l \big\|  \big). \]
Hence, for all $t \in (0,T_e')$, $\partial_t \bm{u}^l$ and $\partial_t \big| \bm{u}^l \big|^\frac{r}{2}$ are uniformly bounded in $L^2\big( 0,t;L^2(\Omega) \big)$.
Moreover, taking $\zeta_1 = - \Delta^2 \mu_\phi^l$ in equation \eqref{eq:second_approximate_chdf2}, we find
\begin{align}
    & \int_\Omega m_\phi \big( \phi^l \big) \big|  \nabla \Delta \mu_\phi^l  \big|^2 \ \mathrm{d} x \nonumber \\
    & \quad = \big(  \nabla \partial_t \phi^l + \nabla \big(  \bm{u}^l \cdot \nabla \phi^l  \big) - \nabla^2 m_\phi \big( \phi^l \big) \cdot \nabla \mu_\phi^l - \nabla^2 \mu_\phi^l \cdot \nabla m_\phi \big( \phi^l \big) - \nabla m_\phi \big( \phi^l \big) \Delta \mu_\phi^l , \nabla \Delta \mu_\phi^l \big) \nonumber \\
    & \quad \le \big(  \big\|  \nabla \partial_t \phi^l  \big\| + \big\| \nabla \bm{u}^l \big\| \big\| \nabla \phi^l \big\|_{L^\infty(\Omega)} + \big\| \bm{u}^l \big\|_{L^3(\Omega)} \big\| \nabla^2 \phi^l \big\|_{L^6(\Omega)}  \big) \big\| \nabla \Delta \mu_\phi^l \big\| \nonumber \\
    & \qquad + C \big(  \big\| \nabla^2 \phi^l \big\|_{L^6(\Omega)} + \big\| \nabla \phi^l \big\|_{L^{12}(\Omega)}^2  \big) \big\| \nabla \mu_\phi^l \big\|_{L^3(\Omega)} \big\| \nabla \Delta \mu_\phi^l \big\| \nonumber \\
    & \qquad + C \big\| \nabla \phi^l \big\|_{L^\infty(\Omega)} \big\| \nabla^2 \mu_\phi^l \big\| \big\| \nabla \Delta \mu_\phi^l \big\|. \label{H3_mu}
\end{align}
Using the Gagliardo--Nirenberg inequality \eqref{gns_inequality}, we obtain
\begin{align*}
    \big\| \nabla \phi^l \big\|_{L^\infty(\Omega)} & \le C \big\| \nabla \phi^l \big\|_{V}^\frac{1}{2} \big\| \nabla \phi^l \big\|_{H^2(\Omega)}^\frac{1}{2} \\
    & \le C \big\| \Delta \phi^l \big\|^\frac{1}{2} \big\| \nabla \Delta \phi^l \big\|^\frac{1}{2} \\
    & \le C_n \big(  \big\| \nabla \phi^l \big\| + \big\| \nabla \psi^l \big\| +\big\| \nabla \mu_\phi^l \big\|  \big)^\frac{1}{4} \big\| \nabla \phi^l \big\|^\frac{1}{4} \big(  \big\| \nabla \mu_\phi^l \big\| + \big\| \nabla \phi^l \big\| + \big\| \nabla \psi^l \big\|  \big)^\frac{1}{2} \\
    & \le C_n \big(  1 + \big\| \nabla \mu_\phi^l \big\|^\frac{3}{4}  \big), \\
    \big\| \nabla \mu_\phi^l \big\|_{L^3(\Omega)} & \le C \big\| \nabla \mu_\phi^l \big\|^\frac{3}{4} \big\| \nabla \Delta \mu_\phi^l \big\|^\frac{1}{4}, \\
    \big\| \nabla^2 \mu_\phi^l \big\| & \le C \big\| \nabla \mu_\phi^l \big\|^\frac{1}{2} \big\| \nabla \Delta \mu_\phi^l \big\|^\frac{1}{2}.
\end{align*}
Combining the above estimates together, we get
\begin{align}
    & \int_\Omega m_\phi \big( \phi^l(\tau) \big) \big|  \nabla \Delta \mu_\phi^l(\tau)  \big|^2 \ \mathrm{d} x \nonumber \\
    & \quad \le C \big\| \nabla \partial_t \phi^l(\tau) \big\|^2 + C \big(  1 + C_1 + \mathcal{E}(\bm{u}^l(\tau),\phi^l(\tau),\psi^l(\tau),\mu_\phi^l(\tau),\mu_\psi^l(\tau))  \big)^\frac{8}{3},
\end{align}
which implies the uniform boundedness of $\nabla \mu_\phi^l$ in $L^2 \big( 0,t;H^2(\Omega) \big)$ with respect to $l$ for all $t \in (0,T_e')$.

Therefore, passing to the limit as $l \to \infty$, along a suitable subsequence, we obtain the desired (approximating) strong solution 
\[(\bm{u}_{m,n},\phi_{m,n},\psi_{m,n},\mu_{\phi,m,n},\mu_{\psi,m,n})\] 
to the system \eqref{eq:semigalerkin_chdf1}--\eqref{eq:semigalerkin_chdf5} on $[0,T_e')$. In particular, we have
\begin{align*}
    & \bm{u}^l \rightharpoonup \bm{u}_{m,n} \ \text{in} \ H^1\big( 0,t;\bm{V}_m \big), \\
    & \phi^l \rightharpoonup \phi_{m,n},\psi^l \rightharpoonup \psi_{m,n} \ \text{in} \ H^1 \big( 0,t;V \big), \\
    & \phi^l \stackrel{*}{\rightharpoonup} \phi_{m,n},\psi^l \stackrel{*}{\rightharpoonup} \psi_{m,n} \ \text{in} \ L^\infty \big( 0,t;H^3(\Omega) \big), \\
    & \mu_\phi^l \rightharpoonup \mu_{\phi,m,n}, \mu_\psi^l \rightharpoonup \mu_{\psi,m,n} \ \text{in} \ H^1 \big( 0,t;V^* \big) \cap L^2\big( 0,t;H^3(\Omega) \big),\\
    & \mu_\phi^l \stackrel{*}{\rightharpoonup }\mu_{\phi,m,n}, \mu_\psi^l \stackrel{*}{\rightharpoonup} \mu_{\psi,m,n} \ \text{in} \ L^\infty \big( 0,t;V \big),
\end{align*}
for all $t \in (0,T_e')$.
The proof is complete.
\qed

\begin{remark} \rm
    The fact that $\phi_{0,k} \in (-1,1)$ and $\psi_{0,k} \in (0,1)$ almost everywhere in $\Omega$ is not used in the above proof.
    Therefore, the well-posedness result holds for any initial data $(\phi_0,\psi_0) \in (H^3(\Omega) \cap W)^2$, a fact we will use to extend the approximate solutions in the following sections.
\end{remark}

\begin{remark} \rm
    The smoothness of $F_{\phi,n,\kappa}$, $F_{\psi,n,\kappa}$, and $G_\kappa$ ensures the uniqueness, which allows us to extend the approximate solutions.
    For the sake of simplicity, we will drop the subscript $\kappa$ and write $F_{\phi,n}$, $F_{\psi,n}$, and $G$ in the following sections. The corresponding estimates can still be carried out for the approximate system with $F_{\phi,n,\kappa}$, $F_{\psi,n,\kappa}$, and $G_\kappa$, and a subsequent compactness argument as $\kappa \to 0$ yields the existence result for the case of the non-mollified potentials $F_{\phi,n}$, $F_{\psi,n}$, and $G$.
\end{remark}

\medskip

\section{Proof of Theorem \ref{ss_chdf_exist_3d}}\label{ss_3d_proof}
\setcounter{equation}{0}

The proof is divided into three parts: (1) local existence in the general case, (2) global existence in the constant-coefficient case, (3) the strict separation property and uniqueness.

\smallskip
\noindent \textbf{Local existence.}
Taking $\bm{\theta}=\bm{u}_{m,n}$ in \eqref{eq:semigalerkin_chdf1} and testing \eqref{eq:semigalerkin_chdf2}--\eqref{eq:semigalerkin_chdf5} with $\mu_{\phi,m,n}$, $\partial_t \phi_{m,n}$, $\mu_{\psi,m,n}$ and $\partial_t \psi_{m,n}$, respectively, we  obtain the basic energy estimate
\begin{align}
    & \frac{\mathrm{d}}{\mathrm{d}t} \Big(  \frac{1}{2} \big\| \nabla \phi_{m,n} \big\|^2 + \frac{\beta}{2} \big\| \nabla \psi_{m,n} \big\|^2 + \int_\Omega F_{\phi,n}\big( \phi_{m,n} \big) + F_{\psi,n}\big( \psi_{m,n} \big) + G(\phi_{m,n},\psi_{m,n}) \ \mathrm{d} x + \frac{\sigma_2}{2} \big\|  \phi_{m,n} - \overline{\phi_{m,n}}  \big\|^2_*  \Big) \nonumber \\ 
    & \qquad + \int_\Omega \nu (\phi_{m,n},\psi_{m,n}) \big| \bm{u}_{m,n} \big|^2 + \eta (\phi_{m,n},\psi_{m,n}) \big| \bm{u}_{m,n} \big|^r \ \mathrm{d} x \nonumber \\
    & \qquad + \frac{1}{2} \int_\Omega m_\phi \big( \phi_{m,n} \big) \big| \nabla \mu_{\phi,m,n} \big|^2 + m_\psi \big( \psi_{m,n} \big) \big| \nabla \mu_{\psi,m,n} \big|^2 \ \mathrm{d} x \nonumber \\
    & \quad \le C \sigma_1 \big|  \overline{\phi_{m,n}}-c  \big| \big(  1 + \big\| \nabla \phi_{m,n} \big\|^2  \big). \label{basic_energy_estimate}
\end{align}
Recall that the construction of $\phi_{0,k}$, $\psi_{0,k}$ and $F_{\phi,n}$, $F_{\psi,n}$ implies
\[ F_{\phi,n}^{(i)} (\phi_{0,k}) = F_\phi^{(i)} (\phi_{0,k}), \ F_{\psi,n}^{(i)} (\psi_{0,k}) = F_\psi^{(i)} (\psi_{0,k}), \]
for all $i \in \{ 0,1,2 \}$.
Without loss of generality, we suppose that 
\begin{align*}
    & \| \nabla \phi_{0,k} \|^2 \le  \| \nabla \phi_0 \|^2 + 1, \ \| \nabla \psi_{0,k} \|^2 \le  \| \nabla \psi_0 \|^2 + 1, \\
    & \big|  \overline{\phi_{0,k}} -\overline{\phi_0}  \big| \le \frac{1}{2} \min \big\{ 1-\overline{\phi_0}, 1+\overline{\phi_0} \big\}, \ \big|  \overline{\psi_{0,k}} -\overline{\psi_0}  \big| \le \frac{1}{2} \min \big\{ \overline{\psi_0}, 1-\overline{\psi_0} \big\}.
\end{align*}
Therefore, we have
\begin{align*}
    \| \Delta \phi_{0,k} \|^2 & \le  (\| \nabla \phi_0 \|+1) \big\| \nabla (-\Delta \phi_0 + F_\phi'(\phi_0)) \big\|, \\
    \beta \| \Delta \psi_{0,k} \|^2 & \le  (\| \nabla \psi_0 \| + 1) \big\| \nabla (-\beta \Delta \psi_0 + F_\psi'(\psi_0)) \big\|, \\
    \| \nabla \mu_{\phi,m,n}(0) \| & = \big\|  \nabla (-\Delta \phi_{0,k} + F_\phi'(\phi_{0,k})) + \sigma_2 \nabla \mathcal{N} \big( \phi_{0,k} - \overline{\phi_{0,k}} \big) + \nabla \partial_\phi G(\phi_{0,k},\psi_{0,k})  \big\| \\
    & \le C \big( \big\| -\Delta \phi_0 + F_\phi'(\phi_0) \big\|_V  + \| \nabla \phi_0 \| + \| \nabla \psi_0 \| + 1 \big), \\
    \| \nabla \mu_{\psi,m,n}(0) \| & \le C \big( \big\| -\beta \Delta \psi_0 + F_\psi'(\psi_0) \big\|_V + \| \nabla \phi_0 \| + \| \nabla \psi_0 \| + 1 \big) , \\
    \nu_* \| \bm{u}_{m,n}(0) \|^2 + \eta_* \| \bm{u}_{m,n}(0) \|^r_{L^r(\Omega)} & \le C \| \phi_{0,k} \|_{L^\infty(\Omega)}^2 \| \nabla \mu_{\phi,m,n}(0) \|^2 + C \| \psi_{0,k} \|_{L^\infty(\Omega)}^2 \| \nabla \mu_{\psi,m,n}(0) \|^2 \\
    & \le C \big(  \big\| -\Delta \phi_0 + F_\phi'(\phi_0) \big\|_V^2 + \big\| -\beta \Delta \psi_0 + F_\psi'(\psi_0) \big\|_V^2  \big) \\
    & \quad + C \big( \| \nabla \phi_0 \|^2 + \| \nabla \psi_0 \|^2 + 1 \big).
\end{align*}
Repeating the derivation of \eqref{higher_order_estimate_approximate} and estimating the terms containing $\partial_t\nu$ and $\partial_t\eta$ using H\"older's inequality and the embedding $V\hookrightarrow L^6(\Omega)$, we obtain
\begin{align}
    & \frac{\mathrm{d}}{\mathrm{d}t} \mathcal{E} ( \bm{u}_{m,n},\phi_{m,n},\psi_{m,n},\mu_{\phi,m,n},\mu_{\psi,m,n} ) + \big\| \nabla \partial_t \phi_{m,n} \big\|^2 + \beta \big\| \nabla \partial_t \psi_{m,n} \big\|^2 \nonumber \\
    & \quad \le C \big(  \big\| \partial_t \phi_{m,n} \big\| \| \nabla \mu_{\phi,m,n} \|^2_{L^4(\Omega)} + \big\| \partial_t \psi_{m,n} \big\| \big\| \nabla \mu_{\psi,m,n} \big\|^2_{L^4(\Omega)}  \big) \nonumber \\
    & \qquad + C \big(  \big\| \partial_t \phi_{m,n} \big\|_V + \big\| \partial_t \psi_{m,n} \big\|_V  \big) \big(  \| \bm{u}_{m,n} \|^2_{L^\frac{12}{5}(\Omega)} + \| \bm{u}_{m,n} \|^r_{L^\frac{6r}{5}(\Omega)}  \big) \nonumber \\
    & \qquad + C \| \bm{u}_{m,n} \|_{L^3(\Omega)} \big(  \big\| \nabla \partial_t \phi_{m,n} \big\| \| \nabla \mu_{\phi,m,n} \| + \big\| \nabla \partial_t \psi_{m,n} \big\| \big\| \nabla \mu_{\psi,m,n} \big\|  \big) \nonumber \\
    & \qquad + C(1 + \| \mu_{\phi,m,n} \|) + C \big(  \big| \overline{\partial_t \phi_{m,n}} \big|^2 + \big\| \partial_t \phi_{m,n} - \overline{\partial_t \phi_{m,n}} \big\|^2 + \big\| \partial_t \psi_{m,n} \big\|^2  \big). \label{higher_order_estimate_uniform_3d}
\end{align}
Next, let us refine the estimates obtained in Section \ref{ss_approximate}.
Firstly, we have
\begin{align*}
    \big\| \partial_t \phi_{m,n} - \overline{\partial_t \phi_{m,n}} \big\| & \le C \big\| \nabla \partial_t \phi_{m,n} \big\|^\frac{1}{2} \big(  \| \nabla \mu_{\phi,m,n} \|^\frac{1}{2} + \| \bm{u}_{m,n} \|^\frac{1}{2} \big\| \nabla \phi_{m,n} \big\|_{L^3(\Omega)}^\frac{1}{2}  \big) \\
    & \le C \big\| \nabla \partial_t \phi_{m,n} \big\|^\frac{1}{2} \big(  \| \nabla \mu_{\phi,m,n} \|^\frac{1}{2} + \| \bm{u}_{m,n} \|^\frac{1}{2} \big(  1+ \| \nabla \mu_{\phi,m,n} \|  \big)^\frac{1}{8}  \big) \\
    & \le \sqrt{\epsilon} \big\| \nabla \partial_t \phi_{m,n} \big\| + C_\epsilon \big( 1 + C_1 + \mathcal{E} \big( \bm{u}_{m,n},\phi_{m,n},\psi_{m,n},\mu_{\phi,m,n},\mu_{\psi,m,n} \big) \big)^\frac{5}{8}, \\
    \big\| \partial_t \psi_{m,n} \big\| & \le C \big\| \nabla \partial_t \psi_{m,n} \big\|^\frac{1}{2} \big(  \big\| \nabla \mu_{\psi,m,n} \big\|^\frac{1}{2} + \| \bm{u}_{m,n} \|^\frac{1}{2} \big(  1+ \big\| \nabla \mu_{\psi,m,n} \big\|  \big)^\frac{1}{8}  \big) \\
    & \le \sqrt{\epsilon} \big\| \nabla \partial_t \psi_{m,n} \big\| + C_\epsilon \big( 1 + C_1 + \mathcal{E} \big( \bm{u}_{m,n},\phi_{m,n},\psi_{m,n},\mu_{\phi,m,n},\mu_{\psi,m,n} \big) \big)^\frac{5}{8}.
\end{align*}
Then, it holds
\begin{align*}
    & \int_\Omega m_\phi \big( \phi_{m,n} \big) \big| \Delta \mu_{\phi,m,n} \big|^2 \ \mathrm{d} x \\
    & \quad = \big(  \partial_t \phi_{m,n} - \overline{\partial_t \phi_{m,n}} + \bm{u}_{m,n} \cdot \nabla \phi_{m,n} - \nabla m_\phi\big( \phi_{m,n} \big) \cdot \nabla \mu_{\phi,m,n} , \Delta \mu_{\phi,m,n} \big)\\
    & \quad \le C \big\| \Delta \mu_{\phi,m,n} \big\| \Big(  \big\| \partial_t \phi_{m,n} - \overline{\partial_t \phi_{m,n}} \big\|^\frac{1}{2}_* \big\| \nabla \partial_t \phi_{m,n} \big\|^\frac{1}{2} \\
    & \qquad + \| \bm{u}_{m,n} \| \big\| \nabla \phi_{m,n} \big\|_{L^\infty(\Omega)} + \big\| \nabla \phi_{m,n} \big\|_{L^6(\Omega)} \| \nabla \mu_{\phi,m,n} \|_{L^3(\Omega)}  \Big) \\
    & \quad \le C \big\| \nabla \partial_t \phi_{m,n} \big\|^\frac{1}{2} \big(  \| \bm{u}_{m,n} \| \big(  1+ \| \nabla \mu_{\phi,m,n} \|  \big)^\frac{1}{4} + \| \nabla \mu_{\phi,m,n} \|  \big)^\frac{1}{2} \big\| \Delta \mu_{\phi,m,n} \big\| \\
    & \qquad + C \| \bm{u}_{m,n} \| \big(  \| \nabla \mu_{\phi,m,n} \|^\frac{1}{2} + 1  \big)^\frac{1}{2} \big(  1 + \big\| \nabla \mu_{\phi,m,n} \big\|  \big)^\frac{1}{2} \big\| \Delta \mu_{\phi,m,n} \big\| \\
    & \qquad + C \big(  \| \nabla \mu_{\phi,m,n} \|^\frac{1}{2} + 1  \big) \| \nabla \mu_{\phi,m,n} \|^\frac{1}{2} \big\| \Delta \mu_{\phi,m,n} \big\|^\frac{3}{2} \\
    & \quad \le \frac{1}{2} \int_\Omega m_\phi \big( \phi_{m,n} \big) \big| \Delta \mu_{\phi,m,n} \big|^2 \ \mathrm{d} x + \frac{1}{2} \big\| \nabla \partial_t \phi_{m,n} \big\|^\frac{3}{2} + C \big(  1 + \| \bm{u}_{m,n} \|^\frac{9}{2} + \| \nabla \mu_{\phi,m,n} \|^\frac{9}{2}  \big),
\end{align*}
and, similarly for $\mu_{\psi,m,n}$, we find
\begin{align*}
    \int_\Omega m_\psi \big( \psi_{m,n} \big) \big| \Delta \mu_{\psi,m,n} \big|^2 \ \mathrm{d} x \le \frac{1}{2} \big\| \nabla \partial_t \psi_{m,n} \big\|^\frac{3}{2} + C \big(  1 + \| \bm{u}_{m,n} \|^\frac{9}{2} + \| \nabla \mu_{\psi,m,n} \|^\frac{9}{2}  \big),
\end{align*}
which implies
\begin{align*}
    & \big\| \partial_t \phi_{m,n} \big\| \| \nabla \mu_{\phi,m,n} \|^2_{L^4(\Omega)} \\
    & \quad \le C \Big(  \sigma_1 e^{-\sigma_1 t} + \big\| \nabla \partial_t \phi_{m,n} \big\|^\frac{1}{2} \big(  \| \nabla \mu_{\phi,m,n} \|^\frac{1}{2} + \| \bm{u}_{m,n} \|^\frac{1}{2} \big(  1+ \| \nabla \mu_{\phi,m,n} \|  \big)^\frac{1}{8}  \big)  \Big) \| \nabla \mu_{\phi,m,n} \|^\frac{1}{2} \big\| \Delta \mu_{\phi,m,n} \big\|^\frac{3}{2} \\
    & \quad \le C \big(  \| \nabla \mu_{\phi,m,n} \|^\frac{1}{2} \big\| \nabla \partial_t \phi_{m,n} \big\|^\frac{9}{8} + \| \nabla \mu_{\phi,m,n} \|^\frac{1}{2} + \| \bm{u}_{m,n} \|^\frac{31}{8} + \| \nabla \mu_{\phi,m,n} \|^\frac{31}{8}  \big) \\
    & \qquad + C \big\| \nabla \partial_t \phi_{m,n} \big\|^\frac{13}{8} \big(  1 + \| \bm{u}_{m,n} \|^\frac{9}{8} + \| \nabla \mu_{\phi,m,n} \|^\frac{9}{8}  \big) + C \big\| \nabla \partial_t \phi_{m,n} \big\|^\frac{1}{2} \big(  1 + \| \bm{u}_{m,n} \|^\frac{9}{2} + \| \nabla \mu_{\phi,m,n} \|^\frac{9}{2}  \big) \\
    & \quad \le \epsilon  \big\| \nabla \partial_t \phi_{m,n} \big\|^2 + C_\epsilon \big(  1 + C_1 + \mathcal{E} \big( \bm{u}_{m,n},\phi_{m,n},\psi_{m,n},\mu_{\phi,m,n},\mu_{\psi,m,n} \big)  \big)^3.
\end{align*}
Testing the approximate momentum equation \eqref{eq:semigalerkin_chdf1} with $ \bm{A} \bm{u}_{m,n} = (I + \bm{A}) \bm{u}_{m,n} - \bm{u}_{m,n} \in \bm{V}_m$, we have
\begin{align*}
    & \int_\Omega \nu(\phi_{m,n},\psi_{m,n}) | \nabla \bm{u}_{m,n} |^2 + \eta(\phi_{m,n},\psi_{m,n}) | \bm{u}_{m,n} |^{r-2} | \nabla \bm{u}_{m,n} |^2 \ \mathrm{d} x \\
    & \qquad + \int_\Omega (r-2) \eta(\phi_{m,n},\psi_{m,n}) | \bm{u}_{m,n} |^{r-2} | \nabla \bm{u}_{m,n} \cdot \bm{e}_{\mathbb{R}^d} (\bm{u}_{m,n}) |^2 \ \mathrm{d} x \\
    & \quad = \int_{\partial \Omega} \big(  \nu(\phi_{m,n},\psi_{m,n}) + \eta(\phi_{m,n},\psi_{m,n}) | \bm{u}_{m,n} |^{r-2}  \big) \bm{u}_{m,n} \cdot \partial_{\bm{n}} \bm{u}_{m,n} \ \mathrm{d} \Gamma \\
    & \qquad - \int_\Omega (\nabla \nu(\phi_{m,n},\psi_{m,n}) \cdot \nabla) \bm{u}_{m,n} \cdot \bm{u}_{m,n} + | \bm{u}_{m,n} |^{r-2} (\nabla \eta(\phi_{m,n},\psi_{m,n}) \cdot \nabla) \bm{u}_{m,n} \cdot \bm{u}_{m,n} \ \mathrm{d} x \\
    & \qquad + (  \nabla \mu_{\phi,m,n} \times \nabla \phi_{m,n} + \nabla \mu_{\psi,m,n} \times \nabla \psi_{m,n} , \nabla \times \bm{u}_{m,n} ).
\end{align*}
First, for any $\bm{u} \in \bm{V}_m$, it holds
\begin{align*}
    \bm{u} \cdot \partial_{\bm{n}} \bm{u} 
    &  = \bm{n} \cdot (\bm{u} \cdot \nabla) \bm{u} + \bm{u} \cdot \big(  (  (\nabla \bm{u})^T - \nabla \bm{u}  ) \cdot \bm{n}  \big) \\
    &  = \bm{n} \cdot (\bm{u} \cdot \nabla_{\partial \Omega}) \bm{u} + \bm{u} \cdot \big(  (\nabla \times \bm{u}) \times \bm{n}  \big) \\
    &  = (\bm{u} \cdot \nabla_{\partial \Omega}) (\bm{u} \cdot \bm{n}) - \bm{u} \cdot (\bm{u} \cdot \nabla_{\partial \Omega}) \bm{n} \\
    &  = \Pi(\bm{u},\bm{u}),
\end{align*}
where we use the boundary conditions of the eigenfunctions of $I + \bm{A}$, that is,
\[ \bm{u} \cdot \bm{n} =0, \quad \bm{n} \times (\nabla \times \bm{u}) = \bm{0}, \quad \text{a.e. on $\partial \Omega$}. \]
Here, $\Pi$ is the second fundamental form of the boundary $\partial \Omega$ (see \cite[Eq.~(3,1,1,8)]{G1985} and \cite[Section 4]{GHW2026} for its properties and the application).
This implies that
\begin{align*}
    & \big| \int_{\partial \Omega} \big(  \nu(\phi_{m,n},\psi_{m,n}) + \eta(\phi_{m,n},\psi_{m,n}) | \bm{u}_{m,n} |^{r-2}  \big) \bm{u}_{m,n} \cdot \partial_{\bm{n}} \bm{u}_{m,n} \ \mathrm{d} \Gamma \big| \\
    & \quad = \Big| \int_{\partial \Omega} \big(  \nu(\phi_{m,n},\psi_{m,n}) + \eta(\phi_{m,n},\psi_{m,n}) | \bm{u}_{m,n} |^{r-2}  \big) \Pi(\bm{u}_{m,n},\bm{u}_{m,n}) \, \mathrm{d} \Gamma \Big| \\
    & \quad \le C \| \Pi \|_{L^\infty(\partial \Omega)}(\| \bm{u}_{m,n} \|^2_{\bm{L}^2 (\partial \Omega)} +  \| \bm{u}_{m,n} \|^r_{\bm{L}^r (\partial \Omega)}) \\
    & \quad = C \| \Pi \|_{L^\infty(\partial \Omega)} (\| \bm{u}_{m,n} \|^2_{\bm{L}^2 (\partial \Omega)} + \| | \bm{u}_{m,n} |^{\frac{r}{2}} \|^2_{L^2 (\partial \Omega)}) \\
    & \quad \le C \| \Pi \|_{L^\infty(\partial \Omega)} (\| \bm{u}_{m,n} \| \| \bm{u}_{m,n} \|_{V} + \| | \bm{u}_{m,n} |^{\frac{r}{2}} \| \| | \bm{u}_{m,n} |^{\frac{r}{2}} \|_{V}) \\
    & \quad \le C \| \Pi \|_{L^\infty(\partial \Omega)} (\| \bm{u}_{m,n} \|^2 + \| | \bm{u}_{m,n} |^{\frac{r}{2}} \|^2) \\
    & \qquad + C \| \Pi \|_{L^\infty(\partial \Omega)} (\| \bm{u}_{m,n} \| \| \nabla \bm{u}_{m,n} \| + \| | \bm{u}_{m,n} |^{\frac{r}{2}} \| \| \nabla | \bm{u}_{m,n} |^{\frac{r}{2}} \|) \\
    & \quad \le \frac{\nu_*}{4} \| \nabla \bm{u}_{m,n} \|^2 + \frac{(r-1)\eta_*}{r^2} \| \nabla | \bm{u}_{m,n} |^{\frac{r}{2}} \|^2 + C \| \bm{u}_{m,n} \|^r_{\bm{L}_\sigma^r(\Omega)} + C \| \bm{u}_{m,n} \|^2.
\end{align*}

On the other hand, it follows from Lemma \ref{chain_rule} that
\begin{align*}
    & \int_\Omega \eta(\phi_{m,n},\psi_{m,n}) | \bm{u}_{m,n} |^{r-2} | \nabla \bm{u}_{m,n} |^2 \ \mathrm{d} x \\
    & \quad \ge \int_\Omega \eta(\phi_{m,n},\psi_{m,n}) | \bm{u}_{m,n} |^{r-2} | \nabla \bm{u}_{m,n} \cdot \bm{e}_{\mathbb{R}^d}(\bm{u}_{m,n}) |^2 \ \mathrm{d} x \\
    & \quad = \frac{4}{r^2} \int_\Omega \eta(\phi_{m,n},\psi_{m,n}) \big|  \nabla | \bm{u}_{m,n} |^\frac{r}{2}  \big|^2 \ \mathrm{d} x, \\
    & \int_\Omega| \bm{u}_{m,n} |^{r-2} (\nabla \eta(\phi_{m,n},\psi_{m,n}) \cdot \nabla) \bm{u}_{m,n} \cdot \bm{u}_{m,n} \ \mathrm{d} x \\
    & \quad = \frac{2}{r} \int_\Omega| \bm{u}_{m,n} |^\frac{r}{2} (\nabla \eta(\phi_{m,n},\psi_{m,n}) \cdot \nabla) | \bm{u}_{m,n} |^\frac{r}{2} \ \mathrm{d} x,
\end{align*}
which implies
\begin{align*}
    & \Big|  \int_\Omega| \bm{u}_{m,n} |^{r-2} (\nabla \eta(\phi_{m,n},\psi_{m,n}) \cdot \nabla) \bm{u}_{m,n} \cdot \bm{u}_{m,n} \ \mathrm{d} x  \Big| \\
    & \quad \le C \big\|  | \bm{u}_{m,n} |^\frac{r}{2}  \big\|_{L^3(\Omega)} \big\|  \nabla | \bm{u}_{m,n} |^\frac{r}{2}  \big\| \big(  \| \nabla \phi_{m,n} \|_{L^6(\Omega)} + \| \nabla \psi_{m,n} \|_{L^6(\Omega)}  \big) \\
    & \quad \le C \big\| | \bm{u}_{m,n} |^\frac{r}{2} \big\|^\frac{1}{2} \big\| | \bm{u}_{m,n} |^\frac{r}{2} \big\|_V^\frac{3}{2} \big(  1 + \| \nabla \mu_{\phi,m,n} \|^\frac{1}{2} + \| \nabla \mu_{\psi,m,n} \|^\frac{1}{2}  \big).
\end{align*}
Therefore, it holds
\begin{align}
    & \| \bm{u}_{m,n} \|_V^2 + \big\|  | \bm{u}_{m,n} |^\frac{r}{2}  \big\|_V^2 \nonumber \\
    & \quad \le C \big(  \| \bm{u}_{m,n} \|_{\bm{L}^r(\Omega)}^r + \| \bm{u}_{m,n} \|_{\bm{L}^2(\Omega)}^2  \big) \big(  1 + \| \nabla \mu_{\phi,m,n} \|^2 + \| \nabla \mu_{\psi,m,n} \|^2  \big) \nonumber \\
    & \qquad + C\big(  1 + \| \nabla \mu_{\phi,m,n} \|^\frac{7}{2} + \| \nabla \mu_{\psi,m,n} \|^\frac{7}{2}  \big) \nonumber \\
    & \quad \le C \big( 1 + C_1 + \mathcal{E} ( \bm{u}_{m,n},\phi_{m,n},\psi_{m,n},\mu_{\phi,m,n},\mu_{\psi,m,n} ) \big)^2,\label{H1_sp_3d}
\end{align}
where we used the estimate
\begin{align*}
    & \| \nabla \phi_{m,n} \|_{\bm{L}^\infty(\Omega)} \le C \| \Delta \phi_{m,n} \|^\frac{1}{2} \big\| \Delta \phi_{m,n} \big\|_{L^6(\Omega)}^\frac{1}{2}  \le C (1 + \| \nabla \mu_{\phi,m,n} \|)^\frac{3}{4}, \\
    & \| \nabla \psi_{m,n} \|_{\bm{L}^\infty(\Omega)} \le C (1 + \| \nabla \mu_{\psi,m,n} \|)^\frac{3}{4}.
\end{align*}
Then, we estimate the second and third terms on the right-hand side of \eqref{higher_order_estimate_uniform_3d}.
It follows from the interpolation that
\begin{align*}
    & \big(  \big\| \partial_t \phi_{m,n} \big\|_V + \big\| \partial_t \psi_{m,n} \big\|_V  \big) \big(  \| \bm{u}_{m,n} \|^2_{L^\frac{12}{5}(\Omega)} + \| \bm{u}_{m,n} \|^r_{L^\frac{6r}{5}(\Omega)}  \big) \\
    & \quad \le C \big(  \sigma_1 e^{-\sigma_1 t} | \overline{\phi_{0,k}} - c | + \big\| \nabla \partial_t \phi_{m,n} \big\| + \big\| \nabla \partial_t \psi_{m,n} \big\|  \big) \big(  \| \bm{u}_{m,n} \|^\frac{3}{2} \| \bm{u}_{m,n} \|_V^\frac{1}{2} + \big\| | \bm{u}_{m,n} |^\frac{r}{2} \big\|^\frac{3}{2} \big\| | \bm{u}_{m,n} |^\frac{r}{2} \big\|_V^\frac{1}{2}  \big) \\
    & \quad \le \epsilon \big(  \| \nabla \partial_t \phi_{m,n} \|^2 + \beta \| \nabla \partial_t \psi_{m,n} \|^2  \big) + C_\epsilon \big( 1 + C_1 + \mathcal{E} ( \bm{u}_{m,n},\phi_{m,n},\psi_{m,n},\mu_{\phi,m,n},\mu_{\psi,m,n} ) \big)^\frac{5}{2}, \\
    & \| \bm{u}_{m,n} \|_{L^3(\Omega)} \big(  \big\| \nabla \partial_t \phi_{m,n} \big\| \| \nabla \mu_{\phi,m,n} \| + \big\| \nabla \partial_t \psi_{m,n} \big\| \big\| \nabla \mu_{\psi,m,n} \big\|  \big) \\
    & \quad \le C \| \bm{u}_{m,n} \|^\frac{1}{2} \| \bm{u}_{m,n} \|_V^\frac{1}{2} \big(  \big\| \nabla \partial_t \phi_{m,n} \big\| \| \nabla \mu_{\phi,m,n} \| + \big\| \nabla \partial_t \psi_{m,n} \big\| \big\| \nabla \mu_{\psi,m,n} \big\|  \big) \\
    & \quad \le \epsilon \big(  \| \nabla \partial_t \phi_{m,n} \|^2 + \beta \| \nabla \partial_t \psi_{m,n} \|^2  \big) + C_\epsilon \big( 1 + C_1 + \mathcal{E} ( \bm{u}_{m,n},\phi_{m,n},\psi_{m,n},\mu_{\phi,m,n},\mu_{\psi,m,n} ) \big)^{\frac{5}{2}}.
\end{align*}
Summing up the above estimates, we have
\begin{align}
    & \frac{\mathrm{d}}{\mathrm{d}t} \mathcal{E} ( \bm{u}_{m,n},\phi_{m,n},\psi_{m,n},\mu_{\phi,m,n},\mu_{\psi,m,n} )  + \frac{1}{2} \big\| \nabla \partial_t \phi_{m,n} \big\|^2 + \frac{\beta}{2} \big\| \nabla \partial_t \psi_{m,n} \big\|^2 \nonumber \\
    & \quad \le \hypertarget{C_2_energy_higher_order}{C_2} \big( 1 + C_1 + \mathcal{E} ( \bm{u}_{m,n},\phi_{m,n},\psi_{m,n},\mu_{\phi,m,n},\mu_{\psi,m,n} ) \big)^\frac{9}{2}, \label{energy_estimate_local_higher_order}
\end{align}
where $C_2 >0$ depends on $E_{\mathrm{free}}(\phi_0,\psi_0)$, $F_\phi$, $F_\psi$, $G$, $m_\phi$, $m_\psi$, $\nu$, $\eta$, and $\Omega$.
Note that the initial value of $\mathcal{E} ( \bm{u}_{m,n},\phi_{m,n},\psi_{m,n},\mu_{\phi,m,n},\mu_{\psi,m,n} )$ is uniformly bounded with respect to $m,n,k$, that is,
\[
\mathcal{E} ( \bm{u}_{m,n}(0),\phi_{0,k},\psi_{0,k},\mu_{\phi,m,n}(0),\mu_{\psi,m,n}(0) ) \le K_0,
\]
where
\begin{align*}
    K_0 & \stackrel{\rm{def}}{=} 3 \big(  \frac{3}{4} m_\phi^* + \frac{\nu^*}{\nu_*^2} + \frac{\eta^*}{\nu_* \eta_*}  \big) \|  -\Delta \phi_0 + F_\phi'(\phi_0)  \|_V^2 + \big(  m_\psi^* + \frac{2 \nu^*}{\nu_*^2} + \frac{2 \eta^*}{\nu_* \eta_*}  \big) \|  -\beta \Delta \psi_0 + F_\psi'(\psi_0)  \|_V^2 \\
    & \quad + \Big(  3 \sigma_2^2 C_P^4 \big(  \frac{3}{4} m_\phi^* + \frac{\nu^*}{\nu_*^2} + \frac{\eta^*}{\nu_* \eta_*}  \big) + \| G \|_{C^2(\mathbb{R}^2)}^2 \big(  \frac{9}{2} m_\phi^* + 2 m_\psi^* + \frac{10 \nu^*}{\nu_*^2} + \frac{10 \eta^*}{\nu_* \eta_*}  \big)  \Big) (1 + \| \nabla \phi_0 \|^2) \\
    & \quad + 4 C_{\mathrm{source}} \sigma_1 (1 + \| \nabla \phi_0 \|^2) \\
    & \quad + \| G \|_{C^2(\mathbb{R}^2)}^2 \big(  \frac{9}{2} m_\phi^* + 2 m_\psi^* + \frac{10 \nu^*}{\nu_*^2} + \frac{10 \eta^*}{\nu_* \eta_*}  \big) (1 + \| \nabla \psi_0 \|^2).
\end{align*}
Thanks to uniqueness, we can extend  $(\bm{u}_{m,n},\phi_{m,n},\psi_{m,n},\mu_{\phi,m,n},\mu_{\psi,m,n})$ to the interval $[0,T_0)$, where
\[ T_0 \stackrel{\rm{def}}{=} \frac{2}{7C_2} \big(  1 + C_1 + K_0  \big)^{-\frac{7}{2}}. \]
In addition, for all $t \in [0,T_0)$, we have the bound
\begin{align*}
    & 1 + C_1 + \mathcal{E}(\bm{u}_{m,n}(t),\phi_{m,n}(t),\psi_{m,n}(t),\mu_{\phi,m,n}(t),\mu_{\psi,m,n}(t)) + \int_0^t \frac{1}{2} \big\| \nabla \partial_t \phi_{m,n} \big\|^2 + \frac{\beta}{2} \big\| \nabla \partial_t \psi_{m,n} \big\|^2 \ \mathrm{d} \tau \\
    & \quad \le \Big(  \big( 1 + C_1 + K_0 \big)^{-\frac{7}{2}} - \frac{7C_2}{2} t  \Big)^{-\frac{2}{7}}.
\end{align*}
Moreover, it follows from \eqref{H3_mu} that
\begin{align*}
    & \int_\Omega m_\phi \big( \phi_{m,n} \big) \big|  \nabla \Delta \mu_{\phi,m,n}  \big|^2 \ \mathrm{d} x \\
    & \quad \le \big(  \big\|  \nabla \partial_t \phi_{m,n}  \big\| + \big\| \nabla \bm{u}_{m,n} \big\| \big\| \nabla \phi_{m,n} \big\|_{L^\infty(\Omega)} + \big\| \bm{u}_{m,n} \big\|_{L^3(\Omega)} \big\| \nabla^2 \phi_{m,n} \big\|_{L^6(\Omega)}  \big) \big\| \nabla \Delta \mu_{\phi,m,n} \big\| \\
    & \qquad + C \big(  \big\| \nabla^2 \phi_{m,n} \big\|_{L^6(\Omega)} + \big\| \nabla \phi_{m,n} \big\|_{L^{12}(\Omega)}^2  \big) \big\| \nabla \mu_{\phi,m,n} \big\|_{L^3(\Omega)} \big\| \nabla \Delta \mu_{\phi,m,n} \big\| \\
    & \qquad + C \big\| \nabla \phi_{m,n} \big\|_{L^\infty(\Omega)} \big\| \nabla^2 \mu_{\phi,m,n} \big\| \big\| \nabla \Delta \mu_{\phi,m,n} \big\| \\
    & \quad \le C \big\|  \nabla \partial_t \phi_{m,n}  \big\|^2 + \frac{1}{4} \int_\Omega m_\phi \big( \phi_{m,n} \big) \big|  \nabla \Delta \mu_{\phi,m,n}  \big|^2 \ \mathrm{d} x \\
    & \qquad + C \big( 1 + C_1 + \mathcal{E} ( \bm{u}_{m,n},\phi_{m,n},\psi_{m,n},\mu_{\phi,m,n},\mu_{\psi,m,n} ) \big)^\frac{11}{4}.
\end{align*}
In conclusion, for all $t \in (0,T_0)$, we have
\begin{align*}
    & \bm{u}_{m,n} \in L^\infty \big( 0,t;\bm{L}^{3r}_\sigma(\Omega) \cap \bm{H}^1(\Omega) \big),  \quad \text{with} \ | \bm{u}_{m,n} |^\frac{r}{2} \in L^\infty \big( 0,t;H^1(\Omega) \big), \\
    & \phi_{m,n},\psi_{m,n} \in L^\infty \big( 0,t;W^{2,6}(\Omega) \big) \cap H^1 \big( 0,t;V \big) \cap W^{1,\infty} \big( 0,t;V^* \big), \\
    & \mu_{\phi,m,n},\mu_{\psi,m,n} \in L^\infty \big( 0,t;V \big) \cap L^2 \big( 0,t;H^3(\Omega) \big),
\end{align*}
where the functions are uniformly bounded in the corresponding spaces with respect to $m,n,k$.
Here the Sobolev embedding $H^1(\Omega) \hookrightarrow L^6(\Omega)$ and estimate \eqref{H1_sp_3d} are used to control the $\bm{L}^{3r}(\Omega)$-norm of $\bm{u}_{m,n}$, i.e., the $L^6(\Omega)$-norm of $| \bm{u}_{m,n} |^{\tfrac{r}{2}}$.
So we can first pass to the limit as $n \to \infty$ and then as $m \to \infty$ and finally as $k \to \infty$ to obtain a strong solution in $[0,T_0)$ using standard compactness results combined with a diagonal argument.
The convergence of the Forchheimer term can be proved by a standard monotonicity argument and the energy identity \eqref{basic_energy_identity}; see \cite{GOPW2026}. 
Indeed, this monotonicity argument verifies
\[
    \limsup \int_{Q_t} \eta (\phi_{m,n},\psi_{m,n}) | \bm{u}_{m,n} |^r \, \mathrm{d} (x,\tau) \le \int_{Q_t} \eta (\widehat{\phi},\widehat{\psi}) | \widehat{\bm{u}} |^r \, \mathrm{d} (x,\tau),
\]
where the limit procedure is conducted as described above and the hat notation is used to indicate the corresponding limit.
It follows from the strong convergence of $\phi_{m,n}$ and $\psi_{m,n}$ in $C([0,t];H^s(\Omega))$ for all $s \in [0,2)$ and weak convergence of $\bm{u}_{m,n}$ in $L^r(Q_t)$ that
\[
    \lim \int_{Q_t} \eta (\phi_{m,n},\psi_{m,n}) | \bm{u}_{m,n} |^r \, \mathrm{d} (x,\tau) = \int_{Q_t} \eta (\widehat{\phi},\widehat{\psi}) | \widehat{\bm{u}} |^r \, \mathrm{d} (x,\tau),
\]
which implies that 
\[
    \eta^\frac{1}{r} (\phi_{m,n},\psi_{m,n}) \bm{u}_{m,n} \to \eta^\frac{1}{r} (\widehat{\phi},\widehat{\psi}) \widehat{\bm{u}} \quad \text{in} \ L^r(Q_t),
\]
because $L^r(Q_t)$ is uniformly convex.
Therefore, using again the strong convergence of $\phi_{m,n}$ and $\psi_{m,n}$, we can conclude that $\bm{u}_{m,n} \to \widehat{\bm{u}}$ in $L^r(Q_t)$ for all $t \in (0,T_0)$.
This implies that the weak-* limit of $| \bm{u}_{m,n} |^\frac{r}{2}$ in $L^\infty(0,t;H^1(\Omega))$ is exactly $| \widehat{\bm{u}} |^\frac{r}{2}$.

In conclusion, we can take $T_e = \frac{2}{3} T_0$. 

\medskip

\noindent \textbf{Global existence.}
In this part, we assume that the coefficients $m_\phi$, $m_\psi$, $\nu$ and $\eta$ are positive constants and $r\ge 3$.  
Then we prove the global existence by refining the estimates in Section \ref{ss_approximate}.

The basic energy estimate is the same as in the previous proof of the local existence.
Since the coefficients are constants, the higher-order estimate \eqref{higher_order_estimate_uniform_3d} reduces to the following form:
\begin{align*}
    & \frac{\mathrm{d}}{\mathrm{d}t} \mathcal{E} ( \bm{u}_{m,n},\phi_{m,n},\psi_{m,n},\mu_{\phi,m,n},\mu_{\psi,m,n} ) + \big\| \nabla \partial_t \phi_{m,n} \big\|^2 + \beta \big\| \nabla \partial_t \psi_{m,n} \big\|^2 \nonumber \\
    & \quad \le C \| \bm{u}_{m,n} \|_{L^3(\Omega)} \big(  \big\| \nabla \partial_t \phi_{m,n} \big\| \| \nabla \mu_{\phi,m,n} \| + \big\| \nabla \partial_t \psi_{m,n} \big\| \big\| \nabla \mu_{\psi,m,n} \big\|  \big) + C(1 + \| \mu_{\phi,m,n} \|) \nonumber \\
    & \qquad + C \big(  \big| \overline{\partial_t \phi_{m,n}} \big|^2 + \big\| \partial_t \phi_{m,n} - \overline{\partial_t \phi_{m,n}} \big\|_* \| \nabla \partial_t \phi_{m,n} \| + \big\| \partial_t \psi_{m,n} \big\|_* \big\| \nabla \partial_t \psi_{m,n} \big\|  \big) \nonumber \\
    & \quad \le \frac{1}{2} \big(  \big\| \nabla \partial_t \phi_{m,n} \big\|^2 + \beta \big\| \nabla \partial_t \psi_{m,n} \big\|^2  \big) \nonumber \\
    & \qquad + C \| \bm{u}_{m,n} \|_{L^3(\Omega)}^2 \big(  \| \nabla \mu_{\phi,m,n} \|^2 + \| \nabla \mu_{\psi,m,n} \|^2  \big) \nonumber \\
    & \qquad + C \big(  \| \bm{u}_{m,n} \|_{L^3(\Omega)}^2 + \| \nabla \mu_{\phi,m,n} \|^2 + \| \nabla \mu_{\psi,m,n} \|^2  \big) + C(1 + \| \mu_{\phi,m,n} \|),
\end{align*}
which implies
\begin{align}
    & \frac{\mathrm{d}}{\mathrm{d}t} \mathcal{E} ( \bm{u}_{m,n},\phi_{m,n},\psi_{m,n},\mu_{\phi,m,n},\mu_{\psi,m,n} ) + \frac{1}{2} \big\| \nabla \partial_t \phi_{m,n} \big\|^2 + \frac{\beta}{2} \big\| \nabla \partial_t \psi_{m,n} \big\|^2 \nonumber \\
    & \quad \le C \big(  1 + C_1 + \mathcal{E} ( \bm{u}_{m,n},\phi_{m,n},\psi_{m,n},\mu_{\phi,m,n},\mu_{\psi,m,n} )  \big)^2. \label{higher_order_estimate_uniform_3d_constant}
\end{align}
This allows us to extend the approximate solution $(\bm{u}_{m,n},\phi_{m,n},\psi_{m,n},\mu_{\phi,m,n},\mu_{\psi,m,n})$ to $[0,\infty)$.
Then, since the energy estimate \eqref{basic_energy_estimate} ensures the uniformly local integrability of the dissipative terms, 
which is independent of the approximation parameters $m,n,k$, 
we have that $\mathcal{E} ( \bm{u}_{m,n},\phi_{m,n},\psi_{m,n},\mu_{\phi,m,n},\mu_{\psi,m,n} )$ is uniformly locally integrable in $[0,\infty)$.
Thus, the Gronwall inequality, the uniform Gronwall inequality and comparison in equations \eqref{eq:semigalerkin_chdf2}--\eqref{eq:semigalerkin_chdf5} yield
\begin{align*} 
    & \bm{u}_{m,n} \in L^\infty\big( 0,\infty;\bm{L}^r_\sigma(\Omega) \big), \\
    & \phi_{m,n}, \psi_{m,n} \in L^\infty \big( 0,\infty;W^{2,6}(\Omega) \cap W \big) \cap H^1_{\mathrm{uloc}}\big( 0,\infty;V \big) \cap W^{1,\infty}\big( 0,\infty;V^* \big),\\
    & \mu_{\phi,m,n},\mu_{\psi,m,n} \in L^\infty\big( 0,\infty;V \big) \cap L^2_{\mathrm{uloc}} \big( 0,\infty;W \big),
\end{align*}
and the above functions are uniformly bounded in the corresponding spaces with respect to $m$, $n$ and $k$.
Therefore, the same calculation as in the previous proof of the local existence gives the uniform $L^2_{\mathrm{uloc}} (0,\infty;H^3(\Omega))$-estimate of  $\mu_{\phi,m,n}$, $\mu_{\psi,m,n}$, the uniform $L^\infty (0,\infty;\bm{H}^1(\Omega))$-estimate of $\bm{u}_{m,n}$, and the uniform $L^\infty(0,\infty;H^1(\Omega))$-estimate of $| \bm{u}_{m,n} |^\frac{r}{2}$, with respect to $m$, $n$, and $k$. 
Then, the existence of a global-in-time strong solution also follows from a standard compactness argument.

\medskip

\noindent \textbf{Separation property and uniqueness.} 
It follows from \eqref{gns_inequality} that there exists some constant $C_{GN}>0$ depending on $\Omega$ such that 
\[
\| f \|_{C( \overline{\Omega} )} \le C_{GN} \| f \|_V^{\frac{3}{4}} \| f \|_{W^{2,6}(\Omega)}^{\frac{1}{4}},
\]
for all $f \in W^{2,6}(\Omega)$.
Let $T_s' \in I$ be a positive time in the existence interval $I$ constructed in the previous proof.
Then it holds
\begin{align*}
    \| \phi(t)-\phi_0 \|_{C( \overline{\Omega} )} & \le C_{GN} \| \phi(t)-\phi_0 \|_V^{\frac{3}{4}} \| \phi(t)-\phi_0 \|_{W^{2,6}(\Omega)}^{\frac{1}{4}} \\
    & \le 2 C_{GN} \| \phi \|_{L^\infty( 0,\tfrac{1}{2} T_s';W^{2,6}(\Omega) )}^{\frac{1}{4}} \| \phi \|_{H^1( 0,\tfrac{1}{2} T_s';V )}^{\frac{3}{4}} t^{\frac{3}{8}}, \\
    \| \psi(t)-\psi_0 \|_{C( \overline{\Omega} )} & \le 2 C_{GN} \| \psi \|_{L^\infty( 0,\tfrac{1}{2} T_s';W^{2,6}(\Omega) )}^{\frac{1}{4}} \| \psi \|_{H^1( 0,\tfrac{1}{2} T_s';V )}^{\frac{3}{4}} t^{\frac{3}{8}}, 
\end{align*}
for all $t \in [0,\tfrac{1}{2} T_s']$.
Therefore, if there exists $\delta \in (0,\tfrac{1}{2})$ such that
\[
\| \phi_0 \|_{L^\infty(\Omega)} \le 1 - \delta, \quad \| \psi_0 - \frac{1}{2} \|_{L^\infty(\Omega)} \le \frac{1}{2} - \delta,
\]
then we have
\[
\| \phi(t) \|_{L^\infty(\Omega)} \le 1 - \frac{9}{10} \delta, \quad \| \psi(t) - \frac{1}{2} \|_{L^\infty(\Omega)} \le \frac{1}{2} - \frac{9}{10} \delta,
\]
for all $t \in [0,T_s]$, where
\begin{align}
    \hypertarget{T_s_delta}{T_s} \stackrel{\mathrm{def}}{=} \min \Big\{  & \frac{T_s'}{2}, \Big(  1 + 2 C_{GN} \| \phi \|_{L^\infty( 0,\tfrac{1}{2} T_s';W^{2,6}(\Omega) )}^{\frac{1}{4}} \| \phi \|_{H^1( 0,\tfrac{1}{2} T_s';V )}^{\frac{3}{4}}  \Big)^{-\frac{8}{3}} \Big( \frac{\delta}{10} \Big)^{\frac{8}{3}}, \nonumber \\
    & \Big(  1 + 2 C_{GN} \| \psi \|_{L^\infty( 0,\tfrac{1}{2} T_s';W^{2,6}(\Omega) )}^{\frac{1}{4}} \| \psi \|_{H^1( 0,\tfrac{1}{2} T_s';V )}^{\frac{3}{4}}  \Big)^{-\frac{8}{3}} \Big( \frac{\delta}{10} \Big)^{\frac{8}{3}}  \Big\} \label{T_s_delta}
\end{align}
only depends on $\delta$, $\Omega$, $K_0$ and the parameters of system \eqref{eq:chdf1}--\eqref{eq:chdf6}.
In this case, uniqueness can be proved in the same way as in Section \ref{ss_approximate} on the time interval $[0,T_s]$ with the help of the strict separation property.
The proof of Theorem \ref{ss_chdf_exist_3d} is complete. \qed
\smallskip 

Before ending this section, we give some comments on the result and its proof.

\begin{remark} \label{rem:variable_df_coefficients}\rm
The global existence result in Theorem \ref{ss_chdf_exist_3d}
also holds for concentration-dependent coefficients $\nu$ and $\eta$
satisfying $(\mathbf{H2})$, provided that $m_\phi,m_\psi$ are positive
constants and $r\ge3$.
Indeed, working at the approximation level and suppressing the indices,
we set
\[
M=1+\|\nabla\mu_\phi\|+\|\nabla\mu_\psi\|,
\qquad
w=|\bm{u}|^{\tfrac{r}{2}},
\qquad
D^2=\|\nabla\partial_t\phi\|^2
+\beta\|\nabla\partial_t\psi\|^2.
\]
The basic energy estimate gives
$\|\phi\|_{L^6(\Omega)}+\|\psi\|_{L^6(\Omega)}\le C$.
Testing the momentum equation with $\bm{u}$ and integrating by parts,
we obtain
\[
\nu_*\|\bm{u}\|^2+\eta_*\|\bm{u}\|_{L^r(\Omega)}^r
\le CM\|\bm{u}\|_{L^3(\Omega)}
\le CM\|\bm{u}\|_{L^r(\Omega)}.
\]
Hence $\|\bm{u}\|_{L^r(\Omega)}^{r-1}\le CM$, and thus
$\|\bm{u}\|^2+\|\bm{u}\|_{L^r(\Omega)}^r\le CM^{\frac{3}{2}}$.
Substituting this bound into \eqref{H1_sp_3d} yields
\[
\|\bm{u}\|_{\bm H^1(\Omega)}^2+\|w\|_V^2\le CM^{\frac{7}{2}}.
\]
Consequently, an interpolation gives
\[
\|\bm{u}\|_{L^{\frac{6r}{5}}(\Omega)}^r
=\|w\|_{L^\frac{12}{5}(\Omega)}^2
\le C\|w\|^\frac{3}{2}\|w\|_V^\frac{1}{2}
\le CM^\frac{9}{8}M^\frac{7}{8}
=CM^2.
\]
Using the mass balance relations and the Poincar\'e--Wirtinger inequality,
we therefore obtain
\[
\left|\int_\Omega
\partial_t\eta(\phi,\psi)|\bm{u}|^r\,\mathrm{d}x\right|
\le C\bigl(\|\partial_t\phi\|_V+\|\partial_t\psi\|_V\bigr)M^2
\le \varepsilon D^2+C_\varepsilon M^4.
\]
The term involving $\partial_t\nu$ and the convective terms admit
analogous bounds.
Since $M^2\le C(1+C_1+\mathcal E)$ and the mobility derivative terms
vanish, estimate \eqref{higher_order_estimate_uniform_3d_constant}
remains valid.
Then the same Gronwall and compactness arguments yield a global strong
solution. All constants above are independent of the approximation
parameters.
\end{remark}

\begin{remark} \label{rem:stopping_time_improvement}\rm
A stopping-time argument allows us to further relax the restriction
on $r$ for the global existence result in
Theorem \ref{ss_chdf_exist_3d}.
More precisely, global strong solutions exist for $r\ge 30/13$
if $m_\phi,m_\psi,\eta$ are positive constants and $\nu$ satisfies
$(\mathbf{H2})$, and for $r\ge 14/5$ if $m_\phi,m_\psi$ are positive
constants and both $\nu$ and $\eta$ satisfy $(\mathbf{H2})$.
In view of Remark \ref{rem:variable_df_coefficients}, it suffices
to consider the case $2<r<3$.
Let us work with the semi-Galerkin approximation, retaining the
mollification parameter $\kappa>0$ and suppressing the approximation
indices. We choose the parameters so that the initial energy and
initial chemical potentials satisfy the uniform bounds used above,
and additionally require
\[
0<\kappa\le
\min\bigl\{
1+s_{-1,n}^{\phi},\,1-s_{1,n}^{\phi},
s_{0,n}^{\psi},\,1-s_{1,n}^{\psi}
\bigr\}.
\]
Let $\tau$ be the first time at which
\[
\max\bigl\{
\|\phi(t)\|_{L^\infty(\Omega)},
\|\psi(t)\|_{L^\infty(\Omega)}
\bigr\}=2,
\]
or the maximal existence time if no such time occurs.
Since the approximate initial data lie in their physical intervals,
we find $\tau>0$.
On $[0,\tau)$, set
\begin{align*}
& M=1+\|\nabla\mu_\phi\|+\|\nabla\mu_\psi\|,
\qquad U_p=\|\bm{u}\|_{L^p(\Omega)},
\qquad 
b=\|\nabla\phi\|_{L^\infty(\Omega)}
  +\|\nabla\psi\|_{L^\infty(\Omega)},\\
& X=\|\bm{u}\|_{\bm{H}^1(\Omega)}^2
  +\bigl\||\bm{u}|^{\tfrac{r}{2}}\bigr\|_V^2,
\qquad
D^2=\|\nabla\partial_t\phi\|^2
   +\beta\|\nabla\partial_t\psi\|^2.
\end{align*}
Testing the momentum equation with $\bm{u}$ and using the stopped
$L^\infty$-bounds gives $U_r^{r-1}\le CM$.
Moreover, the elliptic regularity theory and the Gagliardo--Nirenberg inequality yield
\[
\|\phi\|_{W^{2,6}(\Omega)}
+\|\psi\|_{W^{2,6}(\Omega)}\le CM,
\qquad b\le CM^\frac{2}{3}.
\]
Estimating the coefficient gradients in $L^\infty(\Omega)$
in the derivation of \eqref{H1_sp_3d}, we therefore obtain
\[
X\le C(1+b^2)(U_2^2+U_r^r)+CM^2b^2
 \le CM^\frac{10}{3}.
\]
Consequently, we have
\[
U_3^2
\le U_r^{r-1}U_{3r}^{3-r}
\le CM^{1+\frac{10(3-r)}{3r}}
\le CM^2
\qquad\text{if }r\ge\frac{30}{13}.
\]
This controls both the convective terms and the term involving
$\partial_t\nu$ in \eqref{higher_order_estimate_uniform_3d}
by $\varepsilon D^2+C_\varepsilon M^4$, since
$U_\frac{12}{5}^2\le CU_3^2$.
If $\eta$ is nonconstant, an interpolation also gives
\[
U_\frac{6r}{5}^{2r}
=\bigl\||\bm{u}|^{\tfrac{r}{2}}\bigr\|_{L^\frac{12}{5}(\Omega)}^4
\le C U_r^{\tfrac{3r}{2}}X^\frac{1}{2}
\le CM^{\frac{3r}{2(r-1)}+\frac53}.
\]
Hence, we have
\[
\left|\int_\Omega
\partial_t\eta(\phi,\psi)|\bm{u}|^r\,\mathrm{d}x\right|
\le C(D+1)U_\frac{6r}{5}^{r}
\le\varepsilon D^2+C_\varepsilon M^4\quad \text{for}\ r\ge\frac{14}{5}.
\]
Since the mobility derivative terms vanish, the remaining terms
are estimated as in the constant-coefficient case. Thus, with
$Y=1+C_1+\mathcal E$ and $M^2\le CY$, we find
\[
Y'+cD^2\le CY^2.
\]
The basic energy estimate provides a uniform bound for the integral
of $Y$ over every subinterval of $[0,\tau)$ of length at most one.
Together with the uniform initial bound, the Gronwall and uniform
Gronwall inequalities imply
\[
\sup_{0\le t<\tau}
\bigl(
Y(t)+\|\phi(t)\|_{W^{2,6}(\Omega)}
     +\|\psi(t)\|_{W^{2,6}(\Omega)}
\bigr)\le C.
\]
All constants above are independent of the approximation parameters
and of $\tau$.
It remains to exclude a finite stopping time.
Set $I_\phi=[-1,1]$, $I_\psi=[0,1]$, and
$d_j(s)=\operatorname{dist}(s,I_j)$ for $j\in\{\phi,\psi\}$ and $s \in \mathbb{R}$.
Our choice of $\kappa$ and the normalization of the approximate
potentials imply
\[
F_{j,n,\kappa}(s)
\ge\frac n2\,\operatorname{dist}(s,I_j)^2,
\qquad s\in\mathbb R,\quad j\in\{\phi,\psi\}.
\]
The basic energy estimate therefore gives
\[
\|d_\phi(\phi(t))\|^2+\|d_\psi(\psi(t))\|^2\le Cn^{-1},
\qquad 0\le t<\tau.
\]
Since the preceding $W^{2,6}$-bounds yield uniform spatial
Lipschitz bounds, an interpolation estimate between $L^2$ and $W^{1,\infty}$ in dimension three shows that
\[
\|d_\phi(\phi(t))\|_{L^\infty(\Omega)}
+\|d_\psi(\psi(t))\|_{L^\infty(\Omega)}
\le Cn^{-\frac{1}{5}}.
\]
Thus, for all sufficiently large admissible $n$,
\[
\max\bigl\{
\|\phi(t)\|_{L^\infty(\Omega)},
\|\psi(t)\|_{L^\infty(\Omega)}
\bigr\}
\le 1+Cn^{-\frac{1}{5}}<2
\qquad\text{on }[0,\tau).
\]
For fixed approximation parameters, the chemical potential
relations also give bounded $H^3$-norms of the phase variables,
so Proposition \ref{approximate_solution_exists} allows
continuation. Consequently, a finite $\tau$ would have to be
an actual first exit, contradicting the strict bound above.
Hence $\tau=\infty$ for all sufficiently large admissible $n$,
and the global uniform estimates permit the same limiting
procedure as in the proof of Theorem \ref{ss_chdf_exist_3d}.
\end{remark}

\begin{remark}[The case $\alpha>0$] \label{relax2} \rm
For fixed $\alpha>0$, the above argument for local existence can be adapted to the time-relaxed version of system \eqref{eq:chdf1}--\eqref{eq:chdf6}, obtained by adding $\alpha\partial_t\bm{u}$ to the left-hand side of \eqref{eq:chdf1}, subject to the boundary conditions \eqref{eq:chdfb}.
An additional initial condition $\bm{u}(0)=\bm{u}_0$ is imposed, where
$\bm{u}_0\in\bm{H}^1(\Omega)\cap\bm{L}^r_\sigma(\Omega)$.
The resulting local solution $\bm{u}$ satisfies 
\begin{align}
    \bm{u} \in L^\infty\big(  0,T_\alpha;\bm{H}^1(\Omega)\cap\bm{L}^r_\sigma(\Omega)  \big)
    \cap H^1\big( 0,T_\alpha;\bm{L}^2(\Omega) \big), \quad 
    | \bm{u} |^\frac{r}{2} \in L^2 \big(  0,T_\alpha;H^1(\Omega)  \big),
    \label{alpha-u}
\end{align}
where $T_\alpha$ is a sufficiently short time depending on $\alpha$. 

The proof needs to be modified accordingly. In deriving the higher-order estimate \eqref{higher_order_estimate_approximate}, one encounters the convective term
\begin{equation} \label{convective_phi}
    \int_\Omega \partial_t\mu_\phi\,\bm{u}\cdot\nabla\phi\,\mathrm{d}x,
\end{equation}
and an analogous term involving $\psi$.
For the system without time relaxation, these terms are canceled using the argument in \cite{G2020}.
When $\alpha>0$, however, the calculation leading to \eqref{higher_order_estimate_sp} produces the additional term
$\langle\alpha\partial_t^2\bm{u},\bm{u}\rangle$,
which cannot be controlled by the estimates used for $\alpha=0$. Instead, we test the momentum equation directly with $\partial_t\bm{u}$, retaining the Darcy and Forchheimer contributions as part of the higher-order energy.
This gives
\begin{align}
    & \frac{\mathrm{d}}{\mathrm{d}t}
    \int_\Omega
    \frac{1}{2}\nu(\phi,\psi)| \bm{u} |^2
    + \frac{1}{r}\eta(\phi,\psi)| \bm{u} |^r
    \,\mathrm{d}x
    + \alpha\| \partial_t\bm{u} \|^2
    \notag 
    \\
    & \quad =
    ( 
        \mu_\phi\nabla\phi+\mu_\psi\nabla\psi
    , \partial_t\bm{u} )
    + \int_\Omega
    \frac{1}{2}\partial_t\nu(\phi,\psi)| \bm{u} |^2
    + \frac{1}{r}\partial_t\eta(\phi,\psi)| \bm{u} |^r
    \,\mathrm{d}x.
    \label{high-alpha1}
\end{align}
Next, we write the convective term \eqref{convective_phi} as
\[
    \int_\Omega \partial_t\mu_\phi\,\bm{u}\cdot\nabla\phi
    \,\mathrm{d}x
    =
    - \int_\Omega \mu_\phi\,\partial_t\bm{u}\cdot\nabla\phi
    \,\mathrm{d}x
    - \int_\Omega \mu_\phi\,\bm{u}\cdot\nabla\partial_t\phi
    \,\mathrm{d}x
     + \frac{\mathrm{d}}{\mathrm{d}t}
    \Big(
        \int_\Omega \mu_\phi\,\bm{u}\cdot\nabla\phi
        \,\mathrm{d}x
    \Big).
\]
For fixed $\alpha>0$, the first term is estimated using a small fraction of the dissipation
$\alpha\| \partial_t\bm{u} \|^2$ in \eqref{high-alpha1}. The total time derivative is incorporated into the higher-order energy. Moreover, we have 
\begin{align*}
    \Big| 
        \int_\Omega \mu_\phi\,\bm{u}\cdot\nabla\phi
        \,\mathrm{d}x
     \Big|
    &=
    \Big| 
        \int_\Omega \phi\,\nabla\mu_\phi\cdot\bm{u}
        \,\mathrm{d}x
     \Big|
    \\
    &\le
    \| \phi-\overline{\phi} \|_{L^\infty(\Omega)}
    \| \nabla\mu_\phi \|
    \| \bm{u} \|
    \\
    &\le
    C\| \nabla\phi \|^{\frac12}
    \| \Delta\phi \|^{\frac12}
    \| \nabla\mu_\phi \|
    \| \bm{u} \|
    \\
    &\le
    C\| \nabla\phi \|^{\frac12}
    \big( 
        \| \nabla\phi \|^{\frac12}
        \| \nabla\mu_\phi \|^{\frac12}
        + \| \nabla\phi \|
        + \| \nabla\psi \|
     \big)^{\frac12}
    \| \nabla\mu_\phi \|
    \| \bm{u} \|
    \\
    &\le
    \frac{m_{\phi,*}}{4}\| \nabla\mu_\phi \|^2
    + C\big( 
        \| \nabla\phi \|^2
        + \| \nabla\psi \|^2
     \big)
    \big( 
        \| \bm{u} \|^{\frac83}
        + \| \bm{u} \|^2
     \big).
\end{align*}
The basic energy estimate bounds the last term by a constant depending on $\alpha>0$, the initial data, and the time interval under consideration.
The analogous cross term involving $\psi$ is treated in the same way.
Thus, after adding a sufficiently large constant, the modified higher-order energy controls the gradient of the chemical potentials and the velocity. The spatial velocity estimate also depends on $\alpha$ when $\alpha>0$. 
At the Galerkin level, testing the momentum equation with $\bm{A}\bm{u}$ produces the additional term
\[
    \alpha( \partial_t\bm{u}, \bm{A}\bm{u} )
    =
    \frac{\alpha}{2}\frac{\mathrm{d}}{\mathrm{d}t}
    \| \mathrm{curl} \,\bm{u} \|^2.
\]
The remaining terms are treated as in the derivation of \eqref{H1_sp_3d}.
Combining this estimate with the preceding energy estimates yields the regularity property \eqref{alpha-u}. 
We note that the aforementioned estimates and the local existence time depend on $\alpha$. Hence, the above argument does not directly yield the corresponding existence result for $\alpha=0$ (by taking $\alpha \to 0$).
\end{remark}

%%%%%%%%%%%%%%%%%%%%%%%%%%%%%%%%
\section{Proof of Theorem \ref{ss_chdf_exist_global_small}}
\label{exist_global_small}
\setcounter{equation}{0} 

We first establish local existence of strong solutions satisfying the strict separation property.
For brevity, let us set
\begin{align*}
    P_{\delta,1} & \stackrel{\mathrm{def}}{=} \max \{  \| F_\phi' \|_{L^\infty(-1 + \delta,1 - \delta)}, \| F_\psi' \|_{L^\infty(\delta,1 - \delta)}  \}, \\
    P_{\delta,2} & \stackrel{\mathrm{def}}{=} \max \{  \| F_\phi'' \|_{L^\infty(-1 + \delta,1 - \delta)}, \| F_\psi'' \|_{L^\infty(\delta,1 - \delta)}  \}.
\end{align*}
Then, the proof of the local well-posedness part of Theorem \ref{ss_chdf_exist_3d} gives the following result.
\begin{corollary} \label{ss_local_strong_separate_bounded}
Assume that \textbf{(H1)} and \textbf{(H2)} hold and that $F_\phi\in C^3((-1,1))$, $F_\psi\in C^3((0,1))$, and $G\in C^3((-1,1)\times(0,1))$.
Suppose $\phi_0, \psi_0 \in W \cap H^3(\Omega)$ with 
\[
\| \phi_0 \|_{L^\infty(\Omega)} \le 1 - \delta, \quad \big\| \psi_0 - \frac{1}{2} \big\|_{L^\infty(\Omega)} \le \frac{1}{2} - \delta,
\]
for some $\delta \in (0,\tfrac{1}{2})$.
Suppose that
\[
\| \phi_0 \|_V \le M_1, \quad \| \psi_0 \|_V \le M_1, \quad \| \phi_0 \|_{H^3(\Omega)} \le M_2, \quad \| \psi_0 \|_{H^3(\Omega)} \le M_2,
\]
for some $M_1,M_2 >0$ and set
\begin{align*}
    & E_1 \stackrel{\mathrm{def}}{=} \frac{1}{2} M_1^2 + \frac{\beta}{2} M_1^2 + | \Omega |(\| F_\phi \|_{C([-1,1])} + \| F_\psi \|_{C([0,1])} + \| G \|_{C(\mathbb{R}^2)} + 2 | \sigma_2 | C_P^2), \\
    & E_{2,\delta} \stackrel{\mathrm{def}}{=} \big(  \frac{9}{4} m_\phi^* + \frac{3\nu^*}{\nu_*^2} + \frac{3\eta^*}{\nu_* \eta_*}  \big) \big(  2 M_2^2 + 2 P_{\delta,2}^2 M_1^2  \big) + \big(  m_\psi^* + \frac{2 \nu^*}{\nu_*^2} + \frac{2 \eta^*}{\nu_* \eta_*}  \big) \big(  2 \beta^2 M_2^2 + 2 P_{\delta,2}^2 M_1^2  \big) \\
    & \qquad \quad + \Big(  3 \sigma_2^2 C_P^4 \big(  \frac{3}{4} m_\phi^* + \frac{\nu^*}{\nu_*^2} + \frac{\eta^*}{\nu_* \eta_*}  \big) + \| G \|_{C^2(\mathbb{R}^2)}^2 \big( 9 m_\phi^* + 4 m_\psi^* + \frac{20 \nu^*}{\nu_*^2} + \frac{20 \eta^*}{\nu_* \eta_*}  \big) \Big) (1 + M_1^2) \\
    & \qquad \quad + 4 C_{\mathrm{source}} \sigma_1  (1 + M_1^2) + 4 C_{\mathrm{source}} \sigma_1 e^{4 C_{\mathrm{source}}} \big(  E_1 +  | \Omega | \big(  \| G \|_{C(\mathbb{R}^2)} + 2 | \sigma_2 | C_P^2  \big) + 1   \big) + 1.
\end{align*}
There exists $T_1>0$ such that problem \eqref{eq:chdf1}--\eqref{eq:chdfi} admits a unique strong solution $( \bm{u},\phi,\psi,\mu_\phi,\mu_\psi )$ on $[0,T_1]$ in the sense of Definition \ref{ss_chdf_def}.
Moreover, it holds
\begin{align*}
    & \widetilde{\mathcal{E}}(\bm{u}(t),\phi(t),\psi(t),\mu_\phi(t),\mu_\psi(t)) \le 2 E_{2,\delta}, \quad \forall \ t \in [0,T_1], \\
    & \| \phi(t) \|_{L^\infty(\Omega)} \le 1 - \frac{9}{10} \delta, \quad \big\| \psi(t) - \frac{1}{2} \big\|_{L^\infty(\Omega)} \le \frac{1}{2} - \frac{9}{10} \delta, \quad \forall \ t \in [0,T_1].
\end{align*}
Here, $T_1$ and $\widetilde{\mathcal{E}}$, the modified version of the higher-order energy defined in \eqref{higher_order_energy}, are given by
\begin{align}
    & \hypertarget{modified_higher_order_energy}{\widetilde{\mathcal{E}}}(\bm{u},\phi,\psi,\mu_\phi,\mu_\psi) \stackrel{\mathrm{def}}{=} \mathcal{E}(\bm{u},\phi,\psi,\mu_\phi,\mu_\psi) \\
    & \qquad \qquad \qquad \qquad \quad + 2C_{\mathrm{source}} e^{4 C_{\mathrm{source}}} \big(  E_1 +  | \Omega | \big(  \| G \|_{C(\mathbb{R}^2)} + 2 | \sigma_2 | C_P^2  \big) + 1  \big)  \sigma_1 \big| \overline{\phi}-c \big|, \\
    & \hypertarget{T_1_local_existence}{T_1} \stackrel{\mathrm{def}}{=} \min \Big\{  1, \hyperlink{T_s_delta}{T_s}(E_1,E_{2,\delta},\delta), \frac{8 \sqrt{2} - 1}{ 28 \sqrt{2} \hyperlink{C_2_energy_higher_order}{C_2} E_{2,\delta}^{\frac{7}{2}} }  \Big\},
\end{align}
where $C_2 > 0$ is the coefficient in the estimate \eqref{energy_estimate_local_higher_order}, and $T_s$ is defined by \eqref{T_s_delta}.
\end{corollary}

Next, we establish the existence of a global strong solution for initial data close to a local energy minimizer by iteratively extending the local solution in time. 

As in \cite{HW2024}, using the direct method of the calculus of variations and a dynamic approach, we obtain the following result.
\begin{proposition} \label{minimizer_existence_regularity}
Suppose that (\textbf{H1}) is satisfied and that $\Omega$ is a bounded domain in $\mathbb{R}^d$ with boundary $\partial \Omega$ of class $C^3$, $d\in \{ 2,3 \}$.
Let $m_1 \in (-1,1)$ and $m_2 \in (0,1)$. We consider the following minimization problem:
\begin{equation} \label{eq:minimize_free_energy}
    \mathrm{Find} \ (\phi_*,\psi_*) \in \mathcal{Z}_{m_1,m_2} \quad \mathrm{such} \ \ \mathrm{that} \quad E_{\mathrm{free}}(\phi_*,\psi_*)= \inf\limits_{(\phi,\psi) \in \mathcal{Z}_{m_1,m_2}} E_{\mathrm{free}}(\phi,\psi).
\end{equation}
\begin{itemize}
    \item[(1)] Problem \eqref{eq:minimize_free_energy} admits at least one (global) minimizer $(\phi_*,\psi_*) \in \mathcal{Z}_{m_1,m_2}$.
    \item[(2)] Let $(\phi_*,\psi_*) \in \mathcal{Z}_{m_1,m_2}$ be a (local or global) energy minimizer. 
    Then it belongs to $W \times W$ and there exists a positive constant $\delta_* \in (0,\tfrac{1}{2})$ such that
    \[
    \| \phi_* \|_{L^\infty(\Omega)} \le 1 - \delta_*, \quad \big\| \psi_*-\frac{1}{2} \big\|_{L^\infty(\Omega)} \le \frac{1}{2} - \delta_*.
    \]
    Moreover, $(\phi_*,\psi_*)$ is a solution to the elliptic problem
\begin{align} 
\begin{cases}
    - \Delta \phi_* + \sigma_2 \mathcal{N}\big(  \phi_* - \overline{\phi_*}  \big) + F_\phi'(\phi_*) + \partial_\phi G(\phi_*,\psi_*) = \mu_{\phi,\infty},  \quad \text{a.e. in }\Omega,  
    \\
    - \beta \Delta \psi_* + F_\psi'(\psi_*) + \partial_\psi G(\phi_*,\psi_*) = \mu_{\psi,\infty}, 
    \quad \text{a.e. in }\Omega, 
    \\
    \partial_{\bm{n}} \phi_* = \partial_{\bm{n}} \psi_* = 0,\quad \text{a.e. on } \partial \Omega,
\end{cases}
\label{ch_stationary_P}
\end{align}
where $\mu_{\phi,\infty}, \mu_{\psi,\infty} \in \mathbb{R}$ are constants.
\end{itemize}
\end{proposition}

We also need the following version of the {\L}ojasiewicz--Simon inequality for the coupled elliptic system \eqref{ch_stationary_P} (see \cite{O2025} for a proof).
\begin{lemma}[{\L}ojasiewicz--Simon inequality] \label{ls_inequality}
Let $\Omega \subset \mathbb{R}^d$, $d\in\{2,3\}$, be a bounded domain with a smooth boundary $\partial\Omega$.
Suppose that Assumption (\textbf{H1}) is satisfied, $F_\phi$ is real analytic in $(-1,1)$, $F_\psi$ is real analytic in $(0,1)$, and $G$ is real analytic in $(-1,1) \times (0,1)$.
For any $m_1\in(-1,1)$ and $m_2\in(0,1)$ and any solution $(\phi_*,\psi_*)\in(W\times W)\cap\mathcal{Z}_{m_1,m_2}$ to problem \eqref{ch_stationary_P},
there exist constants $\theta \in (0,\frac12)$ and $C_{\mathrm{LS}},l > 0$ such that for $\widetilde{\phi},\widetilde{\psi} \in V\cap L^\infty(\Omega)$ satisfying
\[
\| \widetilde{\phi}-\phi_* \|_{V} \leq l, \quad \| \widetilde{\psi}-\psi_* \|_{V} \leq l,\quad
\| \widetilde{\phi} \|_{L^{\infty}(\Omega)} \leq 1 -\frac{3}{4}\delta_*,\quad  \big\| \widetilde{\psi}-\frac{1}{2} \big\|_{L^{\infty}(\Omega)} \leq \frac{1}{2} -\frac{3}{4} \delta_*,\]
where $\delta_*\in (0,1)$ is the constant determined in Proposition \ref{minimizer_existence_regularity},
the following inequality holds
\begin{align}
   &  \big|E_{\mathrm{free}}(\widetilde{\phi},\widetilde{\psi})-E_{\mathrm{free}}(\phi_*,\psi_*)\big|^{1-\theta}\notag \\
   &\quad  \le C_{\mathrm{LS}} \left( \|\widetilde{\mu}_\phi-\overline{\widetilde{\mu}_\phi}\|_{V^*_{(0)}}+ \|\widetilde{\mu}_\psi-\overline{\widetilde{\mu}_\psi}\|_{V^*_{(0)}} +  \left(|\overline{\widetilde{\phi}}-\overline{\phi_*}| + |\overline{\widetilde{\psi}}-\overline{\psi_*}|\right)^{1-\theta}\right), \label{eq:ls_inequality}
\end{align}
with $\widetilde{\mu}_\phi=-\Delta \widetilde{\phi} + F_\phi'(\widetilde{\phi}) + \sigma_2 \mathcal{N}(\widetilde{\phi} -\overline{\widetilde{\phi}}) + \partial_\phi G(\widetilde{\phi},\widetilde{\psi})$
and $\widetilde{\mu}_\psi=-\beta \Delta \widetilde{\psi} + F_\psi'(\widetilde{\psi}) + \partial_\psi G(\widetilde{\phi},\widetilde{\psi})$.
\end{lemma}

In what follows, we refine the estimates for strong solutions whose initial data lie near a local energy minimizer $(\phi_*,\psi_*)$. 
Firstly, note that if there exists some $\delta>0$ such that
\[
\| \phi_i \|_{L^\infty(\Omega)} \le 1 - \delta, \quad \Big\| \psi_i-\frac{1}{2} \Big\|_{L^\infty(\Omega)} \le \frac{1}{2} - \delta, \quad i \in \{ 1,2 \},
\]
then it holds
\begin{align}
    & |  E_{\mathrm{free}}(\phi_1,\psi_1) - E_{\mathrm{free}}(\phi_2,\psi_2)  | \nonumber \\
    & \quad \le \frac{1}{2} \| \nabla \phi_1 + \nabla \phi_2 \| \| \nabla \phi_1 - \nabla \phi_2 \| + \frac{\beta}{2} \| \nabla \psi_1 + \nabla \psi_2 \| \| \nabla \psi_1 - \nabla \psi_2 \| \nonumber \\
    & \qquad + | \Omega |^{\frac{1}{2}} (  P_{\delta,1} + 2 \| G \|_{C^1(\mathbb{R}^2)}  ) \big(  \| \phi_1 - \phi_2 \| + \| \psi_1 - \psi_2 \|  \big) \nonumber \\
    & \qquad + \frac{| \sigma_2 |}{2} \big\|  \nabla \mathcal{N} \big(  \phi_1 + \phi_2 - \overline{\phi_1} - \overline{\phi_2}  \big)  \big\| \big\|  \nabla \mathcal{N} \big(  \phi_1 - \phi_2 - \overline{\phi_1} + \overline{\phi_2}  \big)  \big\| \nonumber \\
    & \quad \le \hypertarget{C_free}{C_{\mathrm{free},\delta}} \big(  \| \phi_1 - \phi_2 \|_V + \| \psi_1 - \psi_2 \|_V  \big). \label{energy_difference_estimate}
\end{align}
Here $C_{\mathrm{free},\delta}$ also depends on a common bound for the $V$-norms of $\phi_i$ and $\psi_i$, $i=1,2$.
Then, recall the Sobolev embedding
\[
\| f \|_{C(\overline{\Omega})} \le C_S \| f \|_W,
\]
where $C_S>0$ only depends on $\Omega$.
Assume that \textbf{(H1)} and \textbf{(H2)} hold, and that $F_\phi$ is real analytic in $(-1,1)$, $F_\psi$ is real analytic in $(0,1)$, and $G$ is real analytic in $(-1,1) \times (0,1)$. 
Let $(\phi_*,\psi_*)$ be a given local minimizer of the free energy $E_{\mathrm{free}}$ in $\mathcal{Z}_{m_1,m_2}$ in the sense of Definition \ref{minimizer_def}, with $\sigma_1(m_1 - c)=0$, and $M_2 > 0$ be the prescribed upper bound for the higher-order norm in Corollary \ref{ss_local_strong_separate_bounded}. 
Let $\lambda_*$ corresponding to $(\phi_*,\psi_*)$ be chosen so that for all $(\phi,\psi) \in \mathcal{Z}_{m_1,m_2}$ with $\| \phi-\phi_* \|_V + \| \psi-\psi_* \|_V \le \lambda_*$, it holds
\[
E_{\mathrm{free}}(\phi_*,\psi_*) \le E_{\mathrm{free}}(\phi,\psi).
\]
Let $C_{\mathrm{LS},*}$, $l_*$ and $\theta_*$ be the constants corresponding to $(\phi_*,\psi_*)$ in Lemma \ref{ls_inequality} and $\delta_*$ be the separation constant given by Proposition \ref{minimizer_existence_regularity}.
Let $\epsilon > 0$ be given as in Theorem \ref{ss_chdf_exist_global_small}.
Taking the bound $M_1$ in Corollary \ref{ss_local_strong_separate_bounded} as
\[
M_1 = \max \{  \| \phi_* \|_V + 1, \| \psi_* \|_V + 1  \},
\]
we define
\[
\omega \stackrel{\mathrm{def}}{=} \min \Big\{  1, \epsilon, \frac{\lambda_*}{2}, l_*, \frac{\delta_*}{10(1+C_S)}, \frac{L_{\mathrm{coef}} E_3 T_1}{14 (C_{\mathrm{free},\frac{3}{4}\delta_*} + 1)}  \Big\},
\]
where 
\[
E_3 \stackrel{\mathrm{def}}{=} \max \{ 1, 6 C_{\mathrm{source}} \sigma_1 \}, \quad L_{\mathrm{coef}} \stackrel{\mathrm{def}}{=} \frac{ \min \{  \tfrac{m_{\phi,*}}{2}, m_{\psi,*}, \nu_*, \eta_* \} }{ \max \{  m_\phi^*, m_\psi^*, \nu^*, \eta^*  \} },
\]
and $\hyperlink{T_1_local_existence}{T_1}$ is the length of the existence interval defined in Corollary \ref{ss_local_strong_separate_bounded}.
Defining $E_{2,\delta}$ and $\widetilde{\mathcal{E}}$ as in Corollary \ref{ss_local_strong_separate_bounded}, we have the following uniform-in-time higher-order estimate of strong solutions:
\begin{proposition} \label{uniform_control_small_solution}
    There exist two positive constants $\rho_\phi,\rho_\psi \in (0,\tfrac{1}{3} \omega]$ depending on $\epsilon$, $M_2$, $\delta_*$, $\phi_*$, $\psi_*$, $m_\phi$, $m_\psi$, $\nu$, $\eta$, $P_{\tfrac{3}{4} \delta_*,1}$, and $P_{\tfrac{3}{4} \delta_*,2}$ such that, for any $\phi_0$ and $ \psi_0$ satisfying
    \begin{align*}
        & \| \phi_0 \|_{H^3(\Omega)} \le M_2, \quad \| \psi_0 \|_{H^3(\Omega)} \le M_2, \\
        & \| \phi_0-\phi_* \|_W \le \rho_\phi, \quad \| \psi_0-\psi_* \|_W \le \rho_\psi, \\
        & \overline{\psi_0} = m_2, \\
        & \overline{\phi_0} = m_1, \ \text{if } \sigma_1 = 0,
    \end{align*}
    and any $T \in [T_1,\infty)$, if $( \bm{u},\phi,\psi,\mu_\phi,\mu_\psi )$ is the strong solution to problem \eqref{eq:chdf1}--\eqref{eq:chdfi} in $[0,T]$ originating from $(\phi_0,\psi_0)$ with
    \begin{align*}
        & | \phi(t) | \le 1 - \frac{3}{4} \delta_*, \quad \Big| \psi(t) - \frac{1}{2} \Big| \le \frac{1}{2} - \frac{3}{4} \delta_*, \quad \forall t \in [0,T], \\
        & \widetilde{\mathcal{E}}(\bm{u}(t),\phi(t),\psi(t),\mu_\phi(t),\mu_\psi(t)) \le 2 E_{2,\frac{9}{10} \delta_*}, \quad \forall t \in [0,T],
    \end{align*}
    then it holds
    \begin{align*}
        & \| \phi(t) - \phi_* \|_W + \| \psi(t) - \psi_* \|_W \le \omega, \quad \forall t \in [0,T], \\
        & \int_0^T \widetilde{\mathcal{E}}(\bm{u},\phi,\psi,\mu_\phi,\mu_\psi) \ \mathrm{d} t \ \le \frac{1}{3} E_3 T_1.
    \end{align*}
\end{proposition}

\noindent \textbf{Proof.} We prove this by a contradiction argument like that in \cite[Section 5.2]{HW2024}.
Define
\[
\widehat{T} \stackrel{\mathrm{def}}{=} \sup \big\{  t \in [0,T] \ | \ \forall \tau \in [0,t], \ \| \phi(\tau) - \phi_* \|_W + \| \psi(\tau) - \psi_* \|_W \le \omega  \big\}.
\]
By continuity, we have $\widehat{T}>0$.
Suppose that for all $\rho_\phi, \rho_\psi \in (0,\tfrac{1}{3} \omega)$ there exist $\phi_0$ and $\psi_0$ satisfying
 \begin{align*}
    & \| \phi_0 \|_{H^3(\Omega)} \le M_2, \quad \| \psi_0 \|_{H^3(\Omega)} \le M_2, \\
    & \| \phi_0-\phi_* \|_W \le \rho_\phi, \quad \| \psi_0-\psi_* \|_W \le \rho_\psi, \\
    & \overline{\psi_0} = m_2, \\
    & \overline{\phi_0} = m_1, \ \text{if } \sigma_1 = 0,
\end{align*}
such that $\widehat{T} < T$.
Firstly, recall the differential form of the energy identity 
\begin{align}
    & \frac{\mathrm{d}}{\mathrm{d} t} E_{\mathrm{free}} (\phi(t),\psi(t)) \nonumber \\
    & \qquad + \int_\Omega m_\phi(\phi(t)) | \nabla \mu_\phi |^2 + m_\psi(\psi(t)) | \nabla \mu_\psi |^2 + \nu(\phi(t),\psi(t)) | \bm{u}(t) |^2 + \eta(\phi(t),\psi(t)) | \bm{u}(t) |^r \ \mathrm{d} x \nonumber \\
    & \qquad + \int_\Omega \sigma_1 \big( \overline{\phi}(t)-c \big) \mu_\phi(t) \ \mathrm{d} x = 0. \label{energy_identity} 
\end{align}
Therefore, it holds (cf. \eqref{basic_energy_estimate_approximate})
\begin{align*}
    & \frac{\mathrm{d}}{\mathrm{d} t} E_{\mathrm{free}} (\phi(t),\psi(t)) \nonumber \\
    & \qquad + \int_\Omega \frac{1}{2} m_\phi(\phi(t)) | \nabla \mu_\phi |^2 + m_\psi(\psi(t)) | \nabla \mu_\psi |^2 + \nu(\phi(t),\psi(t)) | \bm{u}(t) |^2 + \eta(\phi(t),\psi(t)) | \bm{u}(t) |^r \ \mathrm{d} x \nonumber \\
    & \quad \le C_{\mathrm{source}}' \sigma_1 e^{-\sigma_1 t} \big| \overline{\phi_0} - c \big| .
\end{align*}
Let $\widetilde{\phi}(t) = \phi(t) - \overline{\phi}(t) + \overline{\phi_*}$. 
Then we have
\begin{align*}
    E_{\mathrm{free}} (\phi(t),\psi(t)) & = E_{\mathrm{free}} (\widetilde{\phi}(t),\psi(t)) + \int_\Omega F_\phi (\phi(t)) + G(\phi(t),\psi(t)) - F_\phi (\widetilde{\phi}(t)) - G(\widetilde{\phi}(t),\psi(t)) \ \mathrm{d} x \nonumber \\
    & \ge E_{\mathrm{free}} (\widetilde{\phi}(t),\psi(t)) - | \Omega | (P_{\frac{2\delta_*}{3},1} + \| G \|_{C^1(\mathbb{R}^2)}) \big|  \overline{\phi} (t) - \overline{\phi_*}  \big| \nonumber \\
    & = E_{\mathrm{free}} (\widetilde{\phi}(t),\psi(t)) - C_{\mathrm{free},\delta_*}' e^{-\sigma_1 t} \big|  \overline{\phi_0} - \overline{\phi_*}  \big|.
\end{align*}
Setting
\[
\widetilde{E}_{\mathrm{free}}(t) \stackrel{\mathrm{def}}{=} E_{\mathrm{free}}(\phi(t),\psi(t)) + (2C_{\mathrm{source}}' + C_{\mathrm{free},\delta_*}') e^{-\sigma_1 t} \big|  \overline{\phi_0} - \overline{\phi_*}  \big|,
\]
we deduce
\begin{align}
    - \frac{\mathrm{d}}{\mathrm{d} t} \widetilde{E}_{\mathrm{free}}(t) & \ge \frac{m_{\phi,*}}{2} \| \nabla \mu_\phi \|^2 + m_{\psi,*} \| \nabla \mu_\psi \|^2 + \nu_* \| \bm{u} \|^2 + \eta_* \| \bm{u} \|^r_{\bm{L}^r(\Omega)} \nonumber \\
    & \quad + (C_{\mathrm{source}}' + C_{\mathrm{free},\delta_*}') \sigma_1 e^{-\sigma_1 t} \big|  \overline{\phi_0} - \overline{\phi_*}  \big| \nonumber \\
    & \ge L_{\mathrm{coef}}' \big(  \| \nabla \mu_\phi \|^2 + \| \nabla \mu_\psi \|^2 + \| \bm{u} \|^2 + \| \bm{u} \|^r_{\bm{L}^r(\Omega)}  \big) \nonumber \\
    & \quad + (C_{\mathrm{source}}' + C_{\mathrm{free},\delta_*}') \sigma_1 e^{-\sigma_1 t} \big|  \overline{\phi_0} - \overline{\phi_*}  \big| \nonumber \\
    & \ge L_{\mathrm{coef}} \big(  m_\phi^* \| \nabla \mu_\phi \|^2 + m_\psi^* \| \nabla \mu_\psi \|^2 + \nu^* \| \bm{u} \|^2 + \eta^* \| \bm{u} \|^r_{\bm{L}^r(\Omega)}  \big) \nonumber \\
    & \quad + (C_{\mathrm{source}}' + C_{\mathrm{free},\delta_*}') \sigma_1 e^{-\sigma_1 t} \big|  \overline{\phi_0} - \overline{\phi_*}  \big|, \quad \forall t \in [0,T], \label{energy_inequality_modified}
\end{align}
where $L_{\mathrm{coef}}' \stackrel{\mathrm{def}}{=} \min \{  \tfrac{m_{\phi,*}}{2}, m_{\psi,*}, \nu_*, \eta_* \} $.
We claim that 
\[
\widetilde{E}_{\mathrm{free}} (t) > E_{\mathrm{free}} (\phi_*,\psi_*),
\]
for all $t \in [0,\widehat{T}]$.
Otherwise, since $(\phi_*,\psi_*)$ is a local minimizer in $\mathcal{Z}_{m_1,m_2}$ and it holds
\begin{align*}
    & \| \widetilde{\phi}(t)-\phi_* \|_W + \| \psi(t)-\psi_* \|_W \\
    & \quad \le \| \phi(t)-\phi_* \|_W + \| \psi(t)-\psi_* \|_W \\
    & \quad \le \omega < \lambda_*, \ \forall t \in [0,\widehat{T}],
\end{align*}
there exists some $t_0 \in [0,\widehat{T}]$ such that
\[
\widetilde{E}_{\mathrm{free}} (t_0) = E_{\mathrm{free}} (\phi_*,\psi_*).
\]
In addition, by the time continuity of $\phi$ and $\psi$ there exists some $\epsilon_T \in (0, T - \widehat{T}]$ such that 
\[
    \| \widetilde{\phi}(t)-\phi_* \|_W + \| \psi(t)-\psi_* \|_W < \lambda_*, \ \forall\, t \in [0,\widehat{T}+\epsilon_T].
\]
Consequently, it follows from \eqref{energy_inequality_modified} that
\begin{align*}
    & \widetilde{E}_{\mathrm{free}} (t) = E_{\mathrm{free}} (\phi_*,\psi_*), \\
    & \| \nabla \mu_\phi(t) \| = \| \nabla \mu_\psi(t) \| = \| \bm{u}(t) \| = \sigma_1 \big| \overline{\phi}(t)-c \big| = 0,
\end{align*}
for all $t \in [t_0,\widehat{T}+ \epsilon_T]$, which implies that $(\bm{u},\phi,\psi,\mu_\phi,\mu_\psi)$ is stationary on $[t_0,\widehat{T}+\epsilon_T]$.
Therefore, the uniqueness of strong solutions to the problem \eqref{eq:chdf1}--\eqref{eq:chdfi} entails that $(\phi(t),\psi(t)) \equiv (\phi(t_0),\psi(t_0))$ for all $t \in [t_0,T]$, which contradicts the assumption $\widehat{T} < T$.
Therefore, under this assumption, the above claim holds.
This enables the application of the {\L}ojasiewicz--Simon inequality on $[0,\widehat{T}]$.
It follows from Lemma \ref{ls_inequality} that
\begin{align}
    & - \frac{\mathrm{d}}{\mathrm{d} t} \big(  \widetilde{E}_{\mathrm{free}}(t) - E_{\mathrm{free}} (\phi_*,\psi_*)  \big)^{\theta_*} \nonumber \\
    & \quad = - \theta_* \big(  \widetilde{E}_{\mathrm{free}}(t) - E_{\mathrm{free}} (\phi_*,\psi_*)  \big)^{\theta_* - 1} \frac{\mathrm{d}}{\mathrm{d} t} \widetilde{E}_{\mathrm{free}} (t) \nonumber \\
    & \quad \ge \frac{ \theta_* \Big(  L_{\mathrm{coef}}' \big(   \| \nabla \mu_\phi \|^2 + \| \nabla \mu_\psi \|^2 + \| \bm{u} \|^2 + \| \bm{u} \|^r_{\bm{L}^r(\Omega)}  \big) + (C_{\mathrm{source}}' + C_{\mathrm{free},\delta_*}') \sigma_1 e^{-\sigma_1 t} \big|  \overline{\phi_0} - \overline{\phi_*}  \big|  \Big) }{ \Big( E_{\mathrm{free}} (\phi,\psi) - E_{\mathrm{free}} (\phi_*,\psi_*) + (2 C_{\mathrm{source}}' + C_{\mathrm{free},\delta_*}') e^{-\sigma_1 t} \big|  \overline{\phi_0} - \overline{\phi_*}  \big|  \Big)^{1-\theta_*} } \nonumber \\
    & \quad \ge \frac{ C \Big(  \| \nabla \mu_\phi \|^2 + \| \nabla \mu_\psi \|^2 + \| \bm{u} \|^2 + \| \bm{u} \|^r_{\bm{L}^r(\Omega)} + \sigma_1 e^{-\sigma_1 t} \big|  \overline{\phi_0} - \overline{\phi_*}  \big|  \Big) }{ \Big(  \| \bm{u} \| + \| \bm{u} \|^{\frac{r}{2}}_{\bm{L}^r(\Omega)} + \| \nabla \mu_\phi \| + \| \nabla \mu_\psi \| + e^{-\sigma_1(1-\theta_*) t} \big|  \overline{\phi_0} - \overline{\phi_*}  \big|^{1-\theta_*}  \Big) } \nonumber \\
    & \quad \ge C_{\mathrm{diss}} \Big(  \| \bm{u} \| + \| \bm{u} \|^{\frac{r}{2}}_{\bm{L}^r(\Omega)} + \| \nabla \mu_\phi \| + \| \nabla \mu_\psi \| + e^{-\frac{\sigma_1}{2} t} \big|  \overline{\phi_0} - \overline{\phi_*}  \big|^{\frac{1}{2}}  \Big), \\
    & C_{\mathrm{diss}} \stackrel{\mathrm{def}}{=} \frac{\theta_* \min \big\{  L_{\mathrm{coef}}', (C_{\mathrm{source}}' + C_{\mathrm{free},\delta_*}')  \big\} \min \{  1, \max \{ \sigma_1, 1 - \mathrm{sgn}(\sigma_1) \}  \} }{ 10 \max \big\{  1, C_{\mathrm{LS},*} + (2 C_{\mathrm{source}}' + C_{\mathrm{free},\delta_*}')^{1-\theta_*}  \big\} } .
\end{align}
Thus, by comparison, we have
\begin{align*}
    & \|  \partial_t \phi - \overline{\partial_t \phi}  \|_{V^*} + \|  \partial_t \psi  \|_{V^*} \\
    & \quad \le \| \bm{u} \| \big(  \| \phi \|_{L^\infty(\Omega)} + \| \psi \|_{L^\infty(\Omega)}  \big) + m_\phi^* \| \nabla \mu_\phi \| + m_\psi^* \| \nabla \mu_\psi \| \\
    & \quad \le - \frac{\max \{ 2, m_\phi^*, m_\psi^* \}}{C_{\mathrm{diss}}} \frac{\mathrm{d}}{\mathrm{d} t} \big(  \widetilde{E}_{\mathrm{free}} - E_{\mathrm{free}}(\phi_*,\psi_*)  \big)^{\theta_*}.
\end{align*}
So it holds
\begin{align*}
    & \int_0^{\widehat{T}} \|  \partial_t \phi  \|_{V^*} + \|  \partial_t \psi  \|_{V^*} \ \mathrm{d} t \\
    & \quad \le \int_0^{\widehat{T}} \|  \overline{\partial_t \phi}  \|_{V^*} - \frac{\max \{ 2, m_\phi^*, m_\psi^* \}}{C_{\mathrm{diss}}} \frac{\mathrm{d}}{\mathrm{d} t} \big(  \widetilde{E}_{\mathrm{free}} - E_{\mathrm{free}}(\phi_*,\psi_*)  \big)^{\theta_*} \ \mathrm{d} t \\
    & \quad \le C_{\Omega} \big|  \overline{\phi_0} - \overline{\phi_*}  \big| + \frac{\max \{ 2, m_\phi^*, m_\psi^* \}}{C_{\mathrm{diss}}} \Big(  E_{\mathrm{free}}(\phi_0,\psi_0) - E_{\mathrm{free}}(\phi_*,\psi_*) + (2 C_{\mathrm{source}}' + C_{\mathrm{free},\delta_*}') \big|  \overline{\phi_0} - \overline{\phi_*}  \big|  \Big)^{\theta_*} \\
    & \quad \le L_1 \big(  \|  \phi_0 - \phi_*  \|_V^{\theta_*} + \|  \psi_0 - \psi_*  \|_V^{\theta_*}  \big),
\end{align*}
where $L_1>0$ depends on $\Omega$, $\| \phi_* \|_V$, $\| \psi_* \|_V$, $\delta_*$, $C_{\mathrm{source}}'$, $C_{\mathrm{free},\delta_*}'$ and coefficients of the system.
Consequently, it follows from \eqref{gns_inequality} that
\begin{align*}
    & \| \phi(t) - \phi_* \|_W + \| \psi(t) - \psi_* \|_W \\
    & \quad \le \| \phi_0 - \phi_* \|_W + \| \psi_0 - \psi_* \|_W \\
    & \qquad + C \big(  \| \phi(t) - \phi_0 \|_{H^3(\Omega)}^{\frac{3}{4}} \| \phi(t) - \phi_0 \|_{V^*}^{\frac{1}{4}} + \| \psi(t) - \psi_0 \|_{H^3(\Omega)}^{\frac{3}{4}} \| \psi(t) - \psi_0 \|_{V^*}^{\frac{1}{4}}  \big) \\
    & \quad \le \| \phi_0 - \phi_* \|_W + \| \psi_0 - \psi_* \|_W \\
    & \qquad + C \big(  M_2 + M_{2,\delta_*}'  \big)^{\frac{3}{4}} \big(  \int_0^{\widehat{T}} \| \partial_t \phi(t) \|_{V^*} + \| \partial_t \psi(t) \|_{V^*} \ \mathrm{d} t  \big)^{\frac{1}{4}} \\
    & \quad \le \| \phi_0 - \phi_* \|_W + \| \psi_0 - \psi_* \|_W \\
    & \qquad + L_2 \big(  M_2 + M_{2,\delta_*}'  \big)^{\frac{3}{4}} \big(  \|  \phi_0 - \phi_*  \|_V^{\frac{\theta_*}{4}} + \|  \psi_0 - \psi_*  \|_V^{\frac{\theta_*}{4}}  \big),
\end{align*}
for all $t \in [0,\widehat{T}]$, where $M_{2,\delta_*}'>0$ is the $L^\infty \big(  0,T;H^3(\Omega)  \big)$-bound of $\phi$ and $\psi$ derived from
\[
\widetilde{\mathcal{E}}\big( \bm{u}(t),\phi(t),\psi(t),\mu_\phi(t),\mu_\psi(t) \big) \le 2 E_{2,\frac{9}{10} \delta_*},
\]
given by
\begin{align}
    & \| \phi \|_{H^3(\Omega)} + \| \psi \|_{H^3(\Omega)} \nonumber \\
    & \quad \le C \big(  \| \phi \|_{H^1(\Omega)} + \| \psi \|_{H^1(\Omega)} + \| \nabla \Delta \phi \| + \| \nabla \Delta \psi \|  \big) \nonumber \\
    & \quad \le ( C + P_{\frac{3}{4} \delta_*, 2} ) \big(  \sqrt{E_1} + 1  \big) + C \big( \| \nabla \mu_\phi \| + \| \nabla \mu_\psi \|  \big) \nonumber \\
    & \quad \le ( C + P_{\frac{3}{4} \delta_*, 2} ) \big(  \sqrt{E_1} + 1  \big) + C \sqrt{\frac{2 E_{2,\frac{9}{10} \delta_*}}{L_{\mathrm{coef}}'}} \nonumber \\
    & \quad \stackrel{\mathrm{def}}{=} M_{2,\delta_*}'. \label{H3_estimate_by_higher_order_energy}
\end{align}
Therefore, taking 
\begin{align*}
    \rho_\phi & \stackrel{\mathrm{def}}{=} \min \Big\{  \frac{L_{\mathrm{coef}} \sqrt{| \Omega |} E_3 T_1}{ 12 C_{\mathrm{source}} e^{4 C_{\mathrm{source}}} \big(  E_1 +  | \Omega | \big(  \| G \|_{C(\mathbb{R}^2)} + 2 | \sigma_2 | C_P^2  \big) + 1  \big) \big( 1 + 2 L_{\mathrm{coef}} \big)}, \\
    & \qquad \qquad \frac{\omega}{6}, \Big(  \frac{\omega}{6 L_2 \big(  M_2 + M_{2,\delta_*}'  \big)^{\frac{3}{4}} + 1 }  \Big)^{\frac{4}{\theta_*}}  \Big\}, \\
    \rho_\psi & \stackrel{\mathrm{def}}{=} \min \Big\{  \frac{\omega}{6}, \Big(  \frac{\omega}{6 L_2 \big(  M_2 + M_{2,\delta_*}'  \big)^{\frac{3}{4}} + 1 }  \Big)^{\frac{4}{\theta_*}}  \Big\},
\end{align*}
we obtain
\[
    \| \phi(t) - \phi_* \|_W + \| \psi(t) - \psi_* \|_W < \frac{2}{3} \omega, \quad \forall t \in [0,\widehat{T}],
\]
which contradicts the definition of $\widehat{T}$.
In conclusion, we have $\widehat{T}=T$ if $\rho_\phi$ and $\rho_\psi$ are chosen as above.
In this case, it holds
\begin{align*}
    & \| \phi(t) \|_{C(\overline{\Omega})} \le \| \phi_* \|_{C(\overline{\Omega})} + C_S \|  \phi(t) - \phi_*  \|_W \le 1 - \frac{9}{10} \delta_*, \\
    & \| \psi(t) - \frac{1}{2} \|_{C(\overline{\Omega})} \le \| \psi_* - \frac{1}{2} \|_{C(\overline{\Omega})} + C_S \|  \psi(t) - \psi_*  \|_W \le \frac{1}{2} - \frac{9}{10} \delta_*, \\
    & \|  \phi(t) - \phi_0  \|_W \le \|  \phi(t) - \phi_*  \|_W + \|  \phi_* - \phi_0  \|_W \le \frac{7}{6} \omega \le \frac{L_{\mathrm{coef}} E_3 T_1}{12(C_{\mathrm{free},\frac{3}{4} \delta_*} + 1 )}, \\
    & \|  \psi(t) - \psi_0  \|_W \le \frac{L_{\mathrm{coef}} E_3 T_1}{12(C_{\mathrm{free},\frac{3}{4} \delta_*} + 1 )},
\end{align*}
for any $t\in [0,T]$. 
Moreover, we have
\begin{align*}
    & \int_0^T \widetilde{\mathcal{E}}(t)\, \mathrm{d} t \\
    & \quad \le \int_{Q_T} \frac{3}{4} m_\phi ( \phi ) \big| \nabla \mu_\phi \big|^2 + \frac{1}{2} m_\psi ( \psi ) \big| \nabla \mu_\psi \big|^2 + \frac{1}{2} \nu ( \phi,\psi ) \big| \bm{u} \big|^2 + \frac{r-1}{r} \eta ( \phi,\psi ) \big| \bm{u} \big|^r \ \mathrm{d} (x,t) \nonumber \\
    & \qquad + 4 C_{\mathrm{source}} e^{4 C_{\mathrm{source}}} \big(  E_1 +  | \Omega | \big(  \| G \|_{C(\mathbb{R}^2)} + 2 | \sigma_2 | C_P^2  \big) + 1  \big) \int_0^T  \sigma_1 \big| \overline{\phi}-c \big| \ \mathrm{d} t \\
    & \quad \le \frac{1}{L_{\mathrm{coef}}} \int_{Q_T} \frac{m_{\phi,*}}{2} \big| \nabla \mu_\phi \big|^2 + m_{\psi,*} \big| \nabla \mu_\psi \big|^2 + \nu_* \big| \bm{u} \big|^2 + \eta_* \big| \bm{u} \big|^r \ \mathrm{d} (x,t) \nonumber \\
    & \qquad + 4 C_{\mathrm{source}} e^{4 C_{\mathrm{source}}} \big(  E_1 +  | \Omega | \big(  \| G \|_{C(\mathbb{R}^2)} + 2 | \sigma_2 | C_P^2  \big) + 1  \big) \big| \overline{\phi_0} - \overline{\phi_*} \big| \\
    & \quad \le \frac{1}{L_{\mathrm{coef}}} \big(  E_{\mathrm{free}} (\phi_0,\psi_0) - E_{\mathrm{free}} (\phi(T),\psi(T))  \big) \\
    & \qquad + 2 C_{\mathrm{source}} e^{4 C_{\mathrm{source}}} \big(  E_1 +  | \Omega | \big(  \| G \|_{C(\mathbb{R}^2)} + 2 | \sigma_2 | C_P^2  \big) + 1  \big) \frac{1 + 2 L_{\mathrm{coef}}}{L_{\mathrm{coef}}} \big| \overline{\phi_0} - \overline{\phi_*} \big| \\
    & \quad \le \frac{C_{\mathrm{free},\frac{3}{4} \delta_*}}{L_{\mathrm{coef}}} \big(  \| \phi(T) - \phi_0 \|_W + \| \psi(T) - \psi_0 \|_W  \big) \\
    & \qquad + 2 C_{\mathrm{source}} e^{4 C_{\mathrm{source}}} \big(  E_1 +  | \Omega | \big(  \| G \|_{C(\mathbb{R}^2)} + 2 | \sigma_2 | C_P^2  \big) + 1  \big)  \frac{1 + 2 L_{\mathrm{coef}}}{L_{\mathrm{coef}} \sqrt{| \Omega |}} \rho_\phi \\
    & \quad \le \frac{E_3 T_1}{6} + 2 C_{\mathrm{source}} e^{4 C_{\mathrm{source}}} \big(  E_1 +  | \Omega | \big(  \| G \|_{C(\mathbb{R}^2)} + 2 | \sigma_2 | C_P^2  \big) + 1  \big)  \frac{1 + 2 L_{\mathrm{coef}}}{L_{\mathrm{coef}} \sqrt{| \Omega |}} \rho_\phi \\
    & \quad \le \frac{1}{3} E_3 T_1.
\end{align*}
The proof is complete.
\qed

\medskip

\noindent \textbf{Proof of Theorem \ref{ss_chdf_exist_global_small}.} Without loss of generality, we take $M_2$ such that
\[
M_2 \ge ( C + P_{\frac{3}{4} \delta_*, 2} ) \big(  \sqrt{E_1'} + 1  \big) + C \sqrt{\frac{2 E_3}{L_{\mathrm{coef}}'}} + M,
\]
where the constant $C$ is the same as in the estimate \eqref{H3_estimate_by_higher_order_energy} and $E_1'$ is defined by
\begin{align*}
    & E_1' \stackrel{\mathrm{def}}{=} \frac{1}{2} ( 1 + \| \phi_* \|_V )^2 + \frac{\beta}{2} ( 1 + \| \psi_* \|_V )^2 + | \Omega |(\| F_\phi \|_{C([-1,1])} + \| F_\psi \|_{C([0,1])} + \| G \|_{C(\mathbb{R}^2)} + 2 | \sigma_2 | C_P^2).
\end{align*}
First, choose the radii so that
\[
    r_\phi \le \frac{\omega^2}{1 + (C_{GN}')^2 (M + \| \phi_* \|_{H^3(\Omega)})}, \quad r_\psi \le \frac{\omega^2}{1 + (C_{GN}')^2 (M + \| \psi_* \|_{H^3(\Omega)})}.
\]
Here $C_{GN}'>0$ is the optimal constant, depending on $\Omega$, such that
\[
    \| f \|_{H^2(\Omega)} \le C_{GN}' \| f \|_V^{\frac{1}{2}} \| f \|_{H^3(\Omega)}^{\frac{1}{2}}, \quad \forall f \in H^3(\Omega).
\]
Corollary \ref{ss_local_strong_separate_bounded}, applied with the separation constant $\delta=\tfrac{9}{10} \delta_*$, yields a strong solution $( \bm{u},\phi,\psi,\mu_\phi,\mu_\psi )$ on $[0,T_1]$ satisfying the separation and higher-order energy bounds assumed in Proposition \ref{uniform_control_small_solution}.
Therefore, by further choosing the smaller radii
\[
    r_\phi = \frac{\rho_\phi^2}{1 + (C_{GN}')^2 (M + \| \phi_* \|_{H^3(\Omega)})}, \quad r_\psi = \frac{\rho_\psi^2}{1 + (C_{GN}')^2 (M + \| \psi_* \|_{H^3(\Omega)})},
\]
where $\rho_\phi$ and $\rho_\psi$ are taken as in the proof of Proposition \ref{uniform_control_small_solution},
we can use this proposition to conclude that
\begin{align*}
    & \| \phi(t) - \phi_* \|_W + \| \psi(t) - \psi_* \|_W \le \omega, \quad \forall t \in [0,T_1], \\
    & \int_0^{T_1} \widetilde{\mathcal{E}}(\bm{u},\phi,\psi,\mu_\phi,\mu_\psi) \ \mathrm{d} t \ \le \frac{1}{3} E_3 T_1.
\end{align*}
Therefore, it follows from the Chebyshev inequality that there exists some $T_2 \in [\tfrac{1}{2} T_1, T_1]$ such that
\[
    \widetilde{\mathcal{E}}(T_2) \le E_3,
\]
which implies that
\[
\| \phi(T_2) \|_{H^3(\Omega)} + \| \psi(T_2) \|_{H^3(\Omega)} \le M_2.
\]
Thus, Corollary \ref{ss_local_strong_separate_bounded} ensures that $( \bm{u},\phi,\psi,\mu_\phi,\mu_\psi )$ can be extended to a solution defined in $[0,\tfrac{3}{2} T_1]$, with
\begin{align*}
    & | \phi(t) | \le 1 - \frac{3}{4} \delta_*, \quad \Big| \psi(t) - \frac{1}{2} \Big| \le \frac{1}{2} - \frac{3}{4} \delta_*, \\
    & \widetilde{\mathcal{E}}(\bm{u}(t),\phi(t),\psi(t),\mu_\phi(t),\mu_\psi(t)) \le 2 E_{2,\frac{9}{10} \delta_*}, 
\end{align*}
for all $t \in [0,\frac{3}{2} T_1]$. 
Then, applying Proposition \ref{uniform_control_small_solution} again, we have
\begin{align*}
    & \| \phi(t) - \phi_* \|_W + \| \psi(t) - \psi_* \|_W \le \omega, \quad \forall t \in [0,\frac{3}{2} T_1], \\
    & \int_0^{\frac{3}{2}T_1} \widetilde{\mathcal{E}}(\bm{u},\phi,\psi,\mu_\phi,\mu_\psi) \ \mathrm{d} t \ \le \frac{1}{3} E_3 T_1.
\end{align*}
Consequently, iterating the above argument, we can extend $( \bm{u},\phi,\psi,\mu_\phi,\mu_\psi )$ to $[0,\infty)$ with
\begin{align*}
    & \| \phi(t) - \phi_* \|_W + \| \psi(t) - \psi_* \|_W \le \omega, \ \forall t \in [0,\infty), \\
    & \widetilde{\mathcal{E}}(\bm{u}(t),\phi(t),\psi(t),\mu_\phi(t),\mu_\psi(t)) \le 2 E_{2,\frac{9}{10} \delta_*}, \ \forall t \in [0,\infty), \\
    & \int_0^\infty \widetilde{\mathcal{E}} (\bm{u},\phi,\psi,\mu_\phi,\mu_\psi) \ \mathrm{d} t \ \le \frac{1}{3} E_3 T_1.
\end{align*}
The proof of Theorem \ref{ss_chdf_exist_global_small} is complete.
\qed

\section{Proof of Theorem \ref{ss_chdf_exist_2d}}
\label{ss_2d}
\setcounter{equation}{0}

In this section, we prove global existence in two dimensions for arbitrarily large initial data, and uniqueness under the additional assumptions stated in Theorem \ref{ss_chdf_exist_2d}.
 \smallskip 
 
\noindent \textbf{Proof of Theorem \ref{ss_chdf_exist_2d}.} We still obtain the energy estimate \eqref{basic_energy_estimate} and a bound for the initial value as in Section \ref{ss_3d_proof}, both of which are uniform with respect to $m$, $n$ and $k$.
Then, we refine estimate \eqref{higher_order_estimate_uniform_3d} as follows:
\begin{align}
    & \frac{\mathrm{d}}{\mathrm{d}t} \mathcal{E} ( \bm{u}_{m,n},\phi_{m,n},\psi_{m,n},\mu_{\phi,m,n},\mu_{\psi,m,n} ) \nonumber \\
    & \qquad + \big\| \nabla \partial_t \phi_{m,n} \big\|^2 + \beta \big\| \nabla \partial_t \psi_{m,n} \big\|^2 \nonumber \\
    & \quad \le C \big(  \big\| \partial_t \phi_{m,n} \big\| \| \nabla \mu_{\phi,m,n} \|^2_{L^4(\Omega)} + \big\| \partial_t \psi_{m,n} \big\| \big\| \nabla \mu_{\psi,m,n} \big\|^2_{L^4(\Omega)}  \big) \nonumber \\
    & \qquad + C \big(  \big\| \partial_t \phi_{m,n} \big\| + \big\| \partial_t \psi_{m,n} \big\|  \big) \big(  \| \bm{u}_{m,n} \|^2_{L^4(\Omega)} + \| \bm{u}_{m,n} \|^r_{L^{2r}(\Omega)}  \big) \nonumber \\
    & \qquad + C \| \bm{u}_{m,n} \|_{L^4(\Omega)} \big(  \big\|  \partial_t \phi_{m,n} - \overline{\partial_t \phi_{m,n}}  \big\|_{L^4(\Omega)} \| \nabla \mu_{\phi,m,n} \| + \big\| \partial_t \psi_{m,n} \big\|_{L^4(\Omega)} \big\| \nabla \mu_{\psi,m,n} \big\|  \big) \nonumber \\
    & \qquad + C(1 + \| \mu_{\phi,m,n} \|) + C \big(  \big| \overline{\partial_t \phi_{m,n}} \big|^2 + \big\| \partial_t \phi_{m,n} - \overline{\partial_t \phi_{m,n}} \big\|^2 + \big\| \partial_t \psi_{m,n} \big\|^2  \big). \label{higher_order_estimate_uniform_2d}
\end{align}
In dimension two, estimate \eqref{H1_sp_3d} can also be refined. Indeed, we have
\begin{align}
    & \nu_* \| \nabla \bm{u}_{m,n} \|^2 + \eta_* \big\|  \nabla | \bm{u}_{m,n} |^\frac{r}{2}  \big\|^2 \nonumber \\
    & \quad \le C \big(  \| \nabla \phi_{m,n} \|_{L^4(\Omega)} + \| \nabla \psi_{m,n} \|_{L^4(\Omega)}  \big) \big(  \| \bm{u}_{m,n} \|_{L^4(\Omega)} \| \nabla \bm{u}_{m,n} \| + \big\|  | \bm{u}_{m,n} |^\frac{r}{2}  \big\|_{L^4(\Omega)} \big\|  \nabla | \bm{u}_{m,n} |^\frac{r}{2}  \big\|  \big) \nonumber \\
    & \qquad + C \| \nabla \bm{u}_{m,n} \| \big(  \| \nabla \phi_{m,n} \|_{L^\infty(\Omega)} \| \nabla \mu_{\phi,m,n} \| + \| \nabla \psi_{m,n} \|_{L^\infty(\Omega)} \| \nabla \mu_{\psi,m,n} \|  \big) \nonumber \\
    & \qquad + C \big(  \| \bm{u}_{m,n} \|^2 + \| \bm{u}_{m,n} \|^r_{L^r(\Omega)}  \big) \nonumber \\
    & \quad \le C \big(  1 + \| \nabla \mu_{\phi,m,n} \|^\frac{1}{4} + \| \nabla \mu_{\psi,m,n} \|^\frac{1}{4}  \big) \big(  \| \bm{u}_{m,n} \|^\frac{1}{2} \|  \bm{u}_{m,n} \|_V^\frac{3}{2} + \big\|  | \bm{u}_{m,n} |^\frac{r}{2}  \big\|^\frac{1}{2} \big\| | \bm{u}_{m,n} |^\frac{r}{2} \big\|_V^\frac{3}{2}  \big) \nonumber \\
    & \qquad + C \| \nabla \bm{u}_{m,n} \| \big( 1 + \| \nabla \mu_{\phi,m,n} \| \big)^\frac{3}{2} \ln^\frac{1}{2} \big(  e + \| \nabla \mu_{\phi,m,n} \|  \big) \nonumber \\
    & \qquad + C \| \nabla \bm{u}_{m,n} \| \big( 1 + \| \nabla \mu_{\psi,m,n} \| \big)^\frac{3}{2} \ln^\frac{1}{2} \big(  e + \| \nabla \mu_{\psi,m,n} \|  \big) \nonumber \\
    & \quad \le \frac{1}{2} \big(  \nu_* \| \nabla \bm{u}_{m,n} \|^2 + \eta_* \big\|  \nabla |  \bm{u}_{m,n} |^\frac{r}{2}  \big\|^2  \big) + C \big(  1 + \| \nabla \mu_{\phi,m,n} \| + \| \nabla \mu_{\psi,m,n} \|  \big) \big(  \| \bm{u}_{m,n} \|^2 + \big\|  | \bm{u}_{m,n} |^\frac{r}{2}  \big\|^2  \big) \nonumber \\
    & \qquad + C \big( 1 + \| \nabla \mu_{\phi,m,n} \| \big)^3 \ln \big(  e + \| \nabla \mu_{\phi,m,n} \|  \big) + C\big( 1 + \| \nabla \mu_{\psi,m,n} \| \big)^3 \ln \big(  e + \| \nabla \mu_{\psi,m,n} \|  \big),\label{H1_sp_2d}
\end{align}
where we used \eqref{gns_inequality} and \eqref{BGW}.
Inequality \eqref{H1_sp_2d}
implies
\begin{align*}
    & \| \bm{u}_{m,n} \|_V^2 + \big\|  | \bm{u}_{m,n} |^\frac{r}{2}  \big\|_V^2 \\
    & \quad \le C \big( 1 + C_1 + \mathcal{E} ( \bm{u}_{m,n},\phi_{m,n},\psi_{m,n},\mu_{\phi,m,n},\mu_{\psi,m,n} ) \big)^\frac{3}{2} \ln \big(  e + C_1 + \mathcal{E} ( \bm{u}_{m,n},\phi_{m,n},\psi_{m,n},\mu_{\phi,m,n},\mu_{\psi,m,n} )  \big).
\end{align*}
Here $C_1$ is defined in \eqref{C_1_define}.
Next, arguing as in \cite{G2020} and using \eqref{BGW} once more, we can improve the estimate of $\partial_t \phi_{m,n}$ and $\partial_t \psi_{m,n}$ as follows:
\begin{align*}
    \big\| \partial_t \phi_{m,n}- \overline{\partial_t \phi_{m,n}} \big\|_* & \le C \big(  \| \nabla \mu_{\phi,m,n} \| + \| \bm{u}_{m,n} \| \| \phi_{m,n} - \overline{\phi_{m,n}} \|_{L^\infty(\Omega)}  \big) \\
    & \le C \big(  \| \nabla \mu_{\phi,m,n} \| + \| \bm{u}_{m,n} \| (1 + \| \nabla \phi_{m,n} \|) \ln^\frac{1}{2} \big( e + \| \nabla \mu_{\phi,m,n} \| \big)  \big) \\
    & \le C \big( 1 + C_1 + \mathcal{E} ( \bm{u}_{m,n},\phi_{m,n},\psi_{m,n},\mu_{\phi,m,n},\mu_{\psi,m,n} ) \big)^\frac{1}{2} \\
    & \quad \ln^\frac{1}{2} \big(  e + C_1 + \mathcal{E} ( \bm{u}_{m,n},\phi_{m,n},\psi_{m,n},\mu_{\phi,m,n},\mu_{\psi,m,n} )  \big),
\end{align*}
\begin{align*}
    \| \partial_t \psi_{m,n} \|_* & \le C \big( 1 + C_1 + \mathcal{E} ( \bm{u}_{m,n},\phi_{m,n},\psi_{m,n},\mu_{\phi,m,n},\mu_{\psi,m,n} ) \big)^\frac{1}{2} \\
    & \quad \ln^\frac{1}{2} \big(  e + C_1 + \mathcal{E} ( \bm{u}_{m,n},\phi_{m,n},\psi_{m,n},\mu_{\phi,m,n},\mu_{\psi,m,n} )  \big).
\end{align*}
Thanks to the above estimates, arguing as in \cite{GOW2026,CGGG2025}, we can obtain the $H^2$-estimate of $\mu_{\phi,m,n}$ and $\mu_{\psi,m,n}$, namely,
\begin{align*}
    \| \mu_{\phi,m,n}-\overline{\mu_{\phi,m,n}} \|_W & = \big\|  \mathcal{N}_{\phi, \phi_{m,n}} \big(  \partial_t \phi_{m,n}- \overline{\partial_t \phi_{m,n}} + \bm{u}_{m,n} \cdot \nabla \phi_{m,n}  \big)  \big\|_W \\
    & \le C \| \nabla \phi_{m,n} \| \| \phi_{m,n} \|_W \big\|  \nabla \mathcal{N}_{\phi, \phi_{m,n}} \big(  \partial_t \phi_{m,n}- \overline{\partial_t \phi_{m,n}} + \bm{u}_{m,n} \cdot \nabla \phi_{m,n}  \big)  \big\| \\
    & \quad + C \|  \partial_t \phi_{m,n}- \overline{\partial_t \phi_{m,n}} + \bm{u}_{m,n} \cdot \nabla \phi_{m,n}  \| \\
    & \le C \| \phi_{m,n} \|_W \| \nabla \mu_{\phi,m,n} \| + C \| \bm{u}_{m,n} \|_{L^4(\Omega)}  \| \nabla \phi_{m,n} \|_{L^4(\Omega)}  \\
    & \quad + C  \big\| \nabla \partial_t \phi_{m,n} \big\|^\frac{1}{2} \| \partial_t \phi_{m,n}- \overline{\partial_t \phi_{m,n}} \|_*^\frac{1}{2} \\
    & \le C \| \nabla \mu_{\phi,m,n} \| \big( 1 + \| \nabla \mu_{\phi,m,n} \| \big)^\frac{1}{2} + C \| \bm{u}_{m,n} \|^\frac{1}{2} \| \bm{u}_{m,n} \|_V^\frac{1}{2} (1 + \| \nabla \mu_{\phi,m,n} \|)^\frac{1}{4} \\ 
    & \quad + C\| \nabla \partial_t \phi_{m,n} \|^\frac{1}{2} \| \bm{u}_{m,n} \|^\frac{1}{2} \ln^\frac{1}{4} (e + \| \nabla \mu_{\phi,m,n} \|) + C \| \nabla \mu_{\phi,m,n} \|^\frac{1}{2} \| \nabla \partial_t \phi_{m,n} \|^\frac{1}{2} \\
    & \le C \| \nabla \partial_t \phi_{m,n} \|^\frac{1}{2} \big( 1 + C_1 + \mathcal{E} ( \bm{u}_{m,n},\phi_{m,n},\psi_{m,n},\mu_{\phi,m,n},\mu_{\psi,m,n} ) \big)^\frac{1}{4} \\
    & \quad \ln^\frac{1}{4} \big(  e + C_1 + \mathcal{E} ( \bm{u}_{m,n},\phi_{m,n},\psi_{m,n},\mu_{\phi,m,n},\mu_{\psi,m,n} )  \big) \\
    & \quad + C \big( 1 + C_1 + \mathcal{E} ( \bm{u}_{m,n},\phi_{m,n},\psi_{m,n},\mu_{\phi,m,n},\mu_{\psi,m,n} ) \big)^\frac{3}{4} \\
    & \quad \ln^\frac{1}{4} \big(  e + C_1 + \mathcal{E} ( \bm{u}_{m,n},\phi_{m,n},\psi_{m,n},\mu_{\phi,m,n},\mu_{\psi,m,n} )  \big).
\end{align*}
Thus, we have
\begin{align*}
    & \| \partial_t \phi_{m,n} \| \| \nabla \mu_{\phi,m,n} \|_{L^4(\Omega)}^2 \\
    & \le C \big(  1 + \| \partial_t \phi_{m,n} - \overline{\partial_t \phi_{m,n}} \|_*^\frac{1}{2} \| \nabla \partial_t \phi_{m,n} \|^\frac{1}{2}  \big) \| \nabla \mu_{\phi,m,n} \| \| \mu_{\phi,m,n} - \overline{\mu_{\phi,m,n}} \|_W \\
    & \le C \| \nabla \partial_t \phi_{m,n} \|^\frac{1}{2} \big( 1 + C_1 + \mathcal{E} ( \bm{u}_{m,n},\phi_{m,n},\psi_{m,n},\mu_{\phi,m,n},\mu_{\psi,m,n} ) \big)^\frac{3}{4} \\
    & \quad \ln^\frac{1}{4} \big(  e + C_1 + \mathcal{E} ( \bm{u}_{m,n},\phi_{m,n},\psi_{m,n},\mu_{\phi,m,n},\mu_{\psi,m,n} )  \big) \\
    & \quad + C \big( 1 + C_1 + \mathcal{E} ( \bm{u}_{m,n},\phi_{m,n},\psi_{m,n},\mu_{\phi,m,n},\mu_{\psi,m,n} ) \big)^\frac{5}{4} \\
    & \quad \ln^\frac{1}{4} \big(  e + C_1 + \mathcal{E} ( \bm{u}_{m,n},\phi_{m,n},\psi_{m,n},\mu_{\phi,m,n},\mu_{\psi,m,n} )  \big) \\
    & \quad + C \| \nabla \partial_t \phi_{m,n} \| \big( 1 + C_1 + \mathcal{E} ( \bm{u}_{m,n},\phi_{m,n},\psi_{m,n},\mu_{\phi,m,n},\mu_{\psi,m,n} ) \big) \\
    & \quad \ln^\frac{1}{2} \big(  e + C_1 + \mathcal{E} ( \bm{u}_{m,n},\phi_{m,n},\psi_{m,n},\mu_{\phi,m,n},\mu_{\psi,m,n} )  \big) \\
    & \quad + C \| \nabla \partial_t \phi_{m,n} \|^\frac{1}{2} \big( 1 + C_1 + \mathcal{E} ( \bm{u}_{m,n},\phi_{m,n},\psi_{m,n},\mu_{\phi,m,n},\mu_{\psi,m,n} ) \big)^\frac{3}{2} \\
    & \quad \ln^\frac{1}{2} \big(  e + C_1 + \mathcal{E} ( \bm{u}_{m,n},\phi_{m,n},\psi_{m,n},\mu_{\phi,m,n},\mu_{\psi,m,n} )  \big) \\
    & \le \epsilon \| \nabla \partial_t \phi_{m,n} \|^2 \\
    & \quad + C_\epsilon \big( 1 + C_1 + \mathcal{E} ( \bm{u}_{m,n},\phi_{m,n},\psi_{m,n},\mu_{\phi,m,n},\mu_{\psi,m,n} ) \big)^2 \\
    & \quad \ln \big(  e + C_1 + \mathcal{E} ( \bm{u}_{m,n},\phi_{m,n},\psi_{m,n},\mu_{\phi,m,n},\mu_{\psi,m,n} )  \big),
\end{align*}
and a similar estimate holds for the term $\| \partial_t \psi_{m,n} \| \| \nabla \mu_{\psi,m,n} \|_{L^4(\Omega)}^2$.
In addition, we have
\begin{align*}
    & \big\| \partial_t \phi_{m,n} \big\| \big(  \| \bm{u}_{m,n} \|^2_{L^4(\Omega)} + \| \bm{u}_{m,n} \|^r_{L^{2r}(\Omega)}  \big) \\
    &\quad \le C \big(  \| \bm{u}_{m,n} \| + \| \bm{u}_{m,n} \|^\frac{r}{2}_{L^{r}(\Omega)}  \big) \big(  \| \bm{u}_{m,n} \|_V + \big\|  | \bm{u}_{m,n} |^\frac{r}{2}  \big\|_V  \big) \big(  1 + \big\| \partial_t \phi_{m,n} - \overline{\partial_t \phi_{m,n}} \big\|_*^\frac{1}{2} \| \nabla \partial_t \phi_{m,n} \|^\frac{1}{2}  \big) \\
    & \quad \le C \big( 1 + C_1 + \mathcal{E} ( \bm{u}_{m,n},\phi_{m,n},\psi_{m,n},\mu_{\phi,m,n},\mu_{\psi,m,n} ) \big)^\frac{5}{4} \\
    & \qquad \ln^\frac{1}{2} \big(  e + C_1 + \mathcal{E} ( \bm{u}_{m,n},\phi_{m,n},\psi_{m,n},\mu_{\phi,m,n},\mu_{\psi,m,n} )  \big) \\
    & \qquad + C \| \nabla \partial_t \phi_{m,n} \|^\frac{1}{2} \big( 1 + C_1 + \mathcal{E} ( \bm{u}_{m,n},\phi_{m,n},\psi_{m,n},\mu_{\phi,m,n},\mu_{\psi,m,n} ) \big)^\frac{3}{2} \\
    & \qquad \ln^\frac{3}{4} \big(  e + C_1 + \mathcal{E} ( \bm{u}_{m,n},\phi_{m,n},\psi_{m,n},\mu_{\phi,m,n},\mu_{\psi,m,n} )  \big) \\
    & \quad \le \epsilon \| \nabla \partial_t \phi_{m,n} \|^2 \\
    & \qquad + C_\epsilon \big( 1 + C_1 + \mathcal{E} ( \bm{u}_{m,n},\phi_{m,n},\psi_{m,n},\mu_{\phi,m,n},\mu_{\psi,m,n} ) \big)^2 \\
    & \qquad \ln \big(  e + C_1 + \mathcal{E} ( \bm{u}_{m,n},\phi_{m,n},\psi_{m,n},\mu_{\phi,m,n},\mu_{\psi,m,n} )  \big),
\end{align*}
and an analogous estimate holds for $\big\| \partial_t \psi_{m,n} \big\| \big(  \| \bm{u}_{m,n} \|^2_{L^4(\Omega)} + \| \bm{u}_{m,n} \|^r_{L^{2r}(\Omega)}  \big)$.
Consider now the third term on the right-hand side of \eqref{higher_order_estimate_uniform_2d} and observe that
\begin{align*}
    & \| \bm{u}_{m,n} \|_{L^4(\Omega)} \big\|  \partial_t \phi_{m,n} - \overline{\partial_t \phi_{m,n}}  \big\|_{L^4(\Omega)} \| \nabla \mu_{\phi,m,n} \| \\
    & \quad \le C \| \bm{u}_{m,n} \|^\frac{1}{2} \| \bm{u}_{m,n} \|_V^\frac{1}{2} \big\|  \partial_t \phi_{m,n} - \overline{\partial_t \phi_{m,n}}  \big\|_*^\frac{1}{4} \big\|  \nabla \partial_t \phi_{m,n}  \big\|^\frac{3}{4} \| \nabla \mu_{\phi,m,n} \| \\
    & \quad \le C \big\|  \nabla \partial_t \phi_{m,n}  \big\|^\frac{3}{4} \big( 1 + C_1 + \mathcal{E} ( \bm{u}_{m,n},\phi_{m,n},\psi_{m,n},\mu_{\phi,m,n},\mu_{\psi,m,n} ) \big)^\frac{5}{4} \\
    & \qquad \ln^\frac{3}{8} \big(  e + C_1 + \mathcal{E} ( \bm{u}_{m,n},\phi_{m,n},\psi_{m,n},\mu_{\phi,m,n},\mu_{\psi,m,n} )  \big) \\
    & \quad \le \epsilon \| \nabla \partial_t \phi_{m,n} \|^2 \\
    & \qquad + C_\epsilon \big( 1 + C_1 + \mathcal{E} ( \bm{u}_{m,n},\phi_{m,n},\psi_{m,n},\mu_{\phi,m,n},\mu_{\psi,m,n} ) \big)^2 \\
    & \qquad \ln^\frac{3}{5} \big(  e + C_1 + \mathcal{E} ( \bm{u}_{m,n},\phi_{m,n},\psi_{m,n},\mu_{\phi,m,n},\mu_{\psi,m,n} )  \big).
\end{align*}
The term $\| \bm{u}_{m,n} \|_{L^4(\Omega)} \big\| \partial_t \psi_{m,n} \big\|_{L^4(\Omega)} \| \nabla \mu_{\psi,m,n} \|$
can be bounded in a similar way.
In conclusion, we obtain
\begin{align*}
    & \frac{\mathrm{d}}{\mathrm{d}t} \mathcal{E} ( \bm{u}_{m,n},\phi_{m,n},\psi_{m,n},\mu_{\phi,m,n},\mu_{\psi,m,n} ) + \frac{1}{2} \| \nabla \partial_t \phi_{m,n} \|^2 + \frac{\beta}{2} \| \nabla \partial_t \psi_{m,n} \|^2 \\
    & \quad \le C \big( 1 + C_1 + \mathcal{E} ( \bm{u}_{m,n},\phi_{m,n},\psi_{m,n},\mu_{\phi,m,n},\mu_{\psi,m,n} ) \big)^2 \ln \big(  e + C_1 + \mathcal{E} ( \bm{u}_{m,n},\phi_{m,n},\psi_{m,n},\mu_{\phi,m,n},\mu_{\psi,m,n} )  \big).
\end{align*}
It then follows from \cite[Lemma A.1]{G2020} that we can extend the solution $(\bm{u}_{m,n},\phi_{m,n},\psi_{m,n},\mu_{\phi,m,n},\mu_{\psi,m,n})$ of problem \eqref{eq:semigalerkin_chdf1}--\eqref{eq:semigalerkin_chdfi} to the whole interval $[0,\infty)$.
The basic energy estimate \eqref{basic_energy_estimate} gives the uniform integrability of $\mathcal{E} ( \bm{u}_{m,n},\phi_{m,n},\psi_{m,n},\mu_{\phi,m,n},\mu_{\psi,m,n} )$ in $[0,\infty)$.
Then, arguing as in Section \ref{ss_3d_proof}, we also deduce the uniform $L^2_{\mathrm{uloc}}(0,\infty;H^3(\Omega))$-bound for $\mu_{\phi,m,n}$ and $\mu_{\psi,m,n}$.
Hence, the combination of \cite[Lemma A.1, Lemma A.2]{G2020} and the estimates obtained above ensures the uniform boundedness in the following spaces with respect to $m$, $n$, $k$:
\begin{align*}
    & \bm{u}_{m,n} \in L^\infty \big( 0,\infty;\bm{L}^{3r}_\sigma(\Omega) \cap \bm{H}^1(\Omega) \big), \quad \text{with} \ | \bm{u}_{m,n} |^\frac{r}{2} \in L^\infty \big( 0,\infty;H^1(\Omega) \big), \\
    & \phi_{m,n}, \psi_{m,n} \in L^\infty \big( 0,\infty;W^{2,p}(\Omega) \big) \cap H^1_{\mathrm{uloc}} \big( 0,\infty;V \big) \cap W^{1,\infty} \big( 0,\infty;V^* \big), \  \forall p \in [2,\infty), \\
    & \mu_{\phi,m,n},\mu_{\psi,m,n} \in L^\infty \big( 0,\infty;V \big) \cap L^2_{\mathrm{uloc}} \big( 0,\infty;H^3(\Omega) \big).
\end{align*}
Therefore, standard compactness results,  combined with a diagonal argument, yield a global-in-time strong solution to the problem \eqref{eq:chdf1}--\eqref{eq:chdfi}.

Under the additional assumptions stated in Theorem \ref{ss_chdf_exist_2d}, the additional regularity can be deduced by arguing as in \cite{GalPoi2025}.
The resulting strict separation property allows us to prove uniqueness (cf. the proof of Proposition \ref{approximate_solution_exists}).
The proof is complete. \qed

\medskip
\appendix

\section{An Auxiliary Lemma}\label{appendix}
\setcounter{equation}{0}

The following lemma helps us to justify some formal differentiation used in previous sections.

\begin{lemma} \label{chain_rule}
Suppose that $\Omega$ is a bounded smooth domain in $\mathbb{R}^n$, $n \in \mathbb{N}_+$,
and set $\partial_0 \stackrel{\rm{def}}{=} \partial_t$ and $\partial_k \stackrel{\rm{def}}{=} \partial_{x_k}$, $k \in [1,n] \cap \mathbb{N}$.
Then for any $T>0$, $p \in (n,\infty)$, $d \in \mathbb{N}_+$ and $f \in C^1 (\mathbb{R})$ with $f(0)=0$,
we have 
\[f(| \bm{u} |) \in H^1 \big(  0,T;L^q(\Omega)  \big) \cap C_w \big(  [0,T];W^{1,p}(\Omega)  \big) \cap C \big( \overline{\Omega} \times [0,T] \big),\]
\[f(|\bm{u}|) \bm{e}_{\mathbb{R}^d}(\bm{u}) \in H^1 \big( 0,T;\bm{L}^q(\Omega) \big) \cap C_w \big( [0,T];\bm{W}^{1,p}(\Omega) \big) \cap C \big( \overline{\Omega} \times [0,T] \big)^d,\]
for all $\bm{u} \in H^1 \big( 0,T;C \big( \overline{\Omega} \big)^d \big) \cap C_w \big( [0,T];\bm{W}^{1,p}(\Omega) \big)$ and $q \in [1,\infty)$, with
\[ \partial_i \big(  f(|\bm{u}|) \bm{e}_{\mathbb{R}^d}(\bm{u})  \big) = \Big(  f'(| \bm{u} |) - \frac{f(| \bm{u} |)}{| \bm{u} |}  \Big) \big(  \bm{e}_{\mathbb{R}^d}(\bm{u}) \cdot \partial_i \bm{u}  \big) \bm{e}_{\mathbb{R}^d}(\bm{u}) + \frac{f(| \bm{u} |)}{| \bm{u} |} \partial_i \bm{u}, \]
\[ \partial_i \big(  f(|\bm{u}|)  \big) = f'(| \bm{u} |) \partial_i | \bm{u} | = f'(| \bm{u} |) \bm{e}_{\mathbb{R}^d}(\bm{u}) \cdot \partial_i \bm{u}, \]
for all $i \in [0,n] \cap \mathbb{N}$,
where we define
\[
    \left.\frac{f(s)}{s}\right|_{s=0} \stackrel{\mathrm{def}}{=} f'(0).
\]
In particular, it holds
\[| \bm{u} | \in H^1 \big(  0,T;L^q(\Omega)  \big) \cap C_w \big( [0,T];W^{1,p}(\Omega) \big) \cap C \big( \overline{\Omega} \times [0,T] \big)\]
for all $q \in [1,\infty)$, with
\[ \partial_i | \bm{u} | = \bm{e}_{\mathbb{R}^d}(\bm{u}) \cdot \partial_i \bm{u}, \]
for all $i \in [0,n] \cap \mathbb{N}$.
\end{lemma}

\noindent \textbf{Proof.} First, for any given $a >0$, we know that the function $h_a:\mathbb{R}^d \to \mathbb{R}$, $h_a(\bm{r})=\sqrt{ a^2 + |\bm{r}|^2 }$ for all $\bm{r} \in \mathbb{R}^d$, is of class $C^\infty(\mathbb{R}^d;\mathbb{R})$,
and $g_a:\mathbb{R}^d \to \mathbb{R}^d$, $g_a(\bm{r})=\tfrac{f \circ h_a (\bm{r})}{h_a(\bm{r})} \bm{r}$ for all $\bm{r} \in \mathbb{R}^d$, is of class $C^1(\mathbb{R}^d;\mathbb{R}^d)$. 
Since $C^\infty(\overline{\Omega})$ is dense in $W^{1,p}(\Omega)$, we have the density of $C^\infty \big(  \overline{\Omega} \times [0,T]  \big)^d$ in $H^1 \big( 0,T;C \big( \overline{\Omega} \big)^d \big) \cap L^\infty \big( [0,T];\bm{W}^{1,p}(\Omega) \big)$ in the sense that for any $\bm{v} \in H^1 \big( 0,T;C \big( \overline{\Omega} \big)^d \big) \cap L^\infty \big( [0,T];\bm{W}^{1,p}(\Omega) \big)$, 
there exists a sequence $\{ \bm{v}^l \}_{l=1}^\infty \subset C^\infty \big(  \overline{\Omega} \times [0,T]  \big)^d$ such that
\begin{align*}
    & \bm{v}^l \stackrel{*}{\rightharpoonup} \bm{v} \ \text{in $L^\infty \big(  0,T;\bm{W}^{1,p}(\Omega)  \big)$},\\
    & \bm{v}^l \to \bm{v} \ \text{in $\bm{W}^{1,p}(\Omega)$ almost everywhere in $(0,T)$}, \\
    & \bm{v}^l \to \bm{v} \ \text{in $H^1 \big(  0,T;C\big( \overline{\Omega} \big)^d  \big)$},
\end{align*}
as $l\to\infty$.
Then it follows from a density argument that $h_a(\bm{u}), f \circ h_a (\bm{u}) \in H^1 \big(  0,T;C \big( \overline{\Omega} \big)  \big) \cap L^\infty \big(  [0,T];W^{1,p}(\Omega)  \big)$ and $ g_a(\bm{u}) \in H^1 \big(  0,T;C \big( \overline{\Omega} \big)^d  \big) \cap L^\infty \big(  [0,T];\bm{W}^{1,p}(\Omega)  \big)$, with
\begin{align*}
    & \partial_i \big(  h_a (\bm{u})  \big) = \frac{\bm{u}}{h_a(\bm{u})} \cdot \partial_i \bm{u}, \\
    & \partial_i \big(  f \circ h_a (\bm{u})  \big) = f' \circ h_a (\bm{u}) \frac{\bm{u}}{h_a(\bm{u})} \cdot \partial_i \bm{u}, \\
    & \partial_i \big(  g_a (\bm{u})  \big) = \Big(  f'\circ h_a(\bm{u}) - \frac{f \circ h_a(\bm{u})}{h_a(\bm{u})}  \Big) \Big(  \frac{\bm{u}}{h_a(\bm{u})} \cdot \partial_i \bm{u}  \Big) \frac{\bm{u}}{h_a(\bm{u})} + \frac{f \circ h_a(\bm{u})}{h_a(\bm{u})} \partial_i \bm{u},
\end{align*}
for all $i \in [0,n] \cap \mathbb{N}$.
So $\partial_t \big(  h_a (\bm{u})  \big) \in L^2 \big( 0,T;C\big( \overline{\Omega} \big) \big)$, $\partial_t \big(  g_a (\bm{u})  \big) \in L^2 \big( 0,T;C\big( \overline{\Omega} \big)^d \big)$, $\partial_i \big(  h_a (\bm{u})  \big) \in L^\infty \big( 0,T;L^p(\Omega) \big)$ and $\partial_i \big(  g_a (\bm{u})  \big) \in L^\infty \big( 0,T;\bm{L}^p(\Omega) \big)$, $1 \le i \le n$, are uniformly bounded with respect to $a \in (0,1]$ in the corresponding spaces. 
Since we have $\bm{u} \in H^1\big( 0,T;C\big( \overline{\Omega} \big)^d \big) \subset C\big( [0,T];C\big( \overline{\Omega} \big)^d \big) = C \big(  \overline{\Omega} \times [0,T]  \big)^d$,
it follows from $h_a(\bm{r}) \to | \bm{r} |$ in $C \big( \mathbb{R}^d;\mathbb{R} \big)$ uniformly with respect to $\bm{r} \in \mathbb{R}^d$ as $a \to 0^+$ that
\[ h_a (\bm{u}) \to | \bm{u} | \ \text{in } C \big( [0,T];C\big( \overline{\Omega} \big) \big), \]
which implies
\[ f \circ h_a (\bm{u}) \to f(| \bm{u} |) \ \text{in } C \big( [0,T];C\big( \overline{\Omega} \big) \big). \]
In addition, for any given $M>0$, we have $g_a(\bm{r}) \to \tfrac{f(| \bm{r} |)}{| \bm{r} |} \bm{r}$ in $\mathbb{R}^d$ uniformly with respect to $\bm{r} \in B_{\mathbb{R}^d}(\bm{0};M)$. This entails
\[ g_a(\bm{u}) \to \frac{f(| \bm{u} |)}{| \bm{u} |} \bm{u} \ \text{in } C \big( [0,T];C\big(  \overline{\Omega} \big)^d  \big), \]
so that
\[ h_a (\bm{u}) \stackrel{*}{\rightharpoonup} | \bm{u} |, \ f \circ h_a (\bm{u}) \stackrel{*}{\rightharpoonup} f (| \bm{u} |) \ \text{in } L^\infty \big( 0,T;W^{1,p}(\Omega) \big), \]
and
\[ g_a (\bm{u}) \stackrel{*}{\rightharpoonup} \frac{f(| \bm{u} |)}{| \bm{u} |} \bm{u} \ \text{in } L^\infty \big( 0,T;\bm{W}^{1,p}(\Omega) \big). \]
In addition, it follows from Lebesgue's dominated convergence theorem that
\begin{align*}
    & \partial_i \big(  h_a (\bm{u})  \big) \stackrel{*}{\rightharpoonup} \bm{e}_{\mathbb{R}^d} (\bm{u}) \cdot \partial_i \bm{u}, \ \partial_i \big(  f \circ h_a (\bm{u})  \big) \stackrel{*}{\rightharpoonup} f'(| \bm{u} |) \bm{e}_{\mathbb{R}^d} (\bm{u}) \cdot \partial_i \bm{u} \ \text{in $L^\infty \big( 0,T;L^p(\Omega) \big)$},\\
    & \partial_i \big( g_a (\bm{u}) \big) \stackrel{*}{\rightharpoonup} \Big(  f'(| \bm{u} |) - \frac{f(| \bm{u} |)}{| \bm{u} |}  \Big) \big(  \bm{e}_{\mathbb{R}^d}(\bm{u}) \cdot \partial_i \bm{u}  \big) \bm{e}_{\mathbb{R}^d}(\bm{u}) + \frac{f(| \bm{u} |)}{| \bm{u} |} \partial_i \bm{u} \ \text{in $L^\infty \big( 0,T;\bm{L}^p(\Omega) \big)$},
\end{align*}
for $1\le i \le n$. This implies
\begin{align*}
    & \partial_i | \bm{u} | = \bm{e}_{\mathbb{R}^d} (\bm{u}) \cdot \partial_i \bm{u}, \quad \partial_i f (| \bm{u} |) = f' (| \bm{u} |) \bm{e}_{\mathbb{R}^d} (\bm{u}) \cdot \partial_i \bm{u} = f' (| \bm{u} |) \partial_i | \bm{u} |,\\
    & \partial_i \Big(  \frac{f(| \bm{u} |)}{| \bm{u} |} \bm{u}  \Big) = \Big(  f'(| \bm{u} |) - \frac{f(| \bm{u} |)}{| \bm{u} |}  \Big) \big(  \bm{e}_{\mathbb{R}^d}(\bm{u}) \cdot \partial_i \bm{u}  \big) \bm{e}_{\mathbb{R}^d}(\bm{u}) + \frac{f(| \bm{u} |)}{| \bm{u} |} \partial_i \bm{u},
\end{align*}
for all $i \in [1,n] \cap \mathbb{N}$.

Similarly, we have
\[ h_a (\bm{u}) \rightharpoonup | \bm{u} |, \ f \circ h_a (\bm{u}) \rightharpoonup f (| \bm{u} |) \ \text{in } H^1 \big( 0,T;L^q(\Omega) \big), \]
and
\[ g_a (\bm{u}) \rightharpoonup \frac{f(| \bm{u} |)}{| \bm{u} |} \bm{u} \ \text{in } H^1 \big( 0,T;\bm{L}^q(\Omega) \big), \]
for all $q \in [1,\infty)$.
Using Lebesgue's dominated convergence theorem once more, we get
\begin{align*}
    & \partial_t \big(  h_a (\bm{u})  \big) \rightharpoonup \bm{e}_{\mathbb{R}^d} (\bm{u}) \cdot \partial_t \bm{u}, \ \partial_t \big(  f \circ h_a (\bm{u})  \big) \rightharpoonup f'(| \bm{u} |) \bm{e}_{\mathbb{R}^d} (\bm{u}) \cdot \partial_t \bm{u} \ \text{in $L^2 \big( 0,T;L^q(\Omega) \big)$},\\
    & \partial_t \big( g_a (\bm{u}) \big) \rightharpoonup \Big(  f'(| \bm{u} |) - \frac{f(| \bm{u} |)}{| \bm{u} |}  \Big) \big(  \bm{e}_{\mathbb{R}^d}(\bm{u}) \cdot \partial_t \bm{u}  \big) \bm{e}_{\mathbb{R}^d}(\bm{u}) + \frac{f(| \bm{u} |)}{| \bm{u} |} \partial_t \bm{u} \ \text{in $L^2 \big( 0,T;\bm{L}^q(\Omega) \big)$},
\end{align*}
for all $q \in [1,\infty)$.
The proof is complete.
\qed

%%%%%%%%%%%%%%%%%%%%%%%%%%%%%%%%
%%%%%%%%%%%%%%%%%%%%%%%%%%%%%%%%
\section*{Acknowledgments}
\noindent M. Grasselli is a member of Gruppo Nazionale per l'Ana\-li\-si Matematica, la Probabilit\`{a} e le loro Applicazioni (GNAMPA), Istituto Nazionale di Alta Matematica (INdAM). His research is part of the activities of ``Dipartimento di Eccellenza 2023--2027'' of Politecnico di Milano.
H. Wu is a member of the Key Laboratory of Mathematics for Nonlinear Sciences (Fudan University), Ministry of Education of China. 

\section*{Statements and Declarations}
\noindent \textbf{Funding.} 
B.-H. Ouyang was supported by the China Scholarship Council Program (Grant No. 202506100096).
H. Wu was supported by the Natural Science Foundation of Shanghai (Grant No. 25ZR1401023).
\smallskip 

\noindent \textbf{Competing interests.} 
The authors have no competing interests to declare that are relevant to the content of this article.
\smallskip 

\noindent \textbf{Data availability.} 
Data sharing is not applicable to this article as no datasets were generated or analyzed during the current study.

\smallskip
\noindent \textbf{Use of artificial intelligence tools.} 
During the preparation of this manuscript, the authors used ChatGPT (GPT-6 Astra, OpenAI) and Gemini (Gemini-4 Argon, Google) to assist in formulating the technical extensions in Remarks \ref{rem:variable_df_coefficients} and \ref{rem:stopping_time_improvement} and for language polishing. 
All mathematical arguments and results were independently verified by the authors, who take full responsibility for the content of the publication.

%%%%%%%%%%%%%%%%%%%%%%%%%%%%%%%%

\end{document}